\documentclass[reqno,centertags]{amsart}
\usepackage{amsmath,amsthm,amscd,amssymb}
\usepackage{latexsym}
\usepackage{upref}
\newcommand{\bbN}{{\mathbb{N}}}
\newcommand{\bbR}{{\mathbb{R}}}

\newcommand{\bbZ}{{\mathbb{Z}}}
\newcommand{\bbC}{{\mathbb{C}}}
\newcommand{\bbQ}{{\mathbb{Q}}}

\newcommand{\calA}{{\mathcal A}}
\newcommand{\calB}{{\mathcal B}}

\newcommand{\calD}{{\mathcal D}}

\newcommand{\calF}{{\mathcal F}}
\newcommand{\calH}{{\mathcal H}}
\newcommand{\calI}{{\mathcal I}}
\newcommand{\calK}{{\mathcal K}}

\newcommand{\calM}{{\mathcal M}}
\newcommand{\calN}{{\mathcal N}}

\newcommand{\calS}{{\mathcal S}}

\newcommand{\calV}{{\mathcal V}}

\newcommand{\no}{\nonumber}
\newcommand{\lb}{\label}
\newcommand{\f}{\frac}
\newcommand{\ul}{\underline}
\newcommand{\ol}{\overline}
\newcommand{\ti}{\tilde  }
\newcommand{\wti}{\widetilde  }

\newcommand{\spec}{\text{\rm{spec}}}

\newcommand{\ran}{\text{\rm{ran}}}
\newcommand{\dom}{\text{\rm{dom}}}

\newcommand{\supp}{\text{\rm{supp}}}

\newcommand{\bi}{\bibitem}
\newcommand{\hatt}{\widehat}
\newcommand{\beq}{\begin{equation}}
\newcommand{\eeq}{\end{equation}}
\newcommand{\ba}{\begin{align}}
\newcommand{\ea}{\end{align}}

\renewcommand{\Re}{\text{\rm Re}}
\renewcommand{\Im}{\text{\rm Im}}

\DeclareMathOperator*{\slim}{s-lim}
\DeclareMathOperator*{\wlim}{w-lim}
\allowdisplaybreaks
\numberwithin{equation}{section}

\newtheorem{theorem}{Theorem}[section]
\newtheorem{lemma}[theorem]{Lemma}
\newtheorem{corollary}[theorem]{Corollary}
\newtheorem{hypothesis}[theorem]{Hypothesis}
\theoremstyle{definition}
\newtheorem{definition}[theorem]{Definition}
\newtheorem{example}[theorem]{Example}

\theoremstyle{remark}
\newtheorem{remark}[theorem]{Remark}

\begin{document}
\title[Herglotz Operators]{Some Applications of
Operator-Valued Herglotz
Functions}
\dedicatory{Dedicated to Moshe Livsic on the occasion
of his 80th birthday}
\author[Gesztesy, Kalton, Makarov, and
Tsekanovskii]{Fritz Gesztesy, Nigel J.~Kalton,
Konstantin A.~Makarov, \\ and Eduard Tsekanovskii}
\address{Department of Mathematics, University of
Missouri, Columbia, MO
65211, USA}
\email{fritz@math.missouri.edu\newline
\indent{\it URL:}
http://www.math.missouri.edu/people/faculty/fgesztesypt.html}
\address{Department of Mathematics, University of
Missouri, Columbia, MO
65211, USA}
\email{nigel@math.missouri.edu}
\address{Department of Mathematics, University of
Missouri, Columbia, MO
65211, USA}
\email{makarov@azure.math.missouri.edu}
\address{Department of Mathematics, University of Missouri,
Columbia, MO
65211, USA}
\email{tsekanov@math.missouri.edu}
\thanks{Research supported by the US National Science
Foundation under
Grant No.~DMS-9623121.}
\date{\today}

\begin{abstract}
We consider operator-valued Herglotz functions and their
applications
 to self-adjoint perturbations of self-adjoint operators
and self-adjoint
 extensions of densely defined  closed symmetric operators.
Our applications
include model operators for both situations, linear
fractional transformations
for Herglotz operators, results on Friedrichs
and Krein extensions, and realization theorems for classes
of Herglotz operators. Moreover, we study the concrete
case of Schr\"odinger
operators on a  half-line and provide two
illustrations of Livsic's
result \cite{Li60} on quasi-hermitian
extensions in the special case of densely
defined symmetric operators with
deficiency indices $(1,1)$.
\end{abstract}

\maketitle
\section{Introduction}\lb{s1}

The principal  purpose of this paper is to extend some of
our recent
 results on matrix-valued Herglotz functions in
\cite{GT97}
to the infinite-dimensional context.

Given a complex Hilbert space $\calK$, a map
$M: \bbC_+\to \calB(\calK)$
 is called a
$\calK$-valued Herglotz function (or simply a Herglotz
operator) if $M$ is analytic on $\bbC_+$ and
$ \Im (M(z))\ge 0$ for all $z\in \bbC_+.$ (We refer
to the end of this introduction for a
glossary on
the notation used in this paper.)
 $\calB(\calK)$-valued Herglotz functions admit the
celebrated
Nevanlinna-Riesz-Herglotz representation
studied, for instance, by Brodskii \cite{Br71},
Sect.~I.4, Krein and Ovcharenko \cite{KO77},
\cite{KO78}, and
Shmulyan \cite{Sh71} in the infinite-dimensional context,
\begin{equation}
\label{1.1}
M(z)=C+Dz+\int_{\bbR}d\Omega(\lambda)
((\lambda-z)^{-1}-\lambda(1+\lambda^2)^{-1}),
\quad z\in \bbC_+,
\end{equation}
where,
\begin{equation}
\label{1.2}
C=C^*\in \calB(\calK), \quad 0\le D \in \calB(\calK),
\end{equation}
 and $\Omega$ is a $\calB(\calK)$-valued
measure satisfying
\begin{equation}
\label{1.3}
\int_{\bbR}d(\xi, \Omega(\lambda)\xi)_{\calK}
(1+\lambda^2)^{-1}< \infty \text{ for all } \xi\in \calK.
\end{equation}
In this paper we study a subclass of
$\calB(\calK)$-valued Herglotz functions where $D=0$
and the Stieltjes integral in \eqref{1.1} is either
understood in
the norm (cf.~Section~\ref{s2}) or the strong operator
topology (cf.~Section~\ref{s3}) in $\calK$. For detailed
discussions of operator-valued Herglotz functions
and their boundary value behavior, see, for instance,
\cite{Ca76}, \cite{RR97}, Ch.~4,
\cite{SF70}, Ch.~V, \cite{TS77}. Throughout this paper
we will adhere to the usual convention
\begin{equation}
\label{1.4}
M(\bar z)=M(z)^*, \quad z \in\bbC_+
\end{equation}
(see, however, Lemma~\ref{l3.12}).

As discussed in some detail in \cite{GT97}, our notion
of Herglotz functions is not without controversy. In
fact, the names Pick, Nevanlinna,
 Nevanlinna-Pick, and $R$-functions (depending on
whether the open
upper half-plane $\bbC_+$ or the open unit disk $D$ are
involved, as well as
depending on the geographical origin of authors) are
also frequently in use.
 Here we follow a tradition in mathematical physics
which appears to favor the
 terminology of Herglotz functions.

A crucial role in our analysis is played by
linear fractional transformations of the type
\begin{equation}
M(z) \longrightarrow M_A(z)=(A_{2,1}+A_{2,2}M(z))(A_{1,1}
+A_{1,2}M(z))^{-1},
\quad z\in \bbC_+,
\lb{1.5}
\end{equation}
where
\begin{align}
&A=\big(A_{p,q}\big)_{1\leq p,q\leq 2} \in{\calA}(\calK
\oplus \calK),
\no \\
&{\calA}(\calK \oplus \calK)=\{A\in \calB(\calK \oplus
\calK)\,|\, A^*JA
=J \},
\quad J=\begin{pmatrix} 0 & -I_{\calK}\\ I_{\calK}&
0\end{pmatrix}.
\lb{1.6}
\end{align}
$M_A$ is a Herglotz operator in $\calK$ whenever $M$ is
one and we refer
to Krein and Shmulyan \cite{KS74} for a detailed
study in connection with
\eqref{1.5}, \eqref{1.6}.

Section~\ref{s5} provides a detailed study of a
model Hilbert
space, variants of which are used in
Sections~\ref{s2} and
\ref{s3}. This construction appears to be of
independent interest.

In Section~\ref{s2} we consider self-adjoint
perturbations
$H_L$ of a self-adjoint (possibly unbounded)
operator $H_0$
in some separable complex Hilbert
space $\calH$
\begin{equation}
H_L=H_0+KLK^*, \quad \dom (H_L)=\dom (H_0),
\lb{1.7}
\end{equation}
where $L=L^* \in \calB(\calK)$ and $K\in
\calB(\calK, \calH)$,
with $\calK$ another separable complex Hilbert space.
 We introduce a model operator
$\hatt H_L$ in $\hatt\calH=L^2(\bbR, \calK;
d\Omega_L)$
for $H_L$ in $\calH$, define the Herglotz operator
\begin{equation}
M_L(z)=K^*(H_L-z)^{-1}K=
\int_{\bbR}d\Omega_L(\lambda)(\lambda-z)^{-1},
\quad z\in \bbC\backslash \bbR,
\lb{1.8}
\end{equation}
where
\begin{equation}
\Omega_L(\lambda)=K^*E_L(\lambda)K, \lb{1.8a}
\end{equation}
with $\{E_L(\lambda)\}_{\lambda\in\bbR}$ the
family of orthogonal spectral projections of $H_L$,
and study the pair $(H_L,H_0)$ in terms of
$(M_L(z), M_0(z))$  following Donoghue's treatment
\cite{Do65}
 of rank-one perturbations of $H_0$.
Moreover, we prove a realization theorem for
  the class of Herglotz operators exemplified by
\eqref{1.8}.

In Section~\ref{s3} we consider self-adjoint extensions
$H$ of a densely defined closed symmetric operator
$\dot H$ with deficiency indices $(k,k)$,
$k\in \bbN\cup \{\infty\}$
in some separable complex Hilbert space $\calH$.
We review our recent note
\cite{GMT97} on Krein's formula relating self-adjoint
extensions of $\dot H$ and introduce the corresponding
Weyl operators $M_{H,\calN}(z)$
\begin{align}
M_{H,\calN}(z)&=zI_{\calN}+
(1+z^2)P_{\calN}(H-z)^{-1}P_{\calN}\big\vert_{\calN} \\
&=\int_{\bbR}d\Omega_{H,\calN}
(\lambda)((\lambda-z)^{-1}-\lambda(1+\lambda^2)^{-1}),
\quad
 z\in \bbC\backslash\bbR, \lb{1.9}
\end{align}
where $\calN$ is a closed linear subspace of the
deficiency subspace
$\calN_+=\ker(\dot H^*-i), $ $ P_{\calN}$ the
orthogonal projection onto $\calN$, and
\begin{equation}
\Omega_{H,\calN}(\lambda)=(1+\lambda^2)
 (P_{\calN}E_H(\lambda)P_{\calN}\big\vert_{\calN}),
\lb{1.10}
\end{equation}
with $\{E_H(\lambda) \}_{\lambda\in \bbR} $ the
family of orthogonal spectral
projections of $H.$ Following \cite{GMT97} we study
linear fractional
 transformation   of $M_{H,\calN_+}(z)$
involving different self-adjoint extensions $H$ of
$\dot H$. Moreover,
following Donoghue \cite{Do65} in the
special case $\dim_{\bbC}(\calN_+)=1$, we consider
a model
$
(\hatt{ \dot H}, \hatt H)
$ in
$\hatt {\calH}=L^2(\bbR, {\calN}_+;
 d\Omega_{H, {\calN}_+})$
for the pair $(\dot H, H)$ in
$\calH$, and discuss Friedrichs and Krein extensions
of $\dot H$  assuming $\dot H$ to be bounded from
below. We conclude
Section~\ref{s3} with realization theorems for  various
classes of Weyl operators of the type \eqref{1.9}.

Section~\ref{s4} provides concrete applications of the
formalism of Section~\ref{s3} specialized to the case
$\dim_{\bbC}(\calN_+)=1.$ We study Schr\"odinger
operators on a  half-line and provide two
illustrations of Livsic's
result \cite{Li60} on quasi-hermitian
extensions in the special case of densely
defined closed prime symmetric operators with
deficiency indices $(1,1)$.

Finally, we briefly introduce some of the notation
used in this paper.
$\bbC_{\pm}=\{z\in \bbC \,
\vert \,\Im (z) \gtrless 0\}$
denote the open upper/lower half-plane, $\bar  z$
the complex conjugate of $z\in \bbC.$ Complex Hilbert
spaces are denoted by $\calH$ or $\calK$, the scalar
product
in $\calH$ (linear in the second factor) by
$(\cdot \, , \cdot)_{\calH}$, with $I_{\calH}$ the
identity operator in $\calH$.
 Direct sums of linear subspaces are indicated by
$\dot +$,
orthogonal direct sums by $\oplus$ (or $\oplus_{\calH}$,
if  necessary).
The Banach space of bounded linear operators from
$\calK$ into $\calH$
is denoted by
$\calB(\calK,\calH)$
(and simply by $\calB(\calH)$ if $\calK=\calH$). The
domain, range, and kernel (null space) of a linear
operator $T$ are denoted by
$\dom (T)$,
$\ran(T)$ and $\ker(T)$, respectively; the resolvent
set and spectrum of $T$
by $\rho(T)$ and $\spec(T)$. The adjoint of $T$
is denoted by $T^*$,
$\Re(T)=(T+T^*)/2$ and
$\Im (T)=(T-T^*)/(2i)$ (assuming $\dom(T)=\dom(T^*))$
abbreviate the real and imaginary part of $T$,
respectively. The symbol $\chi_B$ denotes the
characteristic function of $B\subset\bbR;$ $\Sigma$
denotes the Borel $\sigma-$algebra on $\bbR$.

\section{Construction of a Model Hilbert Space} \lb{s5}

This section describes in some detail the
construction of
a model Hilbert space, variants of which will be of
crucial
importance in Sections~\ref{s2} and \ref{s3}. Rather
than referring to the theory of direct integrals of
Hilbert spaces (see, e.g., \cite{BW83}, Ch.~4,
\cite{BS87}, Ch.~7) we briefly develop the necessary
machinery from scratch and hint at the construction
of related Banach spaces as well.

Let $\mu$ denote a $\sigma-$finite Borel measure
on $\bbR$,
$\Sigma$ the Borel $\sigma-$algebra on $\bbR$,
and suppose for each $\lambda\in\bbR$ we are given a
separable complex Hilbert space $\calK_\lambda$. Let
$\calS(\{\calK_\lambda\}_{\lambda\in\bbR})$
be the vector space associated with the product space
$\prod_{\lambda\in\bbR}\calK_\lambda$ equipped with
the obvious
linear structure. Elements $f$ of
$\calS(\{\calK_\lambda\}_{\lambda\in\bbR})$ are maps
\begin{equation}
\bbR\ni\lambda\to f=\{f(\lambda)\in\calK_\lambda\}
_{\lambda\in\bbR}\in\prod_{\lambda\in\bbR}\calK_\lambda.
\lb{5.1}
\end{equation}

\begin{definition} \lb{d5.1}
A {\it measurable family of Hilbert spaces} $\calM$
modelled
on $\mu$ and
$\{\calK_\lambda\}_{\lambda\in\bbR}$ is a linear
subspace $\calM\subset\calS(\{\calK_\lambda\}
_{\lambda\in\bbR})$ such that $f\in\calM$ if and
only if the map $\bbR\ni\lambda\to (f(\lambda),
g(\lambda))_{\calK_\lambda}\in\bbC$ is
$\mu-$measurable for all
$g\in\calM$. \\
Moreover, $\calM$ is said to be generated by some
subset $\calF$, $\calF\subset\calM$, if for every
$g\in\calM$ we can
find a sequence of functions $h_n\in\text{lin.span}
\{\chi_Bf\in
\calS(\{\calK_\lambda\}\,|\,
B\in\Sigma,f\in\calF\}$ with
$\lim_{n\to\infty}\|g(\lambda)-h_n(\lambda)
\|_{\calK_{\lambda}}=0$
$\mu-$a.e.
\end{definition}
The definition of $\calM$ was chosen with it's maximality
in mind and we refer to Lemma~\ref{l5.3} and
Theorem~\ref{t5.6} for more details
in this respect. An explicit construction of an example of
$\calM$ will be given in Theorem~\ref{t5.5}.

\begin{remark} \lb{r5.2}
The following properties are proved in a standard manner:\\
(i) If $f\in\calM$, $g\in\calS(\{\calK_\lambda\}
_{\lambda\in\bbR})$ and $g=f$ $\mu-$a.e.~then $g\in\calM$.\\
(ii) If $\{f_n\}_{n\in\bbN}\in\calM$,
$g\in\calS(\{\calK_\lambda\}_{\lambda\in\bbR})$
and $f_n(\lambda)\to g(\lambda)$ as $n\to\infty$
$\mu-$a.e.~(i.e., $\lim_{n\to\infty}\|f_n(\lambda)-
g(\lambda)\|
_{\calK_\lambda}=0$ $\mu-$a.e.) then $g\in\calM$.\\
(iii) If $\phi$ is a scalar-valued $\mu$--measurable
function and $f\in\calM$ then $\phi f\in\calM$.\\
(iv) If $f\in\calM$ then $\bbR\ni\lambda\to
\|f(\lambda)\|_{\calK_\lambda}\in [0,\infty)$ is
$\mu$--measurable.
\end{remark}

Let us remark that we shall identify functions in $\calM$
which coincide
$\mu-$a.e.; thus $\calM$ is more precisely a set of
equivalence classes of functions.

\begin{lemma} \lb{l5.3} Let $\{f_n\}_{n\in\bbN}\subset
\calS(\{\calK_\lambda\}_{\lambda\in\bbR})$
such that\\
($\alpha$) $\bbR\ni\lambda\to (f_m(\lambda),f_n(\lambda)
_{\calK_\lambda}\in\bbC$ is $\mu$--measurable for
all $m,n\in\bbN$.\\
($\beta$) For $\mu-$a.e. $\lambda\in\bbR$,
$\ol{\text{lin.span}\{f_n(\lambda)\}}=\calK_\lambda.$ \\
\noindent Then setting
\begin{equation}
\calM=\{g\in\calS(\{\calK_\lambda\}_{\lambda\in\bbR})
 \, |\, (f_n(\lambda),g(\lambda))_{\calK_\lambda}
\text{ is }
\mu\text{--measurable for all }n\in\bbN\}, \lb{5.3}
\end{equation}
one infers\\
(i) $\calM$ is a measurable family of Hilbert spaces.\\
(ii) $\calM$ is generated by $\{f_n\}_{n\in\bbN}.$\\
(iii) $\calM$ is the unique measurable family of
Hilbert spaces
containing the sequence $\{f_n\}_{n\in\bbN}.$\\
(iv) If $\{g_n\}$ is any sequence satisfying
$(\beta)$ then
$\calM$ is generated by $\{g_n\}.$
\end{lemma}
\noindent {\it Sketch of proof. }(i) Without loss of
generality, we
may assume $\{f_n\}_{n\in\bbN}$ contains all rational
linear
combinations, that is, all elements of the type
$\sum_{n=1}^N
\alpha_n f_n,$ with $\alpha_n\in\bbQ,$ $n=1,\dots,N,$
$N\in\bbN$.
For $f\in\calS(\{\calK_\lambda\}_{\lambda\in\bbR})$,
\begin{equation}
\|f(\lambda)\|_{\calK_\lambda} =\sup_{n\in\bbN}
|(f(\lambda),
\chi_{B_n}(\lambda)f_n(\lambda))_{\calK_\lambda}|,
\lb{5.4}
\end{equation}
where $B_n=\{\lambda\in\bbR\, | \, \|f_n(\lambda)
\|_{\calK_\lambda}\leq1\}$. Hence, if $f\in\calM$ then
 the map $\bbR\ni\lambda\to \|f(\lambda)
\|_{\calK_\lambda}
\in[0,\infty)$ is $\mu$--measurable. It then follows
easily that
$\calM$ is a measurable family of Hilbert spaces.\\
(ii) If $g\in\calM$ then
\begin{equation}
\inf_{n\in\bbN} \|g(\lambda)-f_n(\lambda)
\|_{\calK_\lambda}=0
\quad \mu-\text{a.e.} \lb{5.5}
\end{equation}
It follows that if $\varepsilon (\lambda)$ is any
measurable function with $\varepsilon >0$ on $\bbR$,
then one can find a measurable partition
$\{B_n\}_{n\in\bbN}$
of $\bbR$ so that
\begin{equation}
\|g(\lambda)-\sum_{n\in
\bbN}\chi_{B_n}(\lambda)f_n(\lambda)
\|_{\calK_\lambda}\le
\epsilon(\lambda). \lb{5.4a}
\end{equation}
Indeed, for each $\lambda\in\bbR$ let $N(\lambda)$ be
the first
$n$ such that
\begin{equation}
\|g(\lambda)-f_{N(\lambda)}(\lambda)\|_{\calK_\lambda}
<\varepsilon(\lambda).
\lb{5.4b}
\end{equation}
Then $\bbR\ni\lambda\to N(\lambda)\in\bbN$ is
$\mu-$measurable and $B_n=\{\lambda\in\bbR \,|\,
N(\lambda)=n\}$ is the desired partition.
This implies (ii).\\
(iii) If ${\calM}'\subset\calS(\{\calK_\lambda\}
_{\lambda\in\bbR})$ is a measurable family of
Hilbert spaces
containing each $\{f_n\}_{n\in\bbN}$, then
${\calM}'\subseteq\calM$.
However, $\calM\subseteq{\calM}'$ by (ii) then
completes the argument. \\
(iv)This follows immediately from (iii), since
we can define $\calM'$
in a similar way, that is,
\begin{equation}
\calM'=\{h\in\calS(\{\calK_\lambda\}_{\lambda\in\bbR})
 \, |\, (g_n(\lambda),h(\lambda))_{\calK_\lambda}
\text{ is }
\mu\text{--measurable for all }n\in\bbN\}, \lb{5.5a}
\end{equation}
and then $\calM=\calM'$ is clear from (iii).
$\square$\\

Next, let $w$ be a $\mu$--measurable function,
$w>0$ $\mu-$a.e., and consider the space
\begin{equation}
\dot L^2(\calM;wd\mu)=\{f\in\calM \, | \,
\int_\bbR w(\lambda)
d\mu(\lambda) \|f(\lambda)\|^2_{\calK_\lambda}
<\infty\}
\lb{5.6}
\end{equation}
with its obvious linear structure. On
$\dot L^2(\calM;wd\mu)$ one defines a semi-inner
product
$(\cdot,\cdot)_{\dot L^2(\calM;wd\mu)}$ (and hence a
semi-norm $\|\cdot\|_{\dot L^2(\calM;wd\mu)}$) by
\begin{equation}
(f,g)_{\dot L^2(\calM;wd\mu)}=\int_\bbR
w(\lambda)d\mu(\lambda)
(f(\lambda),g(\lambda))_{\calK_\lambda}, \quad
f,g\in \dot L^2(\calM;wd\mu). \lb{5.7}
\end{equation}
That \eqref{5.7} defines a semi-inner product
immediately follows
from the corresponding properties of
$(\cdot,\cdot)_{\calK_\lambda}$ and the linearity of
the integral. Hence $\dot L^2(\calM;wd\mu)$ represents
a pre-Hilbert space and one can complete it in a standard
manner as follows. One defines the equivalence relation
$\sim$, for elements $f,g\in \dot L^2(\calM;wd\mu)$ by
\begin{equation}
f\sim g \text{ if and only if } f=g \quad \mu-\text{a.e.}
\lb{5.7a}
\end{equation}
and hence introduces the set of equivalence classes of
$\dot L^2(\calM;wd\mu)$ denoted by
\begin{equation}
L^2(\calM;wd\mu)=\dot L^2(\calM;wd\mu)/\sim .
\lb{5.7b}
\end{equation}
In particular, introducing the subspace of null functions
\begin{align}
\calN(\calM;wd\mu)&=\{f\in \dot L^2(\calM;wd\mu)
\,|\, \|f(\lambda)\|_{\calK_\lambda}=0 \text{ for }
\mu-\text{a.e. }\lambda\in\bbR\} \no \\
&=\{f\in \dot L^2(\calM;wd\mu) \,|\,
\|f\|_{\dot L^2(\calM;wd\mu)}=0 \}, \lb{5.7c}
\end{align}
$L^2(\calM;wd\mu)$ is precisely the quotient space
$\dot L^2(\calM;wd\mu)/\calN(\calM;wd\mu).$
Denoting the equivalence class of
$f\in \dot L^2(\calM;wd\mu)$ temporarily by $[f]$, the
semi-inner product on $L^2(\calM;wd\mu)$
\begin{equation}
([f],[g])_{L^2(\calM;wd\mu)}=\int_\bbR
w(\lambda)d\mu(\lambda)
(f(\lambda),g(\lambda))_{\calK_\lambda} \lb{5.7d}
\end{equation}
is well-defined (i.e., independent of the chosen
representatives of the equivalence classes) and actually
an inner product. Thus
$L^2(\calM;wd\mu)$ is a normed space and by the usual
abuse of notation we denote its elements in the
following again by
$f,g,$ etc.~The fundamental fact that
$L^2(\calM;wd\mu)$ is also complete is discussed next.

\begin{theorem} \lb{t5.4}
$L^2(\calM;wd\mu)$ is complete and hence a Hilbert space.
\end{theorem}
\begin{proof}It suffices to prove the following fact:
For each $\{f_n\}_{n\in\bbN}\in L^2(\calM;wd\mu)$
with
$\sum_{n\in\bbN} \|f_n\|_{L^2(\calM;wd\mu)}<\infty$,
there is an $f\in L^2(\calM;wd\mu)$ such that
$\sum_{n\in\bbN} f_n=f$. Given such a sequence
$\{f_n\}_{n\in\bbN}$ with $\sum_{n\in\bbN}
\|f_n\|_{L^2(\calM;wd\mu)}=A$ define
\begin{equation}
G(\lambda)=\bigg(\sum_{n\in\bbN}\|f_n(\lambda)\|
_{\calK_\lambda}\bigg)^2. \lb{5.7e}
\end{equation}
Then $G$ is $\mu-$measurable.~From $\sum_{n=1}^N
\|f_n\|_{L^2(\calM;wd\mu)}\leq A$ one computes using
Minkowski's inequality,
\begin{align}
\bigg(\int_\bbR w(\lambda)d\mu(\lambda)
\bigg(\sum_{n=1}^N
\|f_n(\lambda)\|_{\calK_\lambda}\bigg)^2\bigg)^{1/2}
&\leq \sum_{n=1}^N \bigg(\int_\bbR w(\lambda)
d\mu(\lambda)
\|f_n(\lambda)\|_{\calK_\lambda}^2\bigg)^{1/2} \no \\
&=\sum_{n=1}^N \|f_n\|_{L^2(\calM;wd\mu)}\leq A,
\lb{5.7f}
\end{align}
that is,
\begin{equation}
\int_\bbR w(\lambda)d\mu(\lambda) \bigg(\sum_{n=1}^N
\|f_n(\lambda)\|_{\calK_\lambda}\bigg)^2\leq A^2.
\lb{5.7fa}
\end{equation}
Applying the Monotone Convergence Theorem one then
concludes
\begin{equation}
\int_\bbR w(\lambda)d\mu(\lambda) G(\lambda)\le A^2.
\lb{5.7fb}
\end{equation}
Thus $G$ is integrable and hence $\mu-$a.e. finite.
Consequently, we may define
\begin{align}
&f(\lambda)=\begin{cases}\sum_{n\in\bbN} f_n(\lambda),
&\text{if }
\sum_{n\in\bbN} \|f_n(\lambda)\|_{\calK_{\lambda}}
<\infty, \lb{5.7fc} \\
0, & \text{otherwise.} \end{cases}
\end{align}
Then $\|f(\lambda)\|
_{\calK_\lambda}^2\leq G(\lambda)$ for $\mu-$a.e.
$\lambda\in\bbR$ and
\begin{equation}
\sum_{n\in\bbN} f_n(\lambda)=f(\lambda)\quad
\mu-\text{a.e.} \lb{5.7g}
\end{equation}
In particular, $f\in L^2(\calM;wd\mu)$. Finally,
since
\begin{equation}
\bigg\|\sum_{n=1}^N f_n(\lambda)-f(\lambda)\bigg\|
_{\calK_\lambda}\to 0 \text{ as } N\to\infty \quad
\mu-\text{a.e.} \lb{5.7h}
\end{equation}
and
\begin{equation}
\bigg\|f(\lambda)-\sum_{n=1}^N f_n(\lambda)\bigg\|^2
_{\calK_\lambda}
=\bigg\|\sum_{n=N+1}^\infty f_n(\lambda)\bigg\|^2
_{\calK_\lambda}
\leq G(\lambda) \quad \mu-\text{a.e.,} \lb{5.7i}
\end{equation}
the Lebesgue Dominated Convergence theorem yields
\begin{equation}
\lim_{N\to\infty}\bigg\|f-\sum_{n=1}^N f_n\bigg\|
_{L^2(\calM;wd\mu)}=0. \lb{5.7j}
\end{equation}
\end{proof}

Clearly, the analogous construction defines the
Banach spaces $L^p(\calM;wd\mu),$ $p\geq 1$. The
case $p=2$ corresponds precisely to the direct integral
of the Hilbert spaces $\calK_\lambda$ with respect to
the measure $wd\mu$ (see, e.g., \cite{BW83}, Ch.~4,
\cite{BS87}, Ch.~7).

Next, suppose $\calK$ is a separable complex Hilbert
space
and $\Omega:\Sigma\to\calB(\calK)$ is a positive measure
(i.e., countably additive with respect to the strong
operator topology in $\calK$). Assume
\begin{equation}
\Omega(\bbR)=T\geq 0, \quad T\in\calB(\calK). \lb{5.8}
\end{equation}
Moreover, let
$\mu$ be a control measure for $\Omega$, that is,
\begin{equation}
\mu(B)=0 \text{ if and only if } \Omega(B)=0
\text{ for all } B\in\Sigma. \lb{5.8a}
\end{equation}
(E.g., $\mu(B)=\sum_{n\in\calI}2^{-n}(e_n,
\Omega(B)e_n)_\calK$,
with $\{e_n\}_{n\in\calI}$ a complete orthonormal
system in $\calK,$ $\calI\subseteq\bbN$ an appropriate
index set.)

\begin{theorem} \lb{t5.5}
There are separable complex Hilbert spaces
$\calK_\lambda$,
$\lambda\in\bbR$, a measurable family of Hilbert spaces
${\calM}_\Omega(\mu)$
modelled on $\mu$ and $\{\calK_\lambda\}
_{\lambda\in\bbR}$,
and a bounded linear map $\ul \Lambda\in\calB(\calK,
L^2({\calM}_\Omega(\mu);d\mu))$ so that \\
(i) For all $B\in\Sigma$, $\xi,\eta\in\calK$,
\begin{equation}
(\xi,\Omega(B) \eta)_\calK=\int_B d\mu(\lambda)
((\ul \Lambda \xi)(\lambda),(\ul \Lambda \eta)(\lambda))
_{\calK_\lambda}. \lb{5.9}
\end{equation}
(ii) $\ul \Lambda(\{e_n\}_{n\in\calI})$ generates
${\calM}_\Omega(\mu)$, where
$\{e_n\}_{n\in\calI}$ denotes any
sequence of linearly independent elements in $\calK$
with the property $\ol{\text{lin.span}
\{e_n\}_{n\in\calI}}=\calK,$ $\calI\subseteq\bbN$.
In particular, $\ul \Lambda(\calK)$ generates
${\calM}_\Omega(\mu)$.\\
(iii) For all $\xi\in\calK$,
\begin{equation}
\ul \Lambda(\Omega(B) \xi)=\chi_B\ul \Lambda \xi \quad
\mu-\text{a.e.} \lb{5.10}
\end{equation}
\end{theorem}
\begin{proof}
Denote $\calV=\text{lin.span}
\{e_n\}_{n\in\calI}$. By the Radon-Nikodym theorem, there
exist $\mu$--measurable $\phi_{m,n}$ such that
\begin{equation}
\int_B d\mu(\lambda)\phi_{m,n}(\lambda)=
(e_m,\Omega(B)e_n)_\calK. \lb{5.11}
\end{equation}
Next, suppose $v=\sum_{n=1}^N \alpha_n e_n \in\calV$,
$\alpha_n\in\bbC$, $n=1,\dots,N,$ $N\in\calI$. Then
\begin{equation}
(v,\Omega(B)v)_\calK=\int_B d\mu(\lambda)\sum_{m,n=1}^N
\phi_{m,n}(\lambda) \ol{\alpha_m} \alpha_n. \lb{5.12}
\end{equation}
By considering only rational linear combinations
we can deduce that for $\mu-$a.e. $\lambda\in\bbR$,
\begin{equation}
\sum_{m,n} \phi_{m,n}(\lambda)\ol{\alpha_m}
\alpha_n\geq 0
\text{ for all finite sequences }
\{\alpha_n\}\subset\bbC. \lb{5.13}
\end{equation}
Hence we can define a semi-inner product
$(\cdot,\cdot)
_{\lambda}$ on $\calV$ such that
\begin{equation}
(v,w)_{\lambda}=\sum_{m,n}\phi_{m,n}(\lambda)
\ol{\alpha_m} \beta_n \quad \mu-\text{a.e}
\lb{5.14}
\end{equation}
if $v=\sum_n \alpha_n e_n$, $w=\sum_n \beta_n e_n$. \\
Next, let $\calK_\lambda$ be the completion
of $\calV$
with respect $\|\cdot\|_\lambda$ (or, more
precisely the
completion of
$\calV/\calN_{\lambda}$ where
$\calN_{\lambda}=\{\xi\in\calV\,|\,(\xi,\xi)
_{\lambda}=0\}$)
and consider
$\calS(\{\calK_\lambda\}_{\lambda\in\bbR})$. Each
$v\in\calV$ defines an element
$\ul v=\{\ul v(\lambda)\}_{\lambda\in\bbR}
\in\calS(\{\calK_\lambda\}_{\lambda\in\bbR})$
by
\begin{equation}
\ul v (\lambda)=v \text{ for all } \lambda\in\bbR.
\lb{5.15}
\end{equation} Again we identify an element
$v\in \calV$ with an element in
$\calV/{\calN}_{\lambda}\subseteq
\calK_{\lambda}.$ Applying Lemma~\ref{l5.3}, the
collection
$\{\ul e_n\}_{n\in\calI}$ then
generates a measurable family of Hilbert spaces
${\calM}_\Omega(\mu).$ If $v\in\calV$ then
\begin{equation}
\|\ul v\|^2_{L^2({\calM}_\Omega(\mu);d\mu)}=\int_\bbR
d\mu(\lambda) (\ul v(\lambda),\ul v(\lambda))_\lambda
=(v,Tv)_\calK=\|T^{1/2}v\|^2_\calK. \lb{5.16}
\end{equation}
Hence we can define
\begin{equation}
\dot {\ul \Lambda} :\calV\to L^2({\calM}_\Omega(\mu);d\mu),
\quad v\to\dot{\ul \Lambda} v=\ul v
=\{\ul v (\lambda)=v\}_{\lambda\in\bbR}\lb{5.17}
\end{equation}
and denote by $\ul \Lambda\in\calB(\calK,
L^2({\calM}_\Omega(\mu);d\mu)),$ $\|\ul \Lambda\|
_{\calB(\calK,
L^2({\calM}_\Omega(\mu);d\mu))}=
\|T^{1/2}\|_{\calB(\calK)},$ the closure of
$\dot {\ul \Lambda}$.~Then properties
(i)--(iii) hold.
\end{proof}

We now  show that this construction is essentially
unique.

\begin{theorem} \lb{t5.6}
Suppose $\calK'_{\lambda},\ \lambda\in\bbR$ is a
family of separable
complex Hilbert spaces, $\calM'$ is a measurable
family of Hilbert spaces
modelled on $\mu$ and $\{\calK'_{\lambda}\},$ and
$\ul \Lambda'\in\calB(\calK,L^2(\calM';d\mu))$ is
a map satisfying (i),(ii), and (iii) of
the preceding theorem.
Then for $\mu$-a.e. $\lambda\in\bbR$ there is a
unitary operator
$U_{\lambda}:\calK_{\lambda}\to\calK'_{\lambda}$
such that $f=\{f(\lambda)\}_{\lambda\in\bbR}
\in \calM_{\Omega}(\mu)$ if and only if
$U_{\lambda}f(\lambda)\in\calM'$ and for all
$\xi\in\calK,$
\begin{equation}
(\ul \Lambda' \xi)(\lambda)= U_{\lambda}
(\ul \Lambda \xi)(\lambda) \quad \mu-\text{a.e.}
\lb{5.17a}
\end{equation}
\end{theorem}

\begin{proof}We use the notation of the preceding
theorem.~We select
representatives $f'_n\in\calM'$ of $\ul \Lambda'e_n.$
It follows from
condition (i) that for $\mu-$a.e.~$\lambda\in\bbR$
and every $m,n\in\calI$ we
have
\begin{equation}
(f'_m(\lambda),f'_n(\lambda))_{\calK'_{\lambda}}=
(e_m,e_n)_{\lambda}=(\ul e_m(\lambda),
\ul e_n(\lambda))_{\calK_{\lambda}}.
\lb{5.17b}
\end{equation}
Hence we can induce an isometry $U_{\lambda}:
\calK_{\lambda}\to
\calK'_{\lambda}$ such that $U_{\lambda}
\ul e_n(\lambda)=f'_n(\lambda).$

It is easy to see that if $v\in\calV$ we must have
$U_{\lambda}\ul v(\lambda)=
(\ul \Lambda' v)(\lambda)$ $\mu-$a.e.  From
the $L_2-$continuity
of both $\ul \Lambda$ and $\ul \Lambda'$ it follows
that for
every $\xi\in\calK
$ we have
\begin{equation}
(\ul \Lambda' \xi)(\lambda)= U_{\lambda}
(\ul \Lambda \xi)(\lambda) \quad \mu-\text{a.e.}
\lb{5.17c}
\end{equation}

We next observe that if $\ul \Lambda'(\calK)$ generates
$\calM'$ then by a
density argument it must also be true that
$\{f'_n\}_{n\in\calI}$
generates $\calM'.$
It is then immediate that the linear span of
$\{f'_n(\lambda)\}_{n\in\calI}$ must be dense for
$\mu-$a.e.~$\lambda\in\bbR.$  Thus $U_{\lambda}$ is
actually surjective $\mu-$a.e. and so is unitary.

Finally, if $\xi\in\calK$ and $B\in\Sigma$ then
$U_{\lambda}(\chi_B(\lambda)(\ul \Lambda \xi)(\lambda))=
\chi_{B}(\lambda)(\ul \Lambda' \xi)(\lambda)$ $\mu-$a.e.
Thus it follows
by approximation that if $f\in\calM_{\Omega}(\mu)$ then
$U_{\lambda}f(\lambda)\in\calM'$.  Conversely, a similar
argument shows
that if $f\in\calM'$ then $U_{\lambda}^{-1}f(\lambda)\in
\calM_{\Omega}(\mu).$
\end{proof}

Without going into further details, we note that
$\calM_\Omega(\mu)$ depends of course on $\mu$.
However, a change in $\mu$ merely effects a change in
density and so $\calM_\Omega(\mu)$ can essentially be
viewed as $\mu-$independent.

Next, using the notation employed in the proof of
Theorem~\ref{5.5} we recall
\begin{equation}
\calV=\text{lin.span}\{e_n\in\calK \,|\,n\in\calI\}
\lb{5.17d}
\end{equation}
and define
\begin{equation}
{\ul \calV}_\Omega=\text{lin.span}\{\chi_B\ul e_n\in
L^2(\calM_\Omega(\mu);d\mu) \, | \, B\in\Sigma,
\,n\in\calI\}. \lb{5.18}
\end{equation}
The fact that $\{\ul e_n\}_{n\in\calI}$ generates
$\calM_\Omega(\mu)$ implies that ${\ul \calV}_\Omega$
is dense in $L^2({\calM}_\Omega(\mu);d\mu)$,
that is,
\begin{equation}
\ol{{\ul \calV}_{\Omega}}=L^2({\calM}
_\Omega(\mu);d\mu). \lb{5.19}
\end{equation}
The following result will be used in
Section~\ref{s2}.

\begin{lemma} \lb{l5.7}
Suppose $\calK$, $\calH$ are separable complex
Hilbert spaces, $K\in\calB(\calK,\calH)$,
$\{E(B)\}_{B\in\Sigma}$ is a family of orthogonal
projections in $\calH,$ and assume
\begin{equation}
\ol{\text{lin.span}\{E(B)Ke_n\in\calH
\,|\, B\in\Sigma,\, n\in\calI\}}=\calH, \lb{5.20}
\end{equation}
with $\{e_n\}_{n\in\calI},$ $\calI\subset\bbN$ a
complete orthonormal system in $\calK.$ Define
\begin{equation}
\Omega:\Sigma\to\calB(\calK), \quad \Omega(B)=K^*E(B)K,
\lb{5.20a}
\end{equation}
and introduce
\begin{align}
&\dot U:{\ul {\calV}}_\Omega\to\calH, \no \\
&{\ul \calV}_\Omega\ni \sum_{m=1}^M \sum_{n=1}^N
\alpha_{m,n}\chi_{B_m}\ul e_n
\to \dot U\bigg(\sum_{m=1}^M \sum_{n=1}^N \alpha_{m,n}
\chi_{B_m}\ul e_n\bigg) \lb{5.21} \\
&\hspace*{4.35cm}=\sum_{m=1}^M \sum_{n=1}^N \alpha_{m,n}
E(B_m)Ke_n\in\calH, \no \\
&\hspace*{1.75cm} \alpha_{m,n}\in\bbC, \, m=1,\dots,M,
\, n=1,\dots,N, \, M,N\in\calI. \no
\end{align}
Then $\dot U$ extends to a unitary operator
$U:L^2({\calM}_\Omega(\mu);d\mu)\to\calH$.
\end{lemma}

\begin{proof} One computes
\begin{align}
&\bigg\|\dot U\bigg(\sum_{m=1}^M \sum_{n=1}^N \alpha_{m,n}
\chi_{B_m}\ul e_n\bigg)\bigg\|^2_\calH \no \\
&= \sum_{m_1,m_2=1}^M \sum_{n_1,n_2=1}^N
\ol{\alpha_{m_1,n_1}}\alpha_{m_2,n_2}
(e_{n_1},K^* E(B_{m_1}\cap B_{m_2})Ke_{n_2})_\calK \no \\
&=\sum_{m_1,m_2=1}^M  \sum_{n_1,n_2=1}^N
\ol{\alpha_{m_1,n_1}}\alpha_{m_2,n_2}
(e_{n_1},\Omega(B_{m_1}\cap B_{m_2})e_{n_2})_\calK \no \\
&=\sum_{m_1,m_2=1}^M \sum_{n_1,n_2}^N
\ol{\alpha_{m_1,n_1}}\alpha_{m_2,n_2}
\int_{B_{m_1}\cap B_{m_2}} d\mu(\lambda)
(\ul e_{n_1}(\lambda),\ul e_{n_2}
(\lambda))_{\calK_\lambda} \no \\
&=\bigg\|\sum_{m=1}^M \sum_{n=1}^N \alpha_{m,n}
\chi_{B_m}\ul e_n\bigg\|^2_{L^2({\calM}
_\Omega(\mu);d\mu)}. \lb{5.23}
\end{align}
By \eqref{5.19}, $\dot U$ is densely defined and thus
extends to an isometry
$U$ of $L^2({\calM}_\Omega(\mu);d\mu)$ into $\calH$.
In particular, $\ran(U)$ is closed in $\calH$. Thus,
\begin{equation}
\ran(U)\supseteq\ol{\text{lin.span}\{E(B)Ke_n\in\calH
\,|\, B\in\Sigma,\, n\in\calI\}}=\calH \lb{5.24}
\end{equation}
by hypothesis \eqref{5.20a} and hence
$U:L^2({\calM}_\Omega(\mu);d\mu)\to\calH$ is a
unitary operator.
\end{proof}

In view of our comment following Theorem~\ref{t5.6},
concerning the mild dependence on the control measure
$\mu$ of
$\calM_\Omega(\mu)$, we will put more emphasis on the
operator-valued measure $\Omega$ and hence use the
notation $L^2(\bbR,\calK;wd\Omega)$ instead of the more
precise $L^2({\calM}_\Omega(\mu);wd\mu)$ in
Section~\ref{s2}.
\vspace*{2mm}

Finally we adapt Lemma~\ref{l5.7} to the content
of Section~\ref{s3}.

Suppose $\calN$ is a separable complex Hilbert space
and $\wti \Omega:\Sigma\to\calB(\calN)$ a positive
measure. Assume
\begin{equation}
\wti \Omega (\bbR)=\wti T\geq 0, \quad \wti T\in
\calB(\calN) \lb{5.25}
\end{equation}
and let $\ti \mu$ be a control measure for $\wti \Omega$.
Moreover, let $\{u_n\}_{n\in\calI},$ $\calI\subseteq\bbN$
be a sequence of linearly independent elements in
$\calN$ with the property $\ol{\text{lin.span}\{u_n\}_
{n\in\calI}}=\calN.$ As discussed in Theorem~\ref{t5.5},
this yields a measurable family of Hilbert spaces
$\calM_{\ti \Omega}(\ti \mu)$ modelled on $\ti \mu$ and
$\{\calN_\lambda\}_{\lambda\in\bbR}$ and a bounded map
$\ul \Lambda\in\calB(\calN, L^2(\calM_{\ti \Omega}
(\ti \mu);d\ti \mu)),$ $\|\ul
\Lambda\|_{\calB(\calN, L^2(\calM_{\ti \Omega}
(\ti \mu);d\ti \mu))}=
\|{\wti T}^{1/2}\|_{\calB(\calN)}$, such that
$\ul \Lambda (\{u_n\}
_{n\in\calI})$ generates $\calM_{\ti \Omega}(\ti \mu)$
and
\begin{equation}
\ul \Lambda:\calV\to L^2(\calM_{\ti \Omega}
(\ti \mu);d\ti \mu)), \quad v\to\ul \Lambda v =
\ul v =\{\ul v (\lambda)=v\}_{\lambda\in\bbR},
\lb{5.26}
\end{equation}
where
\begin{equation}
\calV = \text{lin.span}\{u_n\}_{n\in\calI}.
\lb{5.27}
\end{equation}
Each $v\in\calV$ defines an element
\begin{equation}
\ul {\ul v}=\{\ul {\ul v}(\lambda)=(\lambda-i)^{-1}v\}
_{\lambda\in\bbR}\in\calS (\{\calN_\lambda\}
_{\lambda\in\bbR}) \lb{5.28}
\end{equation}
and introducing the weight function
\begin{equation}
w_1(\lambda)=1+\lambda^2, \quad \lambda\in\bbR
\lb{5.29}
\end{equation}
and Hilbert space $L^2(\calM_{\ti \Omega}
(\ti \mu);w_1d\ti \mu)$ one computes
\begin{equation}
\|\ul {\ul v}\|^2_{L^2(\calM_{\ti \Omega}
(\ti \mu);d\ti \mu)}=\int_\bbR d\ti \mu (\lambda)
\|\ul v(\lambda)\|^2_{\calN_\lambda}=
(v,\wti T v)_\calN=\|{\wti T}^{1/2} v\|^2_\calN.
\lb{5.30}
\end{equation}
Thus, the linear map
\begin{equation}
\ul {\ul {\dot \Lambda}}:\calV\to L^2(\calM_
{\ti \Omega}(\ti \mu);w_1d\ti \mu), \quad
v\to\ul {\ul {\dot \Lambda}}v=\ul {\ul v}=
\{\ul {\ul v}(\lambda)=(\lambda-i)^{-1} v\}
_{\lambda\in\bbR} \lb{5.31}
\end{equation}
extends to $\ul {\ul \Lambda}\in \calB(\calN,
L^2(\calM_{\ti \Omega} (\ti \mu);w_1d\ti \mu)),$
$\|\ul {\ul \Lambda}\|_{\calB(\calN,
L^2(\calM_{\ti \Omega} (\ti \mu);w_1d\ti \mu))}
=\|{\wti T}^{1/2}\|_{\calB(\calN)}$. Introducing
\begin{equation}
{\ul {\ul \calV}}_{\ti \Omega}=\text{lin.span}
\{\chi_B \ul {\ul v}\in L^2(\calM_{\ti \Omega}
(\ti \mu);w_1d\ti \mu) \,|\,B\in\Sigma,\, n\in\calI\}
\lb{5.32}
\end{equation}
one infers that ${\ul {\ul \calV}}_{\ti \Omega}$ is
dense in $L^2(\calM_{\ti \Omega}(\ti \mu);w_1d\ti \mu),$
that is,
\begin{equation}
\ol{{\ul {\ul \calV}}_{\ti \Omega}}=
L^2(\calM_{\ti \Omega}(\ti \mu);w_1d\ti \mu).
\lb{5.33}
\end{equation}

Given these preliminaries we can state the following
result.

\begin{lemma} \lb{l5.8}
Suppose $\calH$ is a separable complex Hilbert space,
$\calN$ a closed linear subspace of $\calH,$
$P_\calN$ the orthogonal projection in
$\calH$ onto $\calN,$ $\{E(B)\},$ ${B\in\Sigma}$ a
family of orthogonal projections in $\calH,$ and
assume
\begin{equation}
\ol{\text{lin.span}\{E(B)u_n\in\calH\,|\,B\in\Sigma,
\,n\in\calI \}}=\calH, \lb{5.34}
\end{equation}
with $\{u_n\}_{n\in\calI},$ $\calI\subseteq\bbN$ a
complete orthonormal system in $\calN$. Define
\begin{equation}
\wti \Omega:\Sigma\to\calB(\calN), \quad
\wti \Omega (B)=P_\calN E(B) P_\calN\big|_\calN,
\lb{5.35}
\end{equation}
and introduce
\begin{align}
&\dot {\wti U}:{\ul {\ul \calV}}_\Omega\to\calH, \no \\
&{\ul {\ul \calV}}_\Omega\ni \sum_{m=1}^M \sum_{n=1}^N
\alpha_{m,n}\chi_{B_m}\ul {\ul u}_n
\to \dot {\wti U}\bigg(\sum_{m=1}^M \sum_{n=1}^N
\alpha_{m,n}
\chi_{B_m}\ul {\ul u}_n\bigg) \lb{5.36} \\
&\hspace*{4.4cm}=\sum_{m=1}^M \sum_{n=1}^N \alpha_{m,n}
E(B_m)u_n\in\calH, \no \\
&\hspace*{1.55cm} \alpha_{m,n}\in\bbC, \, m=1,\dots,M,
\, n=1,\dots,N, \, M,N\in\calI. \no
\end{align}
Then $\dot {\wti U}$ extends to a unitary operator
$\wti U:L^2({\calM}_{\ti \Omega}(\ti \mu);
w_1d\mu)\to\calH$.
\end{lemma}
\begin{proof}
One computes
\begin{align}
&\bigg\|\dot {\wti U}\bigg(\sum_{m=1}^M \sum_{n=1}^N
\alpha_{m,n}
\chi_{B_m}\ul {\ul u}_n\bigg)\bigg\|^2_\calH \no \\
&= \sum_{m_1,m_2=1}^M \sum_{n_1,n_2=1}^N
\ol{\alpha_{m_1,n_1}}\alpha_{m_2,n_2}
(u_{n_1},E(B_{m_1}\cap B_{m_2})u_{n_2})_\calN
\no \\
&=\sum_{m_1,m_2=1}^M  \sum_{n_1,n_2=1}^N
\ol{\alpha_{m_1,n_1}}\alpha_{m_2,n_2}
(u_{n_1},\wti \Omega(B_{m_1}\cap B_{m_2})u_{n_2})_\calN
\no \\
&=\sum_{m_1,m_2=1}^M \sum_{n_1,n_2}^N
\ol{\alpha_{m_1,n_1}}\alpha_{m_2,n_2}
\int_{B_{m_1}\cap B_{m_2}} d\ti \mu(\lambda)
(\ul u_{n_1}(\lambda),\ul u_{n_2}
(\lambda))_{\calN_\lambda} \no \\
&=\bigg\|\sum_{m=1}^M \sum_{n=1}^N \alpha_{m,n}
\chi_{B_m}\ul {\ul u}_n\bigg\|^2_{L^2({\calM}
_{\ti \Omega}(\ti\mu);w_1d\ti\mu)}. \lb{5.37}
\end{align}
By \eqref{5.33}, $\dot {\wti U}$ is densely defined
and extends to an isometry
$\wti U$ of $L^2({\calM}_{\ti \Omega}(\ti\mu);
w_1d\ti\mu)$ into $\calH$.
In particular, $\ran(\wti U)$ is closed in $\calH$.
Thus,
\begin{equation}
\ran(\wti U)\supseteq\ol{\text{lin.span}
\{E(B)u_n\in\calH
\,|\, B\in\Sigma,\, n\in\calI\}}=\calH \lb{5.38}
\end{equation}
by hypothesis \eqref{5.34} and hence
$U:L^2({\calM}_{\ti \Omega}(\ti\mu);w_1d\ti\mu)
\to\calH$ is a unitary operator.
\end{proof}

Analogous to our comments following Lemma~\ref{l5.7},
in Section~\ref{s3} we will emphasize the
role of $\wti \Omega$ and hence use the somewhat
imprecise notation $L^2(\bbR,\calN;wd\wti\Omega)$,
with various weight functions $w$, as opposed to the
precise notation $L^2(\calM_{\ti \Omega}(\ti\mu);
wd\ti\mu)$.

\section{On Self-Adjoint Perturbations of Self-Adjoint
Operators} \lb{s2}

In this section
we will focus on the following perturbation problem.
Assuming
\begin{hypothesis}\label{h2.1}
  Let $\calH$ and $\calK$ be
separable complex Hilbert spaces,  $H_0$ a self-adjoint
 (possibly unbounded) operator in $\calH$, $L$ a
bounded
self-adjoint operator in $\calK$, and $K:\calK\to
\calH$ a bounded operator,
\end{hypothesis}
\noindent we define the self-adjoint operator
$H_L$ in $\calH$,
\begin{equation}
H_L=H_0+KLK^*, \quad \dom(H_L)=\dom(H_0).
\lb{2.1}
\end{equation}
Given the perturbation $H_L$ of $H_0$,
 we introduce the associated operator-valued Herglotz
function in $\calK$,
\begin{align}
&M_L(z)=K^*(H_L-z)^{-1}K, \quad z\in
\bbC\backslash \bbR,
\lb{2.2} \\
&\frac{1}{\Im (z)} \Im (M_L(z))=
((H_L-z)^{-1}K)^*(H_L-z)^{-1}K\ge0, \quad
z\in \bbC\backslash \bbR,
\lb{2.3}
\end{align}
and study the pair $(H_L, H_0)$ in terms of the
corresponding pair
$(M_L(z), M_0(z))$. In the special case where
$\dim_{\bbC}(\calK)=1$, this
 perturbation
 problem has been studied in detail by
Donoghue \cite{Do65}
and
 later by Simon and Wolf \cite{SW86}
(see also \cite{Si95}). The case
$\dim_{\bbC}(\calK)=n\in \bbN$,
 has recently been treated in depth in  \cite{GT97}.
 In this section we treat the general
case $\dim_{\bbC}(\calK)\in \bbN\cup \{\infty\}$.

Next, let $\{ E_0(\lambda)\}_{\lambda \in \bbR}$ be
the family
of strongly right-continuous orthogonal
spectral projections of $H_0$ in $\calH$ and suppose
 that $K\calK \subseteq \calH$ is a generating subspace
 for $H_0$, that is, one of the following
(equivalent) equations holds in
\begin{hypothesis}\label{h2.1a}
\begin{subequations} \lb{2.5}
\begin{align}
\calH&=\ol{\text{lin.span}\{
(H_0-z)^{-1}Ke_n\in \calH \, \vert\, n \in \calI,
\, z\in \bbC\backslash \bbR\}}
\lb{2.5a} \\
&=\ol{\text{lin.span}\{
E_0(\lambda)Ke_n\in \calH \, \vert \, n \in
\calI, \, \lambda\in \bbR\}},
\lb{2.5b}
\end{align}
\end{subequations}
where $\{e_n\}_{n\in\calI},$ $\calI\subseteq\bbN$
an appropriate index set, represents a complete
orthonormal system in $\calK$.
\end{hypothesis}
Denoting by
$\{ E_L(\lambda)\}_{\lambda\in \bbR}$
the family of strongly right-continuous orthogonal
spectral projections of $H_L$ in $\calH$ one introduces
\begin{equation}
\Omega_L(\lambda)=K^*E_L(\lambda)K, \quad \lambda\in \bbR
\lb{2.6}
\end{equation}
and hence verifies
\begin{align}
M_L(z)&=K^*(H_L-z)^{-1}K=K^*\int_\bbR dE_L(\lambda)
(\lambda-z)^{-1}K \no \\
&=\int_\bbR d\Omega_L(\lambda)(\lambda-z)^{-1},
\quad z\in  \bbC\backslash \bbR,
\lb{2.7}
\end{align}
where the operator Stieltjes integral \eqref{2.7}
converges in
the norm of $\calB(\calK)$ (cf.~ Theorems~I.4.2 and
I.4.8 in \cite{Br71}). Since
$\slim_{z\to i \infty}z(H_L-z)^{-1}=-I_\calH$,
\eqref{2.6} implies
\begin{equation}
\Omega_L(\bbR)=K^*K.
\lb{2.8}
\end{equation}
Moreover, since
$\slim_{\lambda\downarrow -\infty}E_L(\lambda)=0$,
$\slim_{\lambda\uparrow \infty}E_L(\lambda)=I_\calH$,
one infers
\begin{equation}
\slim_{\lambda\downarrow -\infty}\Omega_L(\lambda)=0,
\quad
\slim_{\lambda\uparrow \infty}\Omega_L(\lambda)=K^*K
\lb{2.9}
\end{equation}
and $\{ \Omega_L(\lambda)\}_{\lambda\in \bbR}
\subset\calB(\calK)$ is a family of uniformly
 bounded, nonnegative, nondecreasing, strongly
right-continuous operators from $\calK$ into itself.
Let $\mu_L$ be a $\sigma-$finite control measure on
$\bbR$ defined, for instance, by
\begin{equation}
\mu_L(\lambda)=\sum_{n\in\calI} 2^{-n}
(e_n,\Omega_L(\lambda)e_n)_\calK, \quad \lambda\in\bbR,
\lb{2.10}
\end{equation}
where $\{e_n\}_{n\in\calI}$ denotes a complete
orthonormal
system in $\calK,$ and then introduce
$L^2({\calM}_{\Omega_L}(\mu_L);d\mu_L)$ as in
Section~\ref{s2}, replacing the pair $(\Omega,\mu)$ by
$(\Omega_L,\mu_L)$, etc.~As noted in Section~\ref{s5},
we will actually use the more suggestive notation
$L^2(\bbR,\calK;wd\Omega_L)$ instead of the more
precise $L^2({\calM}_{\Omega_L}(\mu_L);wd\mu_L)$
($w>0$ a
weight function), for the remainder of this section.
Abbreviating $\hatt\calH_L=L^2(\bbR, \calK;
d \Omega_L)$, we
introduce the unitary operator $U_L:\hatt
\calH_L\to\calH$,
as the operator $U$ in Lemma~\ref{l5.7} and define
$\hatt H_L$ in $\hatt{\calH}_L$ by
\begin{equation}
(\hatt H_L \hat f)(\lambda)=\lambda\hat f(\lambda),
\quad
\hat f \in \dom (\hatt H_L)=L^2(\bbR, \calK;
(1+\lambda^2)d\Omega_L).
\lb{2.11}
\end{equation}
\begin{theorem}\label{t2.2}
Assume Hypotheses~\ref{h2.1} and \ref{h2.1a}.~Then
$H_L$ in $\calH$
is unitarily equivalent to $\hatt H_L$ in
$\hatt{\calH}_L,$
\begin{equation}
H_L=U_L\hatt H_LU_L^{-1}.
\lb{2.12}
\end{equation}
The family of strongly right-continuous orthogonal
spectral projections $\{
\hatt E_L(\lambda)\}_{\lambda\in \bbR}$ of
$\hatt H_L$ in $\hatt{\calH}_L$ is given by
\begin{equation}
(\hatt E_L(\lambda)\hat f)(\nu)=
\theta(\lambda-\nu)\hat f(\nu)
\text{ for } \Omega_L-\text{a.e. } \nu\in \bbR,
\quad \hat f \in\hatt{\calH}_L, \,\,
\theta(x)=\begin{cases} 1,&x\geq 0, \\
  0,&x < 0. \end{cases} \lb{2.13}
\end{equation}
\end{theorem}
\begin{proof}Consider
\begin{equation}
\ul e_n=\{\ul e_n(\lambda)=e_n\}
_{\lambda\in\bbR}\in\hatt \calH_L, \quad n\in\calI,
\lb{2.16}
\end{equation}
then
\begin{equation}
U_L\ul e_n=\int_\bbR dE_L(\lambda)Ke_n=Ke_n,
\quad n\in \calI \lb{2.17}
\end{equation}
and
\begin{equation}
((\hatt H_L-z)^{-1}\ul e_n )(\lambda)=
(\lambda-z)^{-1}\ul e_n(\lambda)=(\lambda-z)^{-1}e_n,
\quad n\in \calI, \,\, z\in \bbC\backslash\bbR
\lb{2.18}
\end{equation}
yield
\begin{equation}
U_L(\hatt H_L-z)^{-1}\ul e_n =\int_\bbR
dE_L(\lambda) (\lambda-z)^{-1}Ke_n =(H_L-z)^{-1}Ke_n,
\quad n\in \calI, \, z\in \bbC\backslash\bbR.
\lb{2.19}
\end{equation}
Using the resolvent equation for $H_L$ and $H_0$,
\begin{subequations} \lb{2.20}
\begin{align}
(H_L-z)^{-1}&=(H_0-z)^{-1}-(H_L-z)^{-1}KLK^*(H_0-z)^{-1}
\lb{2.20a} \\
&=(H_0-z)^{-1}-(H_0-z)^{-1}KLK^*(H_L-z)^{-1}, \quad
z\in \bbC\backslash\bbR,
\lb{2.20b}
\end{align}
\end{subequations}
one verifies
\begin{align}
&(I_\calK+LK^*(H_0-z)^{-1}K)
(I_\calK-LK^*(H_L-z)^{-1}K)
\lb{2.21} \\
&=(I_\calK-LK^*(H_L-z)^{-1}K)
(I_\calK+LK^*(H_0-z)^{-1}K)=I_\calK, \quad z\in
\bbC\backslash\bbR
\end{align}
and
\begin{equation}
(H_L-z)^{-1}K=(H_0-z)^{-1}K
(I_\calK+LK^*(H_0-z)^{-1}K)^{-1}, \quad z\in
\bbC\backslash\bbR.
\lb{2.22}
\end{equation}
Since
\begin{equation}
(I_\calK+LK^*(H_0-z)^{-1}K)^{-1}\in\calB(\calK),
\quad z\in \bbC\backslash\bbR
\lb{2.23}
\end{equation}
by \eqref{2.21}, one infers
\begin{equation}
\ran((I_\calK+LK^*(H_0-z)^{-1}K)^{-1})=\calK,
\quad z\in \bbC\backslash\bbR.
\lb{2.24}
\end{equation}
Since by our assumption \eqref{2.5}, finite linear
combinations of $(H_0-z)^{-1}Ke_n$, $n\in \calI$, $z\in
\bbC\backslash\bbR$ are dense in $\calH$,
\eqref{2.22} and \eqref{2.24} then yield the same
assertion for
$(H_L-z)^{-1}Ke_n$. (I.e., \eqref{2.5} is valid
with $H_0$ replaced by any $H_L$.) Since
$U_L$ is unitary by Lemma~\ref{l5.7}, finite
linear combinations of
vectors of the form $(\hatt H_L-z)^{-1}\ul e_n$
(cf. \eqref{2.19}) are also dense in $\hatt{\calH}$.
This fact,
\eqref{2.19}, and the first resolvent equation for
$\hatt H_L$ yield
\begin{align}
&U_L(\hatt H_L-z)^{-1}
U_L^{-1} U_L(\hatt H_L-z')^{-1}
\ul e_n= U_L(\hatt H_L-z)^{-1}U_L^{-1}
(H_L-z')^{-1}Ke_n \no \\
&=(H_L-z)^{-1}(H_L-z')^{-1}Ke_n, \quad n\in\calI,
\,\, z,z'\in\bbC\backslash\bbR.
\lb{2.26}
\end{align}
Since
finite linear combinations of $(H_L-z')^{-1}Ke_n,$
$n\in\calI$ are dense in $\calH$ we get
\begin{equation}
U_L(\hatt H_L-z)^{-1}U_L^{-1}=(H_L-z)^{-1},
\quad z\in\bbC\backslash\bbR
\lb{2.27}
\end{equation}
and hence \eqref{2.12}. Equation
\eqref{2.13} is then obvious from \eqref{5.10} since
$\hatt H_L$ is the
operator of multiplication by
 $\lambda$ in $\hatt{\calH}_L$.
\end{proof}

If $L_\ell$, $\ell=1,2$ are two bounded
self-adjoint operators in
 $\calK$ (with $\calH, \calK, H_0$,
 and $K$ fixed, i.e., independent of $\ell=1,2)$ one
proves the following result
relating  $M_{L_1}(z)$ and $M_{L_2}(z).$
\begin{theorem}\label{t2.3}
Assume Hypothesis \ref{h2.1}. Let $z\in
\bbC\backslash\bbR$ and suppose $H_{L_\ell}$ and
$M_{L_\ell}(z)$ are defined as in \eqref{2.1} and
\eqref{2.2} with
$\calH,\calK, H_0$ and $K$ independent of $\ell=1,2$
and $L_\ell$, $\ell=1,2$
bounded self-adjoint operators in $\calK$. Then
\begin{subequations} \lb{2.28}
\begin{align}
M_{L_2}(z)&=M_{L_1}(z)(I_\calK+(L_2-L_1)
M_{L_1}(z))^{-1}
\lb{2.28a} \\
&=(I_\calK+M_{L_1}(z)(L_2-L_1))^{-1}M_{L_1}(z).
\lb{2.28b}
\end{align}
\end{subequations}
\end{theorem}
\begin{proof}
Using the resolvent equation for
$H_{L_2}$ and $H_{L_1}$,
\begin{subequations} \lb{2.29}
\begin{align}
(H_{L_2}-z)^{-1}&=(H_{L_1}-z)^{-1}-(H_{L_2}-z)^{-1}
K(L_2-L_1)K^*(H_{L_1}-z)^{-1}
\lb{2.29a} \\
&=(H_{L_1}-z)^{-1}-(H_{L_1}-z)^{-1}K(L_2-L_1)K^*
(H_{L_2}-z)^{-1},\lb{2.29b} \\
& \hspace*{7.7cm} z\in\bbC\backslash\bbR \no
\end{align}
\end{subequations}
and applying $K^*$ on the left and $K$ on the right
of both sides of
\eqref{2.29}, results in
\begin{subequations} \lb{2.30}
\begin{align}
K^*(H_{L_1}-z)^{-1}K&=K^*(H_{L_2}-z)^{-1}
K(I+(L_2-L_1)K^*(H_{L_1}-z)^{-1}K)
\lb{2.30a} \\
&=(I+K^*(H_{L_1}-z)^{-1}K(L_2-L_1))K^*(H_{L_2}-z)^{-1}K
\lb{2.30b}
\end{align}
\end{subequations}
and hence in \eqref{2.28}.
\end{proof}

A comparison of \eqref{2.28} and \eqref{1.5}, \eqref{1.6}
then yields
\begin{equation}
A(L_1, L_2)=
\begin{pmatrix}I_\calK &  L_2-L_1\\ 0&I_\calK
\end{pmatrix}
\in\calA(\calK\oplus\calK)
\lb{2.31}
\end{equation}
for the corresponding matrix $A$ in \eqref{1.5},
\eqref{1.6}.

We note that \eqref{2.28} also imply
\begin{subequations} \lb{2.32}
\begin{align}
(L_2-L_1)M_{L_2}(z)-I_\calK&=
-((L_2-L_1)M_{L_1}(z)+I_\calK)^{-1},
\lb{2.32a} \\
M_{L_2}(z)(L_2-L_1)-I_\calK&=
-(M_{L_1}(z)(L_2-L_1)+I_\calK)^{-1}.
\lb{2.32b}
\end{align}
\end{subequations}

If $K\calK$ is not a generating subspace for
$H_0$ (i.e., \eqref{2.5} does not hold) then $\calH$
decomposes into $\calH=\calH_\calK\oplus
\calH_\calK^\bot$, with
\begin{subequations} \lb{2.33}
\begin{align}
\calH_\calK&=\ol{\text{lin.span}
\{(H_0-z)^{-1}Ke_n\in \calH\, \vert\,  n\in \calI, \,
z\in \bbC\backslash\bbR\}}
\lb{2.33a} \\
&=\ol{\text{lin.span}\{E_0(\lambda)Ke_n\in \calH \,
\vert \, n\in \calI, \, \lambda\in \bbR\}}
\lb{2.33b}
\end{align}
\end{subequations}
and $\calH_\calK$, $\calH_\calK^\bot$
both reducing
subspaces for $H_L$ ($\{e_n\}_{n\in\calI}$ a complete
orthonormal system in $\calK$). Moreover, for
all $L_\ell
\in\calB(\calK)$, $\ell=1,2$
self-adjoint,
\begin{equation}
H_{L_1}=H_{L_2} \text{ on } \dom(H_0)\cap
\calH_\calK^\bot
\lb{2.34}
\end{equation}
and
\begin{equation}
H_0=H_{0,\calK}\oplus H_{0,\calK}^\bot, \quad
H_L=H_{L,\calK}\oplus H_{0,\calK}^\bot, \quad
\ran (K)\subseteq\calH_\calK.
\lb{2.35}
\end{equation}
In particular,
\begin{equation}
M_L(z)=K^*(H_L-z)^{-1}K=K^*(H_{L,\calK}-z)^{-1}K,
\quad z\in\bbC\backslash\bbR
\lb{2.36}
\end{equation}
and the $L$-dependent spectral properties of
$H_L$ in $\calH$
are effectively reduced to those of $H_{L,\calK}$
in $\calH_\calK.$

In connection with our choice of $KLK^*$ as a  bounded
self-adjoint perturbation
 of $H_0$, the
following elementary observation might be of interest.
\begin{lemma}\label{l2.4}
Let $V\in \calB(\calH)$ be self-adjoint. Then $V$ and
$\calH$ can be decomposed as
\begin{equation}
V=K_0L_0K_0^*\oplus 0, \quad \calH=
\ol{\ran(V)}\oplus \ker(V),
\lb{2.37}
\end{equation}
where $K_0:\calK\to \calH$,
$L_0=L_0^*\in \calB(\calK)$, and $\calK=\ol{\ran(V)}.$
\end{lemma}
\begin{proof}
Since $\ol{\ran(V)}=\ker(V)^\bot$, consider
$V_0=V\big\vert_{\ol{\ran(V)}}:\calK\to\calK, $
$\calK=\ol{\ran(V)}$. Then
$V_0=V_0^*\in\calB(\calK)$ and $V_0$ admits the
spectral representation $V_0=\int_a^b
dF_0(\lambda)\lambda$ for some $a,b \in \bbR$ and
some family of self-adjoint spectral projections
$\{F_0(\lambda)\}_{\lambda\in \bbR}$ of $V_0$ in $\calK$.
 The decomposition \eqref{2.37} then follows upon
introducing
\begin{equation}
K_0=\vert V_0\vert^{1/2}=\int_a^b dF_0(\lambda) \vert
\lambda\vert^{1/2}, \quad L_0=\text{sgn}(V_0)=
\int_a^b dF_0(\lambda) \text{sgn}(\lambda).
\lb{2.38}
\end{equation}
\end{proof}
In \eqref{2.6}--\eqref{2.9} we showed that every
collection $(H_0, K, L, \calH, \calK)$ gives
 rise to an operator-valued Herglotz function
$M_L(z)=\int_\bbR d\Omega_L(\lambda)(\lambda-z)^{-1}$
with
certain properties recorded in
\eqref{2.8} and \eqref{2.9}. Conversely, introducing
the following class $\calN_1(\calK)$ of
$\calB(\calK)$-valued Herglotz
functions (we use the symbol $\calN_1(\calK)$ in honor
of R. Nevanlinna)
\begin{equation}
\calN_1(\calK)
=\{
M\in \calB (\calK) \,\text{Herglotz}\,\vert \,
 M(z)=\smallint_\bbR d\Omega(\lambda)(\lambda-z)^{-1};
\, 0\leq\Omega(\bbR)\in\calB(\calK) \},
\lb{2.39}
\end{equation}
we shall show in the remainder of this section
that every
element $M$ of $\calN_1(\calK)$ can be realized in
terms of some collection
$(H_0, K, \calH, \calK)$ as in \eqref{2.7}. (The
operator Stieltjes integral in \eqref{2.39} converges
in the norm of $\calB(\calK)$ by Theorem~I.4.2
of \cite{Br71}.) For this purpose we shall use
a version of
 Naimark's dilation theorem \cite{Na40}, \cite{Na43}
as presented in Appendix I of \cite{AG93} and
Appendix I by Brodskii \cite{Br71}.
\begin{theorem}\label{t2.5}
\mbox{\rm (\cite{Br71}, App.~I, \cite{Na40}.)} Suppose
that
$\Omega(\lambda)$, $\lambda\in \bbR$ is
 a strongly right-continuous nondecreasing function
with
values in $\calB (\calK)$, $\calK$ a complex separable
Hilbert space, and assume
$\slim_{\lambda\downarrow -\infty}\Omega(\lambda)=0$.
Then there exists a separable
complex Hilbert space $\calH$, a $K\in
\calB(\calK,\calH)$, and an orthogonal family of strongly
right-continuous spectral projections
$\{E(\lambda) \}_{\lambda\in\bbR}$ in $\calH$ such that
$
\slim_{\lambda\downarrow - \infty}E(\lambda)=0
$,
$
\slim_{\lambda\uparrow \infty}E(\lambda)=I_\calH
$,
\begin{equation}
\Omega(\lambda)=K^*E(\lambda)K, \quad \lambda\in \bbR,
\lb{2.40}
\end{equation}
and
\begin{equation}
\ol{\{E(\lambda)K\xi\in \calH\,\vert \,
\xi\in \calK, \lambda \in \bbR\}
}=\calH.
\lb{2.41}
\end{equation}
Moreover, if for some
$\lambda_1, \lambda_2 \in \bbR$,
$\Omega(\lambda_1)=\Omega(\lambda_2)$, then
$E(\lambda_1)=E(\lambda_2).$
\end{theorem}

The principal realization theorem for Herglotz operators
of the type  \eqref{2.39}
then reads as follows
\begin{theorem}\label{t2.6}

(i) Any $M\in \calN_1(\calK)$ with associated measure
$\Omega$ can be realized in the form
\begin{equation}
M(z)=K^*(H-z)^{-1}K, \quad z\in\bbC\backslash\bbR,
\lb{2.42}
\end{equation}
where $H$ represents a self-adjoint operator in some
separable complex Hilbert
 space $\calH$, $K\in \calB(\calK,\calH)$, and
\begin{equation}
\Omega(\bbR)=K^*K. \lb{2.43}
\end{equation}
(ii) Suppose $M_\ell\in\calN_1(\calK)$ with
corresponding measures
$\Omega_\ell$, $\ell=1,2$ and $M_1\ne M_2$. Then $M_1$
and $M_2$ can be realized as
\begin{equation}
M_\ell(z)=K^*(H_{L_\ell}-z)^{-1}K, \quad z\in
\bbC\backslash\bbR,
\lb{2.44}
\end{equation}
where $H_{L_\ell}$, $\ell=1,2$ are self-adjoint
perturbations
of one and the same self-adjoint operator
$H_0$ in some separable complex Hilbert space
$\calH$
\begin{equation}
H_{L_\ell}=H_0+KL_\ell K^*, \quad \ell=1,2
\lb{2.45}
\end{equation}
for some $L_\ell=L_\ell^*\in\calB(\calK)$,
$\ell=1,2$ and some
$K\in \calB(\calK,\calH)$ if and only if the
following two conditions hold:
\begin{equation}
\Omega_1(\bbR)=K^*K=\Omega_2(\bbR),
\lb{2.46}
\end{equation}
and for all $z\in \bbC\backslash\bbR$,
\begin{equation}
M_2(z)=M_1(z)(I_\calK+(L_2-L_1)M_1(z))^{-1}.
\lb{2.47}
\end{equation}
\end{theorem}
\begin{proof}
Applying Naimark's dilation theorem, Theorem \ref{t2.5},
to
$\Omega(\lambda)$, $\lambda\in \bbR, $ (assuming
$\slim_{\lambda\downarrow - \infty} \Omega(\lambda)=0$
without loss of generality), yields
$\Omega(\lambda)=K^*E(\lambda)K$, $\lambda\in \bbR$
and introducing the self-adjoint
 operator $H=\int_\bbR dE(\lambda)\lambda$  in $\calH$
then proves \eqref{2.42}.
The normalization condition \eqref{2.43} then follows
as discussed in \eqref{2.6}--\eqref{2.8}.
 In exactly the same manner one proves the necessity
of the normalization
\eqref{2.46}. The necessity of \eqref{2.47} was proven
in Theorem \ref{t2.3}.
In order to prove sufficiency of
\eqref{2.46} and \eqref{2.47}
 for \eqref{2.44} and \eqref{2.45} to hold, we argue
as follows. Suppose
$\slim_{\lambda\downarrow -\infty}\Omega_1(\lambda)=0$
(otherwise, replace
$\Omega_1(\lambda)$ by
$\Omega_1(\lambda)-\slim_{\nu\downarrow -\infty}
\Omega_1(\nu)$) and
represent $M_1(z)$ according to part (i) by
\begin{equation}
M_1(z)=K^*(H_1-z)^{-1}K, \quad z\in\bbC\backslash\bbR
\lb{2.48}
\end{equation}
applying Naimark's dilation theorem and
Theorem \ref{t2.5}. Define
\begin{equation}
H_0=H_1-KL_1K^*, \quad\dom(H_0)=\dom(H_1)
\lb{2.49}
\end{equation}
for some $L_1=L_1^*\in \calB(\calK)$. Next, use
$L_2=L_2^*\in \calB(\calK)$ in \eqref{2.47} to define
\begin{equation}
H_2=H_0+KL_2K^*, \quad\dom(H_2)=\dom(H_0)
\lb{2.50}
\end{equation}
and
\begin{equation}
M_{L_2}(z)=K^*(H_2-z)^{-1}K.
\lb{2.51}
\end{equation}
By Theorem \ref{t2.3},
\begin{equation}
M_{L_2}(z)=M_{L_1}(z)(I_\calK+(L_2-L_1)M_1(z))^{-1}
=M_2(z), \quad z\in\bbC\backslash\bbR \lb{2.52}
\end{equation}
and the proof is complete.
\end{proof}

For a variety of results related to realization
theorems of
Herglotz operators we refer, for instance, to \cite{BGK79}
and the literature cited therein. Fundamental results on
nontangential boundary values of $M_L(z)$ as
$z\to x\in\bbR$,
under various conditions on $K$, can be found in
\cite{Na89}--\cite{Na96}.~Additional results on
operators of the type $M_L(z)$ (including cases where $K$
is a suitable unbounded operator) can be found, for
instance,
in \cite{AL95}, {\cite{Ma92b}, \cite{MS96} and the
references
therein.

\section{On Self-Adjoint Extensions of Symmetric
Operators} \lb{s3}

In this section we consider self-adjoint extensions
$H$ of densely
defined closed symmetric operators $\dot H$ with
deficiency indices $(k,k)$,
$k\in\bbN\cup\{\infty\}$. We revisit Krein's formula
relating self-adjoint
extensions of $\dot H$, introduce the corresponding
operator-valued Weyl
 $m$-functions and their linear fractional
transformations, study a model for the pair
$(\dot H, H)$,
 and consider Friedrichs $H_F$ and Krein
extensions $H_K$ of $\dot H$ in the case where
$\dot H$ is bounded from below.

In the special case $k=1$, detailed investigation of this
 type were undertaken by Donoghue \cite{Do65}.
 The case $k\in \bbN$ was recently
discussed in depth in \cite{GT97} (we also refer to
\cite{HKS98} for another comprehensive treatment of
this subject).~Here we treat the general situation
$k\in \bbN\cup\{\infty\}$ utilizing recent results
in \cite{GMT97}.

We start with a bit of notation
and then recall some pertinent results of
\cite{GMT97}. Let
$\calH$ be a separable complex Hilbert space and
$\dot H:\dom(\dot H)\to \calH$,
$\ol{\dom (\dot H)}=\calH$ a densely defined closed
symmetric linear
operator with equal deficiency indices
$\text{def} (\dot H)=(k,k), k\in \bbN\cup
\{\infty \}.$ The deficiency subspaces
$\calN_{\pm}$ of $\dot H$ are given by
\begin{equation}
\calN_{\pm}=\ker({\dot H}^*\mp i), \quad
\dim_\bbC (\calN_{\pm})=k
\lb{3.1}
\end{equation}
and for any self-adjoint extension $H$ of
$\dot H$ in $\calH$ ,
the corresponding Cayley transform $C_H$ in
$\calH$ is defined by
\begin{equation}
C_H=(H+i)(H-i)^{-1},
\lb{3.2}
\end{equation}
implying
\begin{equation}
C_H\calN_-=\calN_+.
\lb{3.3}
\end{equation}
Two self-adjoint extensions $H_1$ and $H_2$
of $\dot H$ are
called {\it  relatively prime} (w.r.t.
$\dot H$) if
$\dom (H_1)\cap\dom(H_2)=\dom(\dot H).$
Associated with $H_1$ and $H_2$ we
introduce $P_{1,2}(z)\in\calB(\calH)$ by
\begin{align}
&P_{1,2}(z)= (H_1-z)(H_1-i)^{-1} ((H_2-z)^{-1}
-(H_1-z)^{-1})
(H_1-z)(H_1+i)^{-1}, \no \\
&\hspace*{8.2cm} z\in \rho(H_1)\cap\rho(H_2).
\label{3.4}
\end{align}
We refer to Lemma 2 of \cite{GMT97} and \cite{Sa65}
for a detailed discussion of
$P_{1,2}(z)$. Here we only mention the
following properties of
$P_{1,2}(z)$,  $z\in \rho(H_1)\cap\rho(H_2)$,
\begin{align}
&P_{1,2}(z)\big|_{\calN_+^\perp} =0,
\quad P_{1,2}(z)\calN_+
\subseteq \calN_+, \label{3.4a} \\
&\ol{\ran (P_{1,2}(i))}= \calN_+,
\quad\ran (P_{1,2}(z)\big|_{\calN_+})\text{ is
independent of }z\in \rho(A_1)\cap\rho(A_2),
\label{3.4b}\\
&P_{1,2}(i)\big|_{\calN_+} =(i/2)
(I-C_{H_2}C_{H_1}^{-1})\big|_{\calN_+}
=(i/2)(I_{\calN_+}+e^{-2i\alpha_{1,2}})
\lb{3.4c}
\end{align}
for some self-adjoint (possibly unbounded)
operator $\alpha_{1,2} $ in $\calN_+$.

Next, given a self-adjoint extension $H$ of
$\dot H$ and a closed
linear subspace $\calN$ of $\calN_+$,
$\calN\subseteq \calN_+,$
the Weyl-Titchmarsh operator $M_{H,\calN}(z)
\in\calB(\calN)$ associated with the pair
$(H,\calN)$  is defined by
\begin{align}
M_{H,\calN}(z)&=P_\calN (zH+I_\calH)(H-z)^{-1}
P_\calN\big\vert_\calN \no \\
&=zI_\calN+(1+z^2)P_\calN(H-z)^{-1}
P_\calN\big\vert_\calN\,, \quad  z\in \bbC\backslash \bbR,
\lb{3.5}
\end{align}
with $I_\calN$ the identity operator in $\calN$ and
$P_\calN$ the orthogonal projection in $\calH$ onto
$\calN$.

One verifies (cf. Lemma 4 in \cite{GMT97}) for
$H_1$ and $H_2$ relatively prime w.r.t. $\dot H$,
\begin{subequations} \lb{3.6}
\begin{align}
(P_{1,2}(z)\big\vert_{\calN_+})^{-1}
&=(P_{1,2}(i)\big\vert_{\calN_+})^{-1}
-(z-i) P_{\calN_+}(H_1+i)(H_1-z)^{-1}
P_{\calN_+} \lb{3.6a} \\
&=\tan (\alpha_{1,2})-M_{H_1, \calN_+}(z),
\quad z\in \rho(H_1),
\lb{3.6b}
\end{align}
\end{subequations}
where
\begin{equation}
C_{H_2}C_{H_1}^{-1}\big|_{\calN_+} =
-e^{-2i\alpha_{1,2}}.
\lb{3.7}
\end{equation}
Following Saakjan \cite{Sa65} (in a version
presented in Theorem 5 and Corollary 6 in \cite{GMT97}),
Krein's formula then can be summarized as follows.
\begin{theorem}\label{t3.1} \mbox{\rm (\cite{GMT97},
\cite{Sa65}.)} Let $H_1$ and $H_2$ be self-adjoint
extensions of $\dot H$ and $z\in \rho(H_1)\cup\rho(H_2)$.
Then
\begin{align}
(H_2-z)^{-1}&=(H_1-z)^{-1}+(H_1-i)(H_1-z)^{-1}P_{1,2}(z)
(H_1+i)(H_1-z)^{-1} \label{3.8a} \\
&=(H_1-z)^{-1}+(H_1-i)(H_1-z)^{-1} P_{\calN_{1,2,+}}
\label{3.8b} \times \\
&\quad \times (\tan (\alpha_{\calN_{1,2,+}})- M_{H_1,
\calN_{1,2,+}}(z))^{-1}
P_{\calN_{1,2,+}}(H_1+i)(H_1-z)^{-1}, \no
\end{align}
where
\begin{align}
&\calN_{1,2,+}=\ker ((H_1\big|_{\calD (H_1)
\cap\calD (H_2)} )^* -i), \label{3.9} \\
&e^{-2i\alpha_{\calN_{1,2,+}}}=-C_{H_2}C_{H_1}^{-1}
\big|_{\calN_{1,2,+}}, \label{3.10}
\end{align}
and
\begin{equation}
P_{1,2}(i)\big\vert_{\calN_{1,2,+}}=(i/2)
( I-C_{H_2}C_{H_1}^{-1})\big\vert_{\calN_{1,2,+}}.
\lb{3.11}
\end{equation}
\end{theorem}

Next we recall that $M_{H,\calN}$ and hence
$P_{1,2}(z)\big\vert_{\calN_+}$
 and $-(P_{1,2}(z)\big\vert_{\calN_+})^{-1}$
(cf.~\eqref{3.6}), if the latter
 exists, are operator-valued Herglotz functions.

\begin{theorem}\label{t3.2}
Let $H$ be a self-adjoint extension of  $\dot H$
with orthogonal family of spectral
projections $\{E_H(\lambda)\}_{\lambda\in \bbR}$,
$\calN$ a closed subspace of $\calN_+$.
Then the Weyl-Titchmarch operator
$M_{H,\calN}(z)$ is analytic for $z\in
\bbC\backslash\bbR$ and
\begin{equation}\label{3.12}
\Im(z) \Im(M_{H,\calN} (z)) \geq
(\max (1,|z|^2)+|\Re(z)|)^{-1},
\quad z\in \bbC\backslash \bbR.
\end{equation}
In particular, $M_{H,\calN}(z)$ is a
$\calB(\calN)$-valued Herglotz function and
admits the representation valid in the strong
operator topology of $\calN$,
\begin{equation}
M_{H,\calN}(z)=
\int_\bbR
d\Omega_{H,\calN}(\lambda)((\lambda-z)^{-1}-
\lambda(1+\lambda^2)^{-1}), \quad
z\in\bbC\backslash\bbR,
\lb{3.13}
\end{equation}
where
\begin{align}
&\Omega_{H,\calN}(\lambda)=(1+\lambda^2)
(P_\calN E_H(\lambda)P_\calN\big\vert_\calN),
\lb{3.14} \\
&\int_\bbR d\Omega_{H,\calN}(\lambda)
(1+\lambda^2)^{-1}=I_\calN,
\lb{3.15} \\
&\int_\bbR d(\xi,\Omega_{H,\calN}
(\lambda)\xi)_\calH=\infty \text{ for all }
 \xi\in \calN\backslash\{0\}.
\lb{3.16}
\end{align}
\end{theorem}
\begin{proof}
 \eqref{3.13} has been derived in Lemma 7 of
\cite{GMT97},
 hence we confine ourselves to a few hints. An
explicit computation yields
\begin{align}
\Im (z) \Im (M_{H,\calN}(z))&=
P_\calN(I_\calH+H^2)^{1/2}((H-\Re (z))^2+
\Im (z))^2)^{-1} \no \\
&\quad \times (I_\calH+H^2)^{1/2}P_\calN\big\vert_\calN,
\quad \, z\in \bbC\backslash\bbR.
\lb{3.17}
\end{align}
Together with
\begin{equation}
\frac{1+\lambda^2}{(\lambda-\Re(z))^2+(\Im(z))^2}
\ge\frac{1}
{\text{max}(1,\vert z \vert^2)+\vert\Re (z) \vert}
\lb{3.18}
\end{equation}
and the Rayleigh-Ritz argument this yields
\eqref{3.12}. The representation \eqref{3.13}
and the fact \eqref{3.14}
 follow from \eqref{3.5} and
$(H-z)^{-1}\xi=\int_\bbR d(E_H(\lambda)\xi)
(\lambda-z)^{-1}, $ $\xi\in \calH$. \eqref{3.15}
 then follows from
\begin{align}
\int_\bbR d(\Omega_{H,\calN}(\lambda)\xi)
(1+\lambda^2)^{-1}&=
\int_\bbR d (P_\calN E_H(\lambda)\xi)=P_\calN
\int_\bbR d (E_H(\lambda)\xi) \no \\
&=P_\calN \xi=\xi \text{ for all } \xi\in \calN.
\lb{3.19}
\end{align}
Finally,
\begin{equation}
\int_\bbR d(\xi,\Omega_{H,\calN}
(\lambda)\xi)_\calH=
\int_\bbR d(\xi, E_H(\lambda)\xi)_\calH (1+\lambda^2)=
\infty \text{ for all } \xi\in
\calN\backslash\{0\}
\lb{3.20}
\end{equation}
since $\calN\subseteq\calN_+$ and
$\calN_+\cap\dom(H)=\{0\}$ by von Neumann's
formula
\begin{equation}
\dom (H)=\dom (\dot H)\dot + \calN_+\dot+
(-C_H)^{-1}\calN_+.
\lb{3.21}
\end{equation}
\end{proof}

We also recall without proof the principal
result of \cite{GMT97}, the linear fractional
transformation relating the Weyl-Titchmarch
operators associated with different self-adjoint
extensions of $\dot H$.

\begin{theorem} \mbox{\rm (\cite{GMT97}.)}
\label{t3.3}
Let $H_1$ and $H_2$ be self-adjoint extensions
of $\dot H$ and $z\in\rho(H_1)\cap\rho(H_2)$. Then
\begin{align}
M_{H_2, \calN_+}(z)&=(P_{1,2}(i)\big|_{\calN_+}
+(I_{\calN_+}+i
P_{1,2}(i)\big|_{\calN_+})M_{H_1, \calN_+}(z))
\times \no \\
& \quad \times((I_{\calN_+}+
iP_{1,2}(i)\big|_{\calN_+})-
P_{1,2}(i)\big|_{\calN_+}M_{H_1, \calN_+}(z))^{-1},
\label{3.22}
\end{align}
where
\begin{align}
P_{1,2}(i)\big\vert_{\calN_+}&=(i/2)
(I_\calH -C_{H_2}C_{H_1}^{-1})\big\vert_{\calN_+},
\lb{3.22a} \\
I_{\calN_+}+iP_{1,2}(i)\big\vert_{\calN_+}&=(1/2)
(I_\calH +C_{H_2}C_{H_1}^{-1})\big\vert_{\calN_+}.
\lb{3.23}
\end{align}
Introducing
\begin{equation}\lb{3.24}
e^{-2i\alpha_{1,2}}=
-C_{H_2}C_{H_1}^{-1}\big|_{\calN_+},
\end{equation}
\eqref{3.22} can be rewritten as
\begin{align}
M_{H_2, \calN_+}(z)
&=e^{-i\alpha_{1,2}}(\cos (\alpha_{1,2})
+\sin (\alpha_{1,2}) M_{H_1, \calN_+}(z))
\times \no \\
& \quad \times (\sin (\alpha_{1,2})-
\cos (\alpha_{1,2})M_{H_1,
\calN_+}(z))^{-1}e^{i\alpha_{1,2}}. \label{3.25}
\end{align}
\end{theorem}

A comparison of \eqref{3.25} and \eqref{1.5}, \eqref{1.6}
then yields
\begin{equation}
A(\alpha_{1,2})=
\begin{pmatrix}
e^{-i\alpha_{1,2}}\sin(\alpha_{1,2}) &
  -e^{-i\alpha_{1,2}}\cos(\alpha_{1,2})\\
 e^{-i\alpha_{1,2}}\cos(\alpha_{1,2}) &
e^{-i\alpha_{1,2}}\sin(\alpha_{1,2})
\end{pmatrix}
\in\calA (\calK\oplus\calK)
\lb{3.25a}
\end{equation}
for the corresponding matrix $A$ in \eqref{1.5},
\eqref{1.6}.

Weyl operators of the type $M_{H,\calN}(z)$ have
attracted considerable attention in the literature.
The interested reader can find a variety of additional
results, for instance, in \cite{Bu97},
\cite{DM87}--\cite{DMT88}, \cite{KO77}, \cite{KO78},
\cite{Ma92a}, \cite{Ma92b}, \cite{Re98}.

Next we will prepare some material that eventually
will lead to a model for the pair
$(\dot H, H).$ Let $\calN$ be a separable complex
Hilbert space, $\{u_n\}_{n\in\calI},$
$\calI\subseteq\bbN$ a complete
orthonormal system in $\calN,$
$\{\wti \Omega (\lambda)\}_{\lambda\in \bbR}$ a
family of strongly right-continuous
nondecreasing $\calB(\calN)$-valued functions
normalized by
\begin{equation}
\wti\Omega (\bbR)=I_\calK,
\lb{3.26}
\end{equation}
with the property
\begin{equation}
\int_\bbR d(\xi,\wti\Omega (\lambda)\xi)_\calN
(1+\lambda^2)=\infty
 \text{ for all }  \xi\in \calN\backslash\{0\}.
\lb{3.27}
\end{equation}
Introducing the control measure $\ti\mu(B)=
\sum_{n\in\calI} 2^{-n}(u_n,\wti\Omega(B)u_n)_\calN,$
$B\in\Sigma$, and $\ul \Lambda$ as in
Theorem~\ref{t5.5}, we may define $L^p(\bbR,
\calN;wd\wti\Omega),$ $p\ge 1,$ $w\geq0$ a weight
function, as in Section~\ref{s5}. Of special
importance in this section are weight functions
of the type $w_r (\lambda)=(1+\lambda^2)^r,$
$r\in\bbR,$ $\lambda\in\bbR$.
In particular, introducing
\begin{equation}
\Omega(B)=\int_B (1+\lambda^2)d\ti\mu(\lambda)
\f{d\wti\Omega}{d\ti\mu}(\lambda), \quad
B\in\Sigma, \lb{3.27a}
\end{equation}
we abbreviate $\hatt{\calH}=L^2(\bbR, \calN;
d\Omega)$ and define the self-adjoint operator
$\hatt H$ in $\hatt{\calH}$,
\begin{equation}
(\hatt H\hat f)(\lambda)=\lambda\hat f(\lambda),
\quad \hat f \in \dom (\hatt H)=
L^2(\bbR, \calN; (1+\lambda^2) d\Omega), \lb{3.28}
\end{equation}
with corresponding family of strongly
right-continuous orthogonal
 spectral projections
\begin{equation}
(E_{\hat H}(\lambda)\hat f)(\nu)=
\theta(\lambda-\nu)\hat f (\nu)
\text{ for } \Omega-\text{a.e. } \nu
\in \bbR, \quad \hat f \in \hatt{\calH}. \lb{3.29}
\end{equation}
Associated with $\hatt H$ we consider the
linear operator $\hatt{\dot H}$ in
$\hatt{\calH}$ defined as
the following restriction of $\hatt H$
\begin{align}
&\dom (\hatt{\dot H})=\{\hat f\in
\dom (\hatt H) \,\vert  \smallint_\bbR
(1+\lambda^2)d\ti\mu(\lambda) (\ul \xi,
\hat f (\lambda))_{\calN_\lambda}=0
 \text{ for all }\ul \xi\in\ul \Lambda(\calN) \}, \no \\
&\hatt{\dot H}=\hatt H
\big\vert_{\dom(\hat{\dot H})}. \lb{3.30}
\end{align}
(The integral in \eqref{3.30} is well-defined,
see the proof of Theorem \ref{t3.4}
 below.) Here we used the notation introduced
in the proof of Theorem~\ref{t5.5},
\begin{equation}
\ul \xi=\ul \Lambda \xi =\{\ul \xi(\lambda)
=\xi\}_{\lambda\in\bbR}. \lb{3.30a}
\end{equation}

Moreover, introducing the scale of Hilbert spaces
${\hatt{\calH}}_{2r}=L^2(\bbR, \calN;
(1+\lambda^2)^r d\Omega)$, $r\in \bbR$,
$\hatt{\calH}_0=\hatt{\calH}$,
we consider the unitary operator
$R$ from $\hatt{\calH}_2$ to
$\hatt{\calH}_{-2},$
\begin{align}
&R:\hatt{\calH}_2
\longrightarrow \hatt{\calH}_{-2},\quad \hat f
\longrightarrow(1+\lambda^2)\hat f,
\lb{3.31} \\
&(\hat f, \hat g)_{\hat{\calH}_2}=
(\hat f, R \hat g)_{\hat{\calH}}=
(R\hat f, \hat g)_{\hat{\calH}}=
(R\hat f, R\hat g)_{{\hat{\calH}}_{-2}},
\quad \hat f, \hat g \in \hatt{\calH} _2,
\lb{3.32} \\
&(\hat u, \ul v)_{{\hat{\calH}}_{-2}}=
(\hat u, R^{-1} \ul v)_{{\hat{\calH}}_{2}}=
(R^{-1}\hat u, \ul v)_{\hat{\calH}}=
(R^{-1}\hat u, R^{-1}\ul v)_{\hat{\calH}_2},
\quad \hat u, \ul v
\in {\hatt{\calH}}_{-2}.
\lb{3.33}
\end{align}
In particular,
\begin{equation}
\ul {\ul \Lambda}(\calN)\subset\hatt\calH, \quad
\ul \Lambda(\calN)\subset \hatt{\calH}_{-2},
\quad \ul \xi \in\ul \Lambda(\calN)\backslash\{0\}
\Rightarrow \ul \xi\not\in \hatt{\calH}
\lb{3.34}
\end{equation}
(cf.~\eqref{5.31} and \eqref{3.26}--\eqref{3.27a}).
\begin{theorem}\label{t3.4}
The operator $\hatt{\dot H}$ in \eqref{3.30} is
densely defined symmetric and
closed in $\hatt{\calH}$. Its deficiency indices
are given by
\begin{equation}
\rm{def} (\hatt{\dot H})=(k,k), \quad k=
\dim_\bbC(\calN)\in \bbN\cup\{\infty\},
\lb{3.35}
\end{equation}
and
\begin{equation}
\ker (\hatt{\dot  H}^*-z)=\ol{\text{lin.span}
\{\{(\lambda-z)^{-1}e_n\}_{\lambda\in\bbR}
\in \hatt{\calH} \, \vert \, n\in \calI \}}, \quad
z\in\bbC\backslash\bbR. \lb{3.36}
\end{equation}
\end{theorem}
\begin{proof}
Writing
 $\vert\vert \hat f(\lambda)\vert\vert
_{\calN_\lambda}=
(1+\lambda^2)^{-1/2}(1+\lambda^2)^{1/2}
\vert\vert \hat f(\lambda)\vert\vert
_{\calN_\lambda}$
one infers that
$\hat f \in L^1(\bbR, \calN; d\Omega)$ for
$\hat f \in \hatt{\calH}_2$.
Thus the integral in \eqref{3.30} and hence
$\dom(\hatt{\dot H})$ is well-defined. As a
restriction of $\hatt H$,
$\hatt{\dot H}$ is clearly symmetric. By
\eqref{3.30} and \eqref{3.31}--\eqref{3.33} one infers
\begin{equation}
\dom (\hatt H)=
\dom(\hatt{\dot H})=\hatt{\calH}_2
\ominus_{\hat{\calH}_2}R^{-1}\ul \Lambda(\calN),
\lb{3.37}
\end{equation}
where, in obvious notation,
$\ominus_{\hat{\calH}_2}$
indicates the orthogonal complement in
$\hatt{\calH}_2$. Thus $\hatt{\dot H}$ has a closed
graph.

Next, to prove that $\hatt{\dot H}$ is densely
defined in $\hatt{\calH}$, suppose there is
a $\hat g\in \hatt{\calH}$ such that
$\hat g\bot \dom (\hatt{\dot H})$. Then
\begin{equation}
0=(\hat f, \hat g)_{\hat{\calH}}=(\hat f,R^{-1}
 \hat g)_{\hat{\calH}_2} \text{ for all }
 \hat f \in \dom(\hatt{\dot H})
\lb{3.38}
\end{equation}
and hence $R^{-1}\hat g\in R^{-1}\ul \Lambda(\calN)$,
that is, there is an $\xi\in\calN$ such that
$\hat g=\ul \Lambda \xi$
$\Omega-$a.e.~by \eqref{3.37}. Since
$\ul \Lambda \xi\in\ul \Lambda(\calN)\backslash\{0\}$
implies $\ul \Lambda \xi\not\in\hatt \calH$
by \eqref{3.34}, $\hat g\in \hatt{\calH}$
 if and only if $\ul \Lambda \xi=\hat g=0$. Finally,
since $\hatt H$ is self-adjoint,
$\ran (\hatt H-z)=\hatt{\calH}$ for all $z\in
\bbC\backslash\bbR$, and
$(\hatt H \pm i):\hatt{\calH}_2 \to\hatt{\calH}$
is unitary,
\begin{equation}
((\hatt H \pm i) \hat f,(\hatt H \pm i)
\hat g)_{\hat{\calH}}= \int_\bbR (1+\lambda^2)^2
d\ti\mu(\lambda)(\hat f(\lambda),\hat
g(\lambda))_{\calN_\lambda}
=(\hat f, \hat g)_{\hat{\calH}_2}, \quad
\hat f, \hat g \in \hatt{\calH}_2. \lb{3.39}
\end{equation}
Thus \eqref{3.37} and \eqref{3.38} yield
\begin{align}
\hatt{\calH}&=(\hatt H\pm i)\hatt{\calH}_2=
(\hatt H\pm i)(\dom (\hatt{\dot H})
\oplus_{\hat{\calH}_2}R^{-1}\ul \Lambda(\calN)) \no \\
&=(\hatt{\dot H}\pm i)\dom (\hatt{\dot H})
\oplus_{\hat{\calH}} \ol{\text{lin.span}
\{\{(\lambda\pm i)(1+\lambda^2)^{-1}u_n\}
_{\lambda\in\bbR} \in
\hatt{\calH}\, \vert \, n\in \calI \}} \no \\
&=\ran(\hatt{\dot H}\pm i)\oplus_{\hat{\calH}}
\ol{\text{lin.span}\{\{(\lambda\mp i)^{-1}u_n\}
_{\lambda\in\bbR}\in \hatt{\calH}
\, \vert \, n\in \calI\}} \lb{3.40}
\end{align}
and hence
\begin{equation}
\ker(\hatt{\dot H}^*\mp i)=\ol{\text{lin.span}
\{(\{\lambda\mp i)^{-1}u_n\}_{\lambda\in\bbR}
\in \hatt{\calH} \, \vert \, n\in \calI \}}.
\lb{3.41}
\end{equation}
Since $(\lambda-z)^{-1}\xi=(\lambda -i)^{-1}\xi+(z-i)
(\lambda-z)^{-1}(\lambda-i)^{-1}\xi$, with
$\{(\lambda -z)^{-1}(\lambda -i)^{-1}\xi\}
_{\lambda\in\bbR}\in
\hatt{\calH}_2=\dom (\hatt H)$ for all $\xi\in \calN$,
$z\in \bbC\backslash\bbR$, \eqref{3.41} yields
\eqref{3.36}.
\end{proof}
\begin{lemma}\label{l3.5}
Let $\dot H$ be a densely defined linear closed
symmetric operator
in a separable
complex  Hilbert space $\calH$ with deficiency
indices $(k,k)$,
 $k\in \bbN\cup \{\infty\}$. Then
 $\calH$ decomposes into the direct orthogonal sum
\begin{equation}
\calH=\calH_0\oplus \calH_0^\bot, \quad
\ker ({\dot H}^*-i)\subset \calH_0, \quad
 z\in \bbC\backslash\bbR,
\lb{3.42}
\end{equation}
where
$\calH_0$ and $ \calH_0^\bot$ are invariant
subspaces for
all self-adjoint extensions of $\dot H$, that is,
\begin{equation}
(H-z)^{-1}\calH_0\subseteq \calH_0 , \quad
(H-z)^{-1}\calH_0^\bot\subseteq \calH_0^\bot,
\quad
 z\in \bbC\backslash\bbR,
\lb{3.43}
\end{equation}
for all self-adjoint extensions $H$ of $\dot H$
in $\calH$. Moreover,
all self-adjoint extensions  $\dot H$ coincide on
$\calH_0^\bot$, that is, if
$\{H_\alpha\}_{\alpha\in \calI}$
($\,\calI$ an appropriate index set)
 denotes the set of all self-adjoint extensions
of $\dot H$, then
 \begin{equation}
H_\alpha=H_{0,\alpha}\oplus H_0^\bot, \quad
\alpha\in \calI \quad  \text{in} \quad
 \calH=\calH_0\oplus\calH_0^\bot,
\lb{3.44}
\end{equation}
where
\begin{equation}
H_0^\bot \text{ is independent of } \alpha\in \calI.
\lb{3.45}
\end{equation}
\end{lemma}
\begin{proof}
Let $H$ be a fixed self-adjoint extension of $\dot H$,
denote $\calN_\pm=\ker({\dot H}^*\mp i),$ and define
\begin{equation}
\calH_H=\ol{\text{lin.span}\{
(H-z)^{-1}u_+\in \calH \, \vert\,
u_+\in \calN_+, \, z \in \bbC\backslash \bbR \}}.
\lb{3.46}
\end{equation}
Since $(H-z_1)^{-1}
(H-z_2)^{-1}=(z_1-z_2)^{-1}((H-z_1)^{-1}-
(H-z_2)^{-1}), $
$\calH_H$ is invariant with respect to $(H-z)^{-1}$,
$(H-z)^{-1}\calH_H\subseteq\calH_H$, and since
$((H-z)^{-1})^*=(H-\bar z)^{-1}$, also $\calH_H^\bot$
is invariant
 under $(H-z)^{-1}$ for all $z\in \bbC\backslash\bbR.$
Since $\wlim_{z\to i \infty}(-z)(H-z)^{-1}f=f$ for all
$f\in \calH$, one concludes
\begin{equation}
\calN_+\subset \calH_H.
\lb{3.47}
\end{equation}
Next, let $v\in\calH_H^\bot$. Then also
\begin{equation}
w=(H-z)^{-1}v\in \calH_H^\bot, \quad z\in
\bbC\backslash\bbR
\lb{3.48}
\end{equation}
and
\begin{equation}
(u_+,v)_{\calH}=(u_+, w)_{\calH}=0,
\quad u_+\in \calN_+.
\lb{3.49}
\end{equation}
Since $w\in \dom(H)$
\begin{equation}
w\notin\calN_\pm
\lb{3.50}
\end{equation}
(otherwise ${\dot H}^*w=\pm iw$ yields $Hw=\pm iw$
which contradicts
 the self-adjointness of $H$). By von Neumann's
formulas
\begin{equation}
\dom({\dot H}^*)=\dom(\dot H)\oplus_{\calH_+}
\calN_+\oplus_{\calH_+}\calN_-,
\lb{3.51}
\end{equation}
where $\oplus_{\calH_+}$ denotes the direct
orthogonal sum in the Hilbert space $\calH_+$
defined by
\begin{equation}
\calH_+=(\dom({\dot H}^*), \,(\cdot, \cdot)_+),
\quad (f,g)_+=({\dot H}^*f,
{\dot H}^* g)_\calH+(f,g)_\calH,
\quad f,g\in \dom({\dot H}^*). \lb{3.52}
\end{equation}
Using \eqref{3.47},  $Hw=zw+v$ (cf. \eqref{3.48}),
\eqref{3.49}, and \eqref{3.52}
one computes
\begin{align}
(u_+,w)_+
&=({\dot H}^*u_+, {\dot H}^* w)_\calH +(u_+, w)_\calH
=-i(u_+, Hw)_\calH +(u_+,w)_\calH \no \\
&=(-iz+1)(u_+,w)_\calH -i(u_+,v)_\calH=0.
\label{3.53}
\end{align}
\eqref{3.50}, \eqref{3.51}, and \eqref{3.53} then
prove $w\in \dom (\dot H)$ and hence
\begin{equation}
Hw=\dot Hw=zw+v.
\lb{3.54}
\end{equation}
If $\wti H$ is any other self-adjoint extension of
$\dot H$, then $w\in \dom (\dot H)$ also yields
\begin{equation}
\wti H w=\dot H w=zw+v
\lb{3.55}
\end{equation}
and hence
\begin{equation}
w=(H-z)^{-1}v=(\wti H-z)^{-1}v, \quad v \in
\calH_H^\bot.
\lb{3.56}
\end{equation}
Thus the resolvents of all self-adjoint extensions of
$\dot H$
 coincide on $\calH_H^\bot$. Moreover,
\begin{equation}
((\wti H-\bar z)^{-1}u_+, v)_\calH=(u_+,
(\wti H - z)^{-1}v)_\calH=
(u_+,w)_\calH=0
\lb{3.57}
\end{equation}
yields
\begin{equation}
(\wti H-z)^{-1}u_+ \,\bot \,\calH_H^\bot, \quad z
\in\bbC\backslash \bbR
\lb{3.58}
\end{equation}
and hence $\calH_{\ti H}\subseteq\calH_H$. By
symmetry
in $H$ and $\wti H$ ,
$\calH_{\ti H}=\calH_H=\calH_0$ completing the proof.
\end{proof}
In the following we call a densely defined closed
symmetric operator $\dot H$
with deficiency indices $(k,k)$, $k\in
\bbN\cup\{\infty\}$
 {\it prime} if
$\calH_0^\bot=\{0\}$ in the decomposition \eqref{3.42}.

\vspace*{2mm}

Given these preliminaries we can now discuss a model
for the pair $(\dot H, H).$
\begin{theorem}\label{t3.6}
Let $\dot H$ be a densely defined closed prime symmetric
operator
in a separable complex
Hilbert space $\calH$. Assume $H$ to be a self-adjoint
extension of $\dot H$ in $\calH$ with
$\{ E_H(\lambda) \}_{\lambda\in \bbR}$ the
associated family of strongly  right-continuous
orthogonal spectral projections
of $H$ and define
the unitary operator
$\wti U:\hatt{\calH}=L^2(\bbR, \calN_+;
d \Omega_{H,\calN_+})\to \calH$ as the operator
$\wti U$ in Lemma~\ref{l5.8}, where
\begin{equation}
\Omega_{H,\calN_+}(\lambda)=(1+\lambda^2)
(P_{\calN_+}E_H(\lambda)P_{\calN_+}
\big\vert_{\calN_+}),
\lb{3.60}
\end{equation}
with $P_{\calN_+}$
 the orthogonal projection onto $\calN_+=
\ker ( {\dot H}^*-i)$.
Then the pair $(\dot H, H)$ is unitarily
equivalent to the pair $\hatt{\dot H}, \hatt H)$,
\begin{equation}
\dot H=\wti U \hatt{\dot H}{\wti U}^{-1},
\quad H=\wti U\hatt H {\wti U}^{-1},
\lb{3.61}
\end{equation}
where $\hatt H$ and $\hatt{\dot H}$ are defined in
\eqref{3.26}--\eqref{3.30}, and Theorem~\ref{t3.4},
and $\calN$ is identified with $\calN_+$, etc. Moreover,
\begin{equation}
\wti U\hatt{\calN}_+=\calN_+,
\lb{3.62}
\end{equation}
where
\begin{equation}
{\hatt \calN}_+=\ol{\text{lin.span}\{\ul{\ul u}_{+,n}
\in \hatt{\calH} \,\vert \, \ul{\ul u}_{+,n}(\lambda)=
(\lambda-i)^{-1}u_{+,n}, \, \lambda\in\bbR,
\, n\in\calI \}},
\lb{3.63}
\end{equation}
with $\{u_{+,n}\}_{n\in\calI}$ a complete
orthonormal system in ${\calN}_+=\ker({\dot H}^*-i).$
\end{theorem}
\begin{proof}Consider $\ul{\ul u}_{+,n}(\lambda)=
(\lambda-i)^{-1}u_{+,n}$, $n\in\calI,$ then
\begin{equation}
\wti U\ul{\ul u}_{+,n}=\int_\bbR d E_H (\lambda)u_{+,n}
=u_{+,n}, \quad n\in\calI \lb{3.66}
\end{equation}
proves \eqref{3.62}. Moreover,
\begin{equation}
((\hatt H-z)^{-1}\ul{\ul u}_{+,n})(\lambda)=
(\lambda-z)^{-1}(\lambda-i)^{-1}
u_{+,n}, \quad n\in\calI, \, z \in \bbC\backslash \bbR
\lb{3.67}
\end{equation}
yields
\begin{equation}
\wti U(\hatt H-z)^{-1}\ul{\ul u}_{+,n}=\int_\bbR
dE_H(\lambda)(\lambda-z)^{-1}u_{+,n}
=(H-z)^{-1}u_{+,n}, \quad n\in\calI.
\lb{3.68}
\end{equation}
Since by hypothesis $\dot H$ is a prime symmetric
operator, finite
linear combinations of the right-hand side in
\eqref{3.68} are dense in $\calH$. Since $\wti U$ is
unitary, also finite linear
combinations of $(\hatt H-z)^{-1}\ul{\ul u}_{+,n}$ on
the left-hand side of \eqref{3.68} are dense in
 $\hatt{\calH}$. Using the first resolvent
equation one computes from \eqref{3.68}
\begin{align}
\wti U(\hatt H-z)^{-1}{\wti U}^{-1}\wti
U(\hatt H-z')^{-1}\ul{\ul u}_{+,n}
&=\wti U(\hatt H-z)^{-1}{\wti U}^{-1}(H-z')^{-1}
u_{+,n} \no \\
&=(H-z)^{-1}(H-z')^{-1} u_{+,n}.
\lb{3.70}
\end{align}
Since finite linear combinations of the form
$(H-z')^{-1}u_{+,n}$ are dense in $\calH$ we get
\begin{equation}
\wti U(\hatt H-z)^{-1}{\wti U}^{-1}=(H-z)^{-1},
\quad z\in\bbC\backslash \bbR.
\lb{3.71}
\end{equation}
\eqref{3.62} and \eqref{3.71} then yield
$\wti U\hatt{\dot H}{\wti U}^{-1}=\dot H.$
\end{proof}
If $\dot H$ is a densely defined closed non-prime
symmetric operator in $\calH$, then in addition to
\eqref{3.42}, \eqref{3.44}, and \eqref{3.45} one
obtains
\begin{equation}
\dot H={\dot H}_0\oplus H_0^\bot,
\quad \calN_+=\calN_{0,+}\oplus \{ 0 \}
\lb{3.72}
\end{equation}
with respect to the decomposition
$\calH=\calH_0\oplus\calH_0^\bot$.
In particular, the part $H_0^\bot$ of $\dot H$ in
$\calH_0^\bot$ is self-adjoint.
 For any closed linear subspace $\calN$ of
$\calN_+$, $\calN\subseteq \calN_+$,
one then infers $\calN=\calN_0\oplus \{0\}, $
$P_\calN=P_{\calN_0}\oplus 0$ and hence
\begin{equation}
M_{H, \calN}(z)=M_{H_0, \calN_0}(z), \quad z\in
\bbC\backslash \bbR.
\lb{3.73}
\end{equation}
This reduces the $H$-dependent spectral
 properties of the Weyl-Titchmarch operator
effectively to that of $H_0$,
 where $H=H_0\oplus H_0^\bot$ is a self-adjoint
extension of $\dot H$ in $\calH$.

 Next we digress a bit to the special case where
$\dot H\ge 0$
and characterize Friedrichs and Krein extensions,
$H_F$ and $H_K$, of $\dot H$ in $\calH$. Assuming
$\dot H$ to be densely defined in $\calH$ we recall
the definition of $H_F$ and $H_K$ (cf., e.g.,
\cite{AN70}),
\begin{align}
&\dom(H_F^{1/2})=\{f\in\calH\,|\,\text{there is a }
\{f_n\}_{n\in\bbN}\subset\dom(\dot H) \text{ s.t. }
\lim_{n\to\infty}\|f_n-f\|_\calH =0 \no \\
&\hspace*{5.35cm} \text{ and }
\lim_{m,n\to\infty}((f_n-f_m),
\dot H (f_n-f_m))_\calH = 0\}, \no \\
&H_F={\dot H}^*\big|_{\dom({\dot H}^*)\cap
\dom(H_F^{1/2})}, \lb{3.73a} \\
&\dom(H_K)=\{f\in\dom({\dot H}^*)\,|\,\text{there is a }
\{f_n\}_{n\in\bbN}\subset\dom(\dot H) \text{ s.t. } \no \\
&\hspace*{1.8cm} \lim_{n\to\infty}\|\dot Hf_n-
{\dot H}^* f\|_\calH =0
\text{ and } \lim_{m,n\to\infty}((f_n-f_m),\dot H
(f_n-f_m))_\calH = 0\}, \no \\
&H_K={\dot H}^*\big|_{\dom(H_K)}. \lb{3.73b}
\end{align}
Moreover, we recall that
\begin{align}
&\inf\spec(H_F)=\inf\{(g,\dot Hg)_\calH\in\bbR\,|\,
g\in\dom(\dot H),\, \|g\|_\calH=1\}\geq 0, \lb{3.73c} \\
&\inf\spec(H_K)=0, \lb{3.73d}
\end{align}
and
\begin{equation}
0\leq (H_F-\mu)^{-1}\leq (\wti H-\mu)^{-1}\leq
(H_K-\mu)^{-1}, \quad \mu<0 \lb{3.73e}
\end{equation}
for any nonnegative self-adjoint extension $\wti H\geq 0$
of $\dot H$.

Next we discuss a slight refinement of
a result of Krein \cite{Kr47}
(see also \cite{AT82}, \cite{Ts80}, \cite{Ts81}).
We will use an efficient summary of Krein's result
due to Skau \cite{Sk79} (cf.~also \cite{LT77}), which
appears most relevant in our context.
\begin{theorem}\label{t3.7}
Let $\dot H\ge 0$ be a densely defined closed
nonnegative operator in $\calH$ with
 deficiency subspaces $\calN_\pm=
\ker({\dot H}^*\mp i)$.
Suppose $H$ is a self-adjoint extension of $\dot H$
in $\calH $ with corresponding family
 of orthogonal spectral projection
$\{E_H(\lambda)\}_{\lambda\in \bbR}$ and define
\begin{equation}
\Omega_{H,\calN_+} (\lambda)=(1+\lambda^2)
(P_{\calN_+}E_H(\lambda)P_{\calN_+}\big\vert_{\calN_+}).
\lb{3.74}
\end{equation}

\noindent Denote by $H_F$ and $H_K$ the Friedrichs
and Krein extension of $\dot H$, respectively. Then

\noindent (i) $H=H_F$ if and only if $\int_R^\infty
d \vert\vert E_H(\lambda)u_+
\|^2_\calH \lambda=\infty$, or equivalently,
if and only if
$\int_R^\infty d(u_+, \Omega_{H,\calN_+}
(\lambda) u_+)_{\calN_+}\lambda^{-1}=\infty$
 for all $R>0$ and all $u_+\in
\calN_+\backslash\{0\}.$

\noindent (ii) $H=H_K$ if and only if $\int_0^R
d \vert\vert
 E_H(\lambda)u_+\vert\vert^2_\calH \lambda^{-1}
=\infty$, or equivalently,
if and only if
 $\int_0^R d(u_+, \Omega_{H,\calN_+}
(\lambda) u_+)_{\calN_+}\lambda^{-1}=\infty$
 for all $R>0$ and all $u_+\in
\calN_+\backslash\{0\}.$

\noindent (iii) $H=H_F=H_K $ if and only if
$\int_R^\infty d\vert\vert E_H(\lambda)u_+
\vert\vert^2_\calH\lambda=
\int_0^R d\vert\vert E_H(\lambda)u_+
\vert\vert^2_\calH\lambda^{-1}=\infty$ , or
equivalently, if and only if for all $R>0$
and all $u_+\in \calN_+\in
\calN_+\backslash\{0\}$,
$\int_R^\infty d(u_+, \Omega_{H,\calN_+}
(\lambda) u_+)_{\calN_+}\lambda^{-1}=
\int_R^\infty d(u_+, \Omega_{H,\calN_+}
(\lambda) u_+)_{\calN_+}
\lambda^{-1}=\infty$.

\end{theorem}
\begin{proof}
By Lemma \ref{l3.5} and \eqref{3.72}
 we may assume that $\dot H$ is a prime symmetric
operator.
 Moreover, by Theorem \ref{t3.6} we may
identify $(\dot H, H)$ in $\calH$ with the model
pair $(\hatt{\dot H}, \hatt H)$ in
$\hatt{\calH}=L^2(\bbR, \calN_+;
d \Omega_{H,\calN_+}).$ Since by \eqref{3.63},
\begin{equation}
\hatt \calN_+=\ol{\text{lin.span}\{\ul{\ul u}_{+,n}
=\{(\lambda -i)^{-1}u_{+,n}\}_{\lambda\in\bbR}\in
\hatt{\calH} \, \vert \, n\in\calI \}},
\lb{3.75}
\end{equation}
statements (i)--(iii) are reduced to those in
Krein \cite{Kr47}, respectively Skau \cite{Sk79},
who use
$\ker ({\dot H}^*+1)$ instead of $\calN_+=
\ker({\dot H}^*-i)$, by  utilizing the elementary
identity
$(\lambda+1)^{-1}=(\lambda-i)^{-1} -(1+i)
(\lambda +1)^{-1} (\lambda-i)^{-1}$
and the fact that $\{(\lambda+1)^{-1}
(\lambda -i)^{-1}u_{+,n}\}_{\lambda\in\bbR} \in
\hatt {\calH}=
L^2(\bbR,\calN_+; d \Omega_{H,\calN_+})$ for all
$n\in\calI$.
\end{proof}
\begin{corollary} \label{c3.8} \mbox{\rm (\cite{DM91},
\cite{DM95}, \cite{DMT88}, \cite{KO78}, \cite{Ts92}.)}

\noindent (i) $H=H_F$ if and only if
$
\lim_{\lambda\downarrow -\infty}(u_+, M_{H,\calN_+}
(\lambda)u_+)_{\calN_+}=-\infty
$
for all $u_+\in \calN_+\backslash\{0\}.$

\noindent (ii)  $H=H_K$ if and only if
$\lim_{\lambda\uparrow 0}(u_+, M_{H,\calN_+}
(\lambda)u_+)_{\calN_+}=\infty$
for all $u_+\in \calN_+\backslash\{0\}.$

\noindent (iii)  $H=H_F=H_K$ if and only if
$\lim_{\lambda\downarrow -\infty}(u_+,
M_{H,\calN_+}(\lambda)u_+)_{\calN_+}=-\infty$
and
 $\lim_{\lambda\uparrow 0}(u_+, M_{H,\calN_+}
(\lambda)u_+)_{\calN_+}=\infty$
for all $u_+\in \calN_+\backslash \{0\}.$
\end{corollary}
\begin{proof}
Since
\begin{align}
M_{H,\calN_+}(z)&=zI_{\calN_+}
+(1+z^2)P_{\calN_+}(H-z)^{-1}P_{\calN_+}
\big\vert_{\calN_+} \no \\
&=\int_\bbR d \Omega_{H,\calN_+}(\lambda)
((\lambda-z)^{-1}-\lambda(1+\lambda^2)^{-1}),
\quad z\in \bbC \backslash [0, \infty)
\lb{3.76}
\end{align}
by \eqref{3.74}, it suffices to  involve
Theorem~\ref{t3.7}\,(i)--(iii) and the monotone
convergence theorem.
\end{proof}

As a simple illustration we mention the following

\begin{example} \lb{e.3.8a}
Consider the following operator $\dot H$ in
$L^2(\bbR^n;d^nx)$,
\begin{equation}
\dot H=\ol{-\Delta\big|_{C^\infty_0
(\bbR^n\backslash\{0\})}}\geq 0, \quad n=2,3.
\lb{3.76a}
\end{equation}
Then
\begin{equation}
H_F=H_K=-\Delta, \quad \dom(-\Delta)=H^{2,2}(\bbR^2)
\text{ if } n=2 \lb{3.76b}
\end{equation}
is the unique nonnegative self-adjoint
extension of $\dot H$ in $L^2(\bbR^2;d^2x)$ and
\begin{align}
&H_F=-\Delta, \quad \dom(-\Delta)=H^{2,2}(\bbR^3)
\text{ if } n=3, \lb{3.76c} \\
&H_K=Uh_0^NU^{-1}\oplus\bigoplus_{\ell\in\bbN}
Uh_\ell U^{-1} \text{ if } n=3. \lb{3.76d}
\end{align}
Here $H^{p,q}(\bbR^n),$ $p,q\in\bbN$ denote the usual
Sobolev spaces,
\begin{align}
&h_0^N=-\f{d^2}{dr^2}, \quad r>0, \lb{3.76e} \\
&\dom(h_0^N)=\{f\in L^2((0,\infty);dr)\,|\,
f,f'\in AC([0,R]) \text{ for all } R>0; \, f'(0_+)=0;
\no \\
& \hspace*{8.8cm} f''\in L^2((0,\infty);dr)\}, \no \\
&h_\ell=-\f{d^2}{dr^2}+\f{\ell(\ell +1)}{r^2}, \quad r>0,
\, \ell\in\bbN, \lb{3.76f} \\
&\dom(h_\ell)=\{f\in L^2((0,\infty);dr)\,|\,
f,f'\in AC([0,R]) \text{ for all } R>0; \, f(0_+)=0;
\no \\
&\hspace*{5.8cm} -f''+\ell(\ell +1)r^{-2}f\in
L^2((0,\infty);dr)\}, \no
\end{align}
and $U$ denotes the unitary operator,
\begin{equation}
U:L^2((0,\infty);dr)\to L^2((0,\infty);r^2dr),
\quad f(r)\to r^{-1}f(r). \lb{3.76g}
\end{equation}
Equations \eqref{3.76b}--\eqref{3.76d} follow from
Corollary~\ref{c3.8} and the facts
\begin{equation} \lb{3.76h}
(u_+,M_{H_F,\calN_+}(z)u_+)_{L^2(\bbR^n;d^nx)} =
\begin{cases}
-(2/\pi) \ln(z) +2i, & n=2, \\
i(2z)^{1/2} +1, & n=3,
\end{cases}
\end{equation}
and
\begin{equation}
(u_+,M_{H_K,\calN_+}(z)u_+)_{L^2(\bbR^3;d^3x)} =
i(2/z)^{1/2} - 1. \lb{3.76i}
\end{equation}
Here
\begin{equation}
\calN_+=\text{lin.span}\{u_+\}, \quad u_+ (x)=
G_0(i,x,0)/\|G_0(i,\cdot,0)\|_{L^2(\bbR^n;d^nx)},
\,\, x\in\bbR^n\backslash\{0\},
\lb{3.76j}
\end{equation}
where
\begin{equation} \lb{3.76k}
G_0(z,x,y) =
\begin{cases}
\f{i}{4}H_0^{(1)}(z^{1/2}|x-y|), &x\neq y, \, n=2, \\
e^{iz^{1/2}|x-y|}/(4\pi |x-y|), &x\neq y, \, n=3
\end{cases}
\end{equation}
denotes the Green's function of $-\Delta$ on
$H^{2,2}(\bbR^n),$ $n=2,3$ (i.e., the
integral kernel of the resolvent $(-\Delta -z)^{-1}$)
and $H_0^{(1)}(\zeta)$ abbreviates the Hankel function
of the first kind and order zero (cf., \cite{AS72},
Sect.~9.1).
Equation \eqref{3.76h} then immediately follows from
repeated use of the
identity (the first resolvent equation),
\begin{align}
&\int_{\bbR^n} d^nx' G_0(z_1,x,x')G_0(z_2,x',0) =
(z_1-z_2)^{-1}(G_0(z_1,x,0)-G_0(z_2,x,0)), \no \\
&\hspace*{7.4cm} x\neq 0, \, z_1\neq z_2, \, n=2,3
\lb{3.76l}
\end{align}
and its limiting case as $x\to 0$. Finally, \eqref{3.76i}
follows from the following arguments. First one notices
that
$(-(d^2/dr^2)+\nu r^{-2})\big|_{C_0^\infty((0,\infty))}$
is essentially self-adjoint if and only if $\nu\geq 3/4$.
Hence it suffices to consider the restriction of
$\dot H$ to
the centrally symmetric subspace of $L^2(\bbR^3;d^3x)$
corresponding to angular momentum $\ell=0$. But then
it is a
well-known fact (cf.~Lemma~\ref{l4.3}) that the Dirichlet
Donoghue $m$-function
$(u_+,M_{H_F,\calN_+}(z)u_+)_{L^2(\bbR^n;d^nx)}$
corresponding to
\begin{align}
&h_0^D=-\f{d^2}{dr^2}, \quad r>0, \lb{3.76m} \\
&\dom(h_0^N)=\{f\in L^2((0,\infty);dr)\,|\,
f,f'\in AC([0,R]) \text{ for all } R>0; \, f(0_+)=0;
\no \\
& \hspace*{8.7cm} f''\in L^2((0,\infty);dr)\}, \no
\end{align}
and the Neumann Donoghue $m$-function
$(u_+,M_{H_N,\calN_+}(z)u_+)_{L^2(\bbR^n;d^nx)}$
corresponding to $h_0^N$ in \eqref{3.76e} are related
to each other by \eqref{4.25}, with $\alpha=\pi/2,$
$\beta=\pi/4,$ proving \eqref{3.76i}.
\end{example}

Further explicit examples of Krein extensions can be
found in \cite{AS80} and the references therein. All
self-adjoint
extensions of $\dot H$ are described in \cite{AGHKH88},
Section~I.1.1 and Ch.1.5. Generalized Friedrichs and
Krein extensions
in the case where $\dot H$ has deficiency indices $(1,1)$
and $\dot H$ is not necessarily assumed to be bounded
from below,
are studied in detail in
\cite{HS97a}--\cite{HKS97b}.~Interesting
inverse spectral problems associated with self-adjoint
extensions of symmetric operators with gaps were studied
in the series of papers \cite{ABN96},
\cite{BN94}--\cite{BNW93}.

Finally we discuss some realization theorems
for Herglotz
operators of the form \eqref{3.76}. For this
purpose introduce the following set of Herglotz
operators,
\begin{align}
&\calN_0(\calN)=\{M\in\calB(\calN)
\text{ Herglotz}\,
\vert \, M(z)=\smallint_\bbR d \Omega (\lambda)
((\lambda-z)^{-1}-\lambda(1+\lambda^2)^{-1});
\no \\
&\hspace*{1.7cm} \wti \Omega (\bbR)=I_\calN;
\text{ for all }
 \xi\in\calN\backslash \{0\}, \,
\smallint_\bbR d(\xi, \Omega (\lambda)\xi)_\calN=\infty \},
\lb{3.77} \\
&\calN_{0,F }(\calN)=\{M\in \calN_0(\calN)\,\vert \,
\text{supp} (\Omega)\subseteq [0,\infty);
\text{ for all } \xi \in \calN\backslash \{ 0 \},
\no \\
& \hspace*{3.5cm} \smallint_R^\infty d(\xi,
\Omega(\lambda) \xi)_\calN\lambda^{-1}=\infty
\text{ for some }  R>0 \},
\lb{3.78} \\
&\calN_{0, K}(\calN)=\{M\in \calN_0(\calN) \,
\vert \,
\text{supp}(\Omega)\subseteq [0,\infty);
\text{ for all } \xi \in \calN\backslash \{0\},
\no \\
& \hspace*{3.6cm} \smallint_0^R d(\xi,
\Omega(\lambda) \xi)_\calN\lambda^{-1}=\infty
\text{ for some }  R>0 \}, \lb{3.79} \\
&\calN_{0,F,K }(\calN)=\{M\in \calN_0(\calN)\,
\vert\,
 \text{supp} (\Omega)\subseteq [0,\infty);
\text{ for all } \xi \in \calN\backslash \{0\},
\no \\
& \hspace*{2.1cm} \smallint_R^\infty d(\xi,
\Omega(\lambda) \xi)_\calN\lambda^{-1}=
\smallint_0^R d(\xi,
\Omega(\lambda) \xi)_\calN\lambda^{-1}
=\infty \text{ for some }  R>0 \} \no \\
& \hspace*{1.6cm} =\calN_{0,F}(\calN)
\cap\calN_{0,K}(\calN),
\lb{3.80}
\end{align}
where $\calN$ is a separable complex Hilbert space,
$\text{supp} (\Omega)$ denotes the topological
support of $\Omega$, and $\wti \Omega(\lambda)
=(1+\lambda^2)^{-1}\Omega(\lambda),$ $\lambda\in\bbR$.
\begin{theorem}\label{t3.9}

(i) Any $M\in \calN_0(\calN)$ can be realized in
the form
\begin{equation}
M(z)=V^*(zI_{\calN_+} + (1+z^2)P_{\calN_+}
(H-z)^{-1} P_{\calN_+}\big\vert_{\calN_+} ) V, \quad
z\in \bbC\backslash \bbR,
\lb{3.81}
\end{equation}
where $H$ denotes a self-adjoint extension of some
densely defined closed symmetric
operator $\dot H$ with deficiency subspaces
$\calN_\pm$ in some separable Hilbert
space $\calH$.

\noindent (ii)
Any $M\in \calN_{0, F (\text{resp.} K)} (\calN)$ can
be realized in the form
\begin{equation}
M(z)=V^*(zI_{\calN_+} + (1+z^2)P_{\calN_+}
(H_{F (\text{resp.} K)}-z)^{-1} P_{\calN_+}
\big\vert_{\calN_+} ) V, \quad
z\in \bbC\backslash \bbR,
\lb{3.82}
\end{equation}
where $H_{F (\text{resp.} K)} \ge 0$ denotes the
Friedrichs (respectively, Krein) extension
 of some densely defined closed symmetric operator
$\dot H$ with
 deficiency subspaces $\calN_{\pm}$ in some
separable complex Hilbert space $\calH$.

\noindent (iii)
Any $M\in \calN_{0, F , K} (\calN)$ can be realized
in the form
\begin{equation}
M(z)=V^*(zI_{\calN_+} + (1+z^2)P_{\calN_+}
(H_{F, K}-z)^{-1} P_{\calN_+}\big\vert_{\calN_+} ) V,
\quad
z\in \bbC\backslash \bbR,
\lb{3.83}
\end{equation}
where $H_{F,K}\ge 0$ denotes  the unique nonnegative
self-adjoint
 extension of some densely defined closed symmetric
operator $\dot H$
  with deficiency
subspaces $\calN_\pm$ in some separable complex
Hilbert space $\calH$.

In all cases (i)--(iii), $V$ denotes a unitary
operator from $\calN$  to $\calN_+$.
\end{theorem}
\begin{proof}
(i) Define
\begin{equation}
V:\calN\to \hatt{\calN}_+, \quad \xi \longrightarrow
(\cdot -i)^{-1}\xi
\lb{3.84}
\end{equation}
and use the notation developed for the model pair
$(\hatt{\dot H}, \hatt H)$ in
\eqref{3.26}--\eqref{3.30}, Theorem \ref{t3.4},
and Theorem \ref{t3.6}. Then
\begin{equation}
(V\xi,V\eta)_{{\hat{\calN}}_+}
= \int_\bbR d(\xi,\Omega (\lambda)\eta)_{\calN}
(1+\lambda^2)^{-1}
=(\xi,\eta)_{\calN}, \quad \xi, \eta \in \calN \lb{3.85}
\end{equation}
shows that $V$ is a linear isometry from $\calN$
into $\hatt{\calH}_+$,
\begin{equation}
V^*V=I_\calN, \quad \ran(V^*)=\calN.
\lb{3.86}
\end{equation}
By \eqref{3.75} (identifying $\calN_+$ and
$\calN$),
\begin{equation}
V^{-1}:\hatt{\calN}_+\to \calN, \quad
(\cdot -i)^{-1}\xi \longrightarrow \xi
\lb{3.87}
\end{equation}
is also a linear isometry from $\hatt{\calN}_+$
into $\calN$, implying
\begin{equation}
VV^*=I_{\hat{\calN}_+}, \quad\ran(V)=\hatt{\calN}_+.
\lb{3.88}
\end{equation}
Thus $V$ is unitary and one computes
\begin{align}
&(\xi,V^*(zI_{\hat{\calN}_+}+(1+z^2)
P_{\hat{\calN}_+} (\hatt H-z)^{-1}
P_{\hat{\calN}_+}
\big\vert_{\hat{\calN}_+})V\eta)_\calN \no \\
&=(V\xi,(zI_{\hat{\calN}_+}+(1+z^2)
P_{\hat{\calN}_+}
(\hatt H-z)^{-1}P_{\hat{\calN}_+}
\big\vert_{\hat{\calN}_+})V\eta)_{\hat{\calN}_+}
\no \\
&=((\cdot -i)^{-1}\xi,(zI_{\hat{\calN}_+}
+(1+z^2)P_{\hat{\calN}_+} (\hatt H-z)^{-1}
P_{\hat{\calN}_+}
\big\vert_{\hat{\calN}_+})
(\cdot -i)^{-1}\eta)_{\hat{\calN}_+} \no \\
&=\int_\bbR d(\xi, \Omega (\lambda)\eta)_{\calN}
z (1+\lambda^2)^{-1}+
\int_\bbR d(\xi, \Omega (\lambda)\eta)_{\calN}
(1+z^2) (1+\lambda^2)^{-1}(\lambda-z)^{-1}
\no \\
&=\int_\bbR d(\xi,\Omega (\lambda)\eta)_{\calN}
((\lambda-z)^{-1}-\lambda(1+\lambda^2)^{-1})
\no \\
&=(\xi, M(z) \eta)_{\calN}, \quad \xi,\eta,\in \calN,
\,z\in\bbC\backslash \bbR.
\lb{3.89}
\end{align}
(ii) and (iii) then follow in the same way
using Theorem \ref{t3.7}.
\end{proof}

For a whole scale of Nevanlinna classes in the case
where $\dot H$ has deficiency indices $(1,1)$ we refer
to \cite{HSW98}.

\begin{remark}\label{r3.10}
In the special case where $\dim_\bbC
(\calN)\in \bbN$, treated in detail in
\cite{GT97},
 we also considered at length the case where
$H$ and $H_F$
 (respectively, $H_K$) were relatively prime
operators with respect
 to $\dot H$.
 In this case the limiting behavior of $M(z)$
as $\lambda\downarrow -\infty$
(respectively, $\lambda \uparrow 0$) crucially
entered the corresponding results in Theorems
 7.5--7.7 of  \cite{GT97}.
 These limits are given in terms of
$\Re ((P_{1,2} (i)\big\vert_{\calN_+})^{-1})$
 (cf. \eqref{3.11}) identifying $H_1=H$,
$H_2=H_F$ or $H_K$ , etc. In the present
infinite-dimensional case,
$(P_{1,2} (i)\big\vert_{\calN_+})^{-1}$
exists if $H_1$ and $H_2$ are relatively prime
with respect to $\dot H$.
However, $(P_{1,2} (i)\big\vert_{\calN_+})^{-1}$
is not necessarily a bounded operator in
$\calN_+.$
In fact,
\begin{align}
&\Im ((P_{1,2} (i)\big\vert_{\calN_+})^{-1})
=-I_{\calN_+},
\lb{3.90} \\
&\Re ((P_{1,2} (i)\big\vert_{\calN_+})^{-1})
 \in \calB (\calN_+) \text{ if and
only if } \ran(P_{1,2} (i)) =\calN_+
\lb{3.91}
\end{align}
as shown in Lemma 2 of  \cite{GMT97}. This
complicates matters since now the
limits of $M(\lambda)$ as $\lambda\downarrow
-\infty$
(or $\lambda \uparrow 0$)
 may exist but possibly represent unbounded
self-adjoint
operators in $\calN_+$ and thus convergence
of $M(\lambda)$ as $\lambda\downarrow -\infty$
(or $\lambda\uparrow 0$) in these
cases is understood in the strong resolvent sense.
A detailed treatment of this topic goes beyond
the scope of this paper and is thus postponed.
\end{remark}
\begin{theorem}\label{t3.11}
Suppose $M_\ell\in \calN_0(\calN)$, $\ell=1,2$ and
$M_1\ne M_2$. Then $M_1$ and $M_2$ can be realized
as
\begin{equation}
M_\ell(z)=V^*(zI_{\calN_+}+(1+z^2)P_{\calN_+}
(H_\ell-z)^{-1}P_{\calN_+}
\big\vert_{\calN_+})V, \quad \ell=1,2 ,\,\, z\in
\bbC\backslash \bbR,
\lb{3.92}
\end{equation}
where $H_\ell$, $\ell=1,2$ are distinct
self-adjoint extensions of one and
 the same densely defined closed symmetric
operator
$\dot H$ with deficiency subspaces $\calN_\pm$ in
some separable complex Hilbert space $\calH$, and
$V$ denotes a
unitary operator from $\calN$ to $\calN_+$, if
and only if,
\begin{equation}
M_2(z)=
e^{-i\alpha}
(\cos(\alpha)+\sin(\alpha)M_1(z))
(\sin(\alpha)-\cos(\alpha)M_1(z))^{-1}e^{i\alpha},
\quad
z\in \bbC\backslash \bbR
\lb{3.93}
\end{equation}
for some self-adjoint operator $\alpha$ in $\calN$.
\end{theorem}
\begin{proof}
Assuming
\eqref{3.92}, \eqref{3.93} is clear from
\eqref{3.25}. Conversely,
assume
\eqref{3.93}. By Theorem \ref{3.9}\,(i), we may
realize $M_1(z)$ as
\begin{equation}
M_1(z)=V^*(zI_{\calN_+}+(1+z^2)P_{\calN_+}
(H_1-z)^{-1}P_{\calN_+}
\big\vert_{\calN_+})V, \quad z\in \bbC\backslash \bbR.
\lb{3.94}
\end{equation}
If $\wti H\ne H_1$ is another self-adjoint
extension of $\dot H$
 we introduce
\begin{equation}
\wti M(z)=V^*(zI_{\calN_+}+(1+z^2)P_{\calN_+}
(\wti H -z)^{-1}P_{\calN_+}
\big\vert_{\calN_+})V, \quad z\in  \bbC\backslash \bbR,
\lb{3.95}
\end{equation}
and infer from Theorem \ref{t3.3},
\begin{equation}
\wti M(z)=
e^{-i\wti \alpha}
(\cos(\wti \alpha)+\sin(\wti \alpha)M_1(z))
(\sin(\wti \alpha)-\cos(\wti \alpha)M_1(z))^{-1}
e^{i\wti \alpha}, \quad
z\in \bbC\backslash \bbR
\lb{3.96}
\end{equation}
for some $\wti \alpha={\wti \alpha}^*$ in $\calN$.

Since $(H_1-z)(H_1\pm i)^{-1}$ are bounded and
boundedly invertible,
$P_{1,2}(z)$ in \eqref{3.4} uniquely  characterizes
all self-adjoint extensions
$H_2 \ne H_1$ of $\dot H$. Moreover, by
\eqref{3.4a}--\eqref{3.4c}
and von Neumann's representation of self-adjoint
 extensions in terms of Cayley transforms,
all self-adjoint extensions $H_2\ne H_1$ of
$\dot H$ are in a bijective correspondence
to all self-adjoint (possibly unbounded) operators
 $\alpha_{1,2} $
($\alpha_{1,2}\ne \pi /2$) in $\calN_+$.
 Hence we may choose $\wti H$ such that
$\wti\alpha$ equals  $\alpha $
in \eqref{3.93} implying $\wti M(z)=M_2(z).$
\end{proof}

We conclude with a result on analytic continuations
of general Herglotz operators from $\bbC_+$ into a
subset of $\bbC_-$ through an interval of the real
line, which is independent of our emphasis of
perturbation problems in Section~\ref{s2} and
self-adjoint extensions in the present Section~\ref{s3}.
As is well-known, the usual convention for
$M\big|_{\bbC_-}$ by means of reflection as in
\eqref{1.4}, in general, does not
represent the analytic continuation of $M\big|_{\bbC_+}$.
The following result is an adaptation of a theorem of
Greenstein \cite{Gr60} for scalar Herglotz functions to
the present operator-valued context.

\begin{lemma} \lb{l3.12}
Let $\calK$ be a separable
complex Hilbert space and $M$ be a Herglotz operator in
$\calK$ with representation
\eqref{1.1}--\eqref{1.3}. Suppose that the operator
Stieltjes integral
in \eqref{1.1} converges in the strong operator
topology of $\calK$  and let $(\lambda_1,
\lambda_2)\subseteq\bbR$, $\lambda_1 <\lambda_2$. Then
a necessary condition for
$M$ to have an analytic continuation from
$\bbC_+$ into a subset of $\bbC_-$ through the interval
$(\lambda_1,\lambda_2)$ is that for all
$\xi\in\calK$, the associated scalar measures $\omega_\xi=
(\xi,\Omega \,\xi)_\calK$ are purely
absolutely continuous on
$(\lambda_1,\lambda_2)$,
$\omega_\xi \big|_{(\lambda_1,\lambda_2)}
=\big(\omega_\xi \big|_{(\lambda_1,\lambda_2)}\big)_{ac}$,
and
the corresponding density $\omega'_\xi \geq 0$ of
$\omega_\xi$
is real-analytic on $(\lambda_1,\lambda_2)$. If $\calK$ is
finite-dimensional, this condition is also sufficient. If
$M$ has such an analytic continuation into some domain
$\calD_-\subseteq\bbC_-$, then it is given by
\begin{equation}
M(z)=M(\overline z)^* + 2\pi i \Omega'(z),
\quad z\in\calD_-, \lb{3.97}
\end{equation}
where $\Omega'(z)$ denotes the complex-analytic
extension of
$\Omega'(\lambda)$ for $\lambda\in(\lambda_1,\lambda_2)$.
In
particular, $M$ can be analytically continued through
$(\lambda_1,\lambda_2)$ by reflection, that is,
$M(z)=M(\overline z)^*$ for all $z\in\bbC_-$
if $\Omega$ has no support in $(\lambda_1,\lambda_2)$.
\end{lemma}
\begin{proof}
Suppose $M$ has an analytic continuation from
$\bbC_+$ into a subset of $\bbC_-$ through the interval
$(\lambda_1,\lambda_2)$. Then for all $\xi\in\calK$,
Greenstein's result \cite{Gr60}
applies to the scalar Herglotz function $m_\xi (z)=
(\xi,M(z)\xi)_\calK$, $\xi\in\calK$
associated to the measure $\omega_\xi=
(\xi,\Omega \,\xi)_\calK$. Consequently, $m_\xi$ has
an analytic
continuation from $\bbC_+$ into a subset of $\bbC_-$
through
the interval $(\lambda_1,\lambda_2)$ if and only if the
associated scalar measure
$\omega_\xi=(\xi,\Omega \,\xi)_\calK$
is purely absolutely continuous on
$(\lambda_1,\lambda_2)$,
$\omega_\xi \big|_{(\lambda_1,\lambda_2)}
=\big(\omega_\xi \big|_{(\lambda_1,\lambda_2)}\big)_{ac}$,
and
the corresponding density $\omega'_\xi \geq 0$ of
$\omega_\xi$
is real-analytic on $(\lambda_1,\lambda_2)$. In
this case
the analytic continuation of
$m_\xi$ into some domain $\calD_{-,\xi}\subseteq\bbC_-$
is given by
\begin{equation}
m_\xi(z)=m_\xi(\overline z)^* + 2\pi i \omega'_\xi(z),
\quad z\in\calD_{-,\xi}, \lb{3.98}
\end{equation}
where $\omega'_\xi(z)$ denotes the complex-analytic
extension of
$\omega'_\xi(\lambda)$ for $\lambda\in
(\lambda_1,\lambda_2)$. This can be seen as follows: If
$m_x$ can be analytically continued through $(\lambda_1,
\lambda_2)$ into some region $\calD_-\subseteq\bbC_-$,
then $\wti m_\xi(z):=m_\xi(z)-\pi i \omega'_\xi(z)$ is
real-analytic on
$(\lambda_1,\lambda_2)$ and hence can be continued
through
$(\lambda_1,\lambda_2)$ by reflection. Similarly,
$\omega'_\xi(z)$,
being real-analytic, can be continued through
$(\lambda_1,
\lambda_2)$
by reflection. Hence \eqref{3.98} follows from
\begin{equation}
m_\xi(z)-\pi i \omega'_\xi(z)=\wti m_\xi(z)
=\overline{{\wti m_\xi}(\overline z)}
=\overline{m_\xi(\overline z)} + \pi i \omega'_\xi(z),
\quad z\in\calD_-. \lb{3.99}
\end{equation}
Applying a standard polarization
argument, we obtain that the analytic continuation of
$m_{\xi,\eta}(z)=(\xi,M(z)\eta)_\calK$,
$\xi,\eta\in\calK$ into
some domain $\calD_{-,\xi,\eta}\subseteq\bbC_-$ is
given by
\begin{equation}
m_{\xi,\eta}(z)=m_{\xi,\eta}(\overline z)^*
+ 2\pi i \omega'_{\xi,\eta}(z),
\quad z\in\calD_{-,\xi,\eta},
\lb{3.100}
\end{equation}
where $\omega'_{\xi,\eta}(z)=(\xi,
\Omega'(z)\,\eta)_\calK$ is
related to $\omega'_{\xi\pm \eta}(z)$ and
$\omega'_{\xi\pm i\eta}(z)$ by
polarization. In particular, if $M(z)$ has such
an analytic
continuation through the interval
$(\lambda_1,\lambda_2)$
it is necessarily of the form stated in
\eqref{3.97}. If
$\dim_\bbC(\calK)<\infty$, then \eqref{3.98}
and \eqref{3.99}
yield the weak and hence $\calB(\calK)$-analytic
continuation
of $M$ through the interval $(\lambda_1,\lambda_2)$.
\end{proof}
Formula \eqref{3.97} shows that any possible singularity
behavior
of $M\big|_{\bbC_-}$ is determined by that of
$\Omega'\big|_{\bbC_-}$ since $M$, being Herglotz, has
no singularities in $\bbC_+$. Moreover,
analytic continuations through different intervals on
$\bbR$ in general, will lead to different
$\Omega'(z)$ and
hence to branch cuts of $M\big|_{\bbC_-}$.

\section{One-Dimensional Applications} \lb{s4}

In our final section we consider concrete
applications of the
formalism of Section~\ref{s3} in the special case
$\dim_{\bbC}(\calN_+)=1$. We study
Schr\"odinger
operators on a  half-line, compare the corresponding
Donoghue and Weyl-Titchmarsh $m$-functions, and
prove some estimates on linear functionals
associated with these Schr\"odinger operators. We
conclude with two illustrations of Livsic's
result \cite{Li60} on quasi-hermitian
extensions in the special case of densely
defined closed prime symmetric operators with
deficiency indices $(1,1)$ in connection with
first-order differential expressions $-id/dx$.

First we specialize some of the abstract
material in Section~\ref{s3} to the case of a
 densely defined closed prime symmetric
operator $\dot H$ in
 a separable complex Hilbert
space $\calH$ with deficiency indices $(1,1)$.
This case has
been studied in detail by Donoghue  \cite{Do65}
 (see also  \cite{GT97}) and we partly follow
his analysis.

Choose $u_\pm \in \ker ({\dot H}^*\mp i)$ with
$\vert\vert u_\pm \vert\vert_{\calH}=1$ and
introduce the one-parameter
 family $H_\alpha$, $\alpha\in [0,\pi )$ of
self-adjoint
extensions $\dot H$ in $\calH$ by
\begin{align}
&H_\alpha (f+c(u_++e^{2i\alpha} u_-))
=\dot H f +c(iu_+ - i e^{2i\alpha} u_-), \no \\
&\dom (H_\alpha)=\{
 (f+c(u_++e^{2i\alpha} u_-)) \in
\dom ({\dot H}^*)\,\vert\, f \in
 \dom ({\dot H}), \, c\in\bbC \},\quad
\alpha\in [0,\pi ).
\lb{4.1}
\end{align}
Let $\{E_{H_\alpha}(\lambda) \}
_{\lambda\in \bbR}$ be the family
of orthogonal spectral projections of
$H_\alpha$ and suppose that $H_\alpha$ has
simple spectrum for one
(and hence for all) $\alpha\in [0, \pi).$
(This is equivalent
 to the assumption that $\dot H$ is a prime
symmetric operator
 and also equivalent to the fact that
$u_+$ is a cyclic vector for $H_\alpha$ for
all $\alpha\in [0, \pi )$.)
Next we introduce the model representation
$(\hatt{\dot H}_\alpha, \hatt H_\alpha)$ for
$(\dot H, H_\alpha)$ discussed in
\eqref{3.26}--\eqref{3.30}, Theorem \ref{t3.4},
and Theorem \ref{t3.6}. However, since in the
present context $\calN_+$ is a one-dimensional
subspace of $\calH$,
\begin{equation}
\calN_+=\text{lin.span}\{u_+\}, \lb{4.1a}
\end{equation}
the model Hilbert space
$\hatt {\calH}_\alpha=L^2(\bbR,\calN_+;
d\Omega_{H_\alpha,\calN_+})$, $\alpha\in [0,\pi)$
with the operator (in fact, rank-one) valued
measure $\Omega_{H_\alpha,\calN_+}$,
\begin{align}
&\Omega_{H_\alpha,\calN_+}(\lambda)
=\omega_\alpha(\lambda)P_{\calN_+}|_{\calN_+},
\quad P_{\calN_+}=(u_+,\cdot)u_+, \lb{4.1b} \\
&\omega_\alpha(\lambda)=(1+\lambda^2)
\|E_{H_\alpha}(\lambda)u_+\|^2_\calH, \quad
\alpha\in [0,\pi), \no
\end{align}
can be replaced by the model space
$\wti{\calH}_\alpha =L^2(\bbR;d\omega_\alpha)$
with scalar measure $\omega_\alpha$. In
particular, $\omega_\alpha(\lambda)$ can be taken as
the control measure in this special case and
\begin{align}
&V:\hatt {\calH}_\alpha=L^2(\bbR,\calN_+;
d\Omega_{H_\alpha,\calN_+})\to
\wti{\calH}_\alpha =L^2(\bbR;d\omega_\alpha)
\no \\
&\hat f=\{\hat f (\lambda)=
\ti f(\lambda)u_+\}_{\lambda\in\bbR}
\to V\hat f=\ti f=\{\ti f (\lambda)\}
_{\lambda\in\bbR} \lb{4.1d}
\end{align}
represents the corresponding unitary operator from
$\hatt {\calH}_\alpha=L^2(\bbR,\calN_+;
d\Omega_{H_\alpha,\calN_+})$ to
$\wti{\calH}_\alpha =L^2(\bbR;d\omega_\alpha)$.~Hence
we translate in the following some of the
results of Theorems~\ref{t3.4} and \ref{t3.6}
from $\hatt {\calH}_\alpha$ to
$\wti {\calH}_\alpha$. However, due to the trivial
nature of the unitary operator $V$ in \eqref{4.1d}, we
will ignore this additional isomorphism and
simply keep using our $\hatt{\,\,\,\,}$-notation of
Section~\ref{s3} instead of the new
$\wti{\,\,\,\,}$-notation.
Thus, we consider the model Hilbert space
$\hatt{\calH}_\alpha=L^2(\bbR; d \omega_\alpha)$,
$\alpha \in [0, \pi )$, where
\begin{align}
&\omega_\alpha (\lambda)=
(1+\lambda^2)
 \vert\vert E_{H_\alpha}(\lambda)
 u_+ \vert\vert_{\calH}^2, \quad\alpha
\in [0, \pi ), \lb{4.2} \\
&\int_\bbR d \omega_\alpha (\lambda)
(1+\lambda^2)^{-1}=1,
\quad \int_\bbR d \omega_\alpha(\lambda)
=\infty, \quad \alpha \in [0, \pi )
\lb{4.3}
\end{align}
and define in $\hatt{\calH}_\alpha$ the
self-adjoint operator $\hatt H_\alpha$,
\begin{equation}
(\hatt H_\alpha
\hat f)(\lambda)=\lambda \hat f(\lambda),
\quad
\hat f \in \dom (\hatt H_\alpha)=
L^2 (\bbR ; (1+\lambda^2)d \omega_\alpha )
\lb{4.4}
\end{equation}
and its densely defined and closed restriction
$\hatt{\dot H}_\alpha$,
\begin{equation}
\dom (\hatt{\dot H}_\alpha)=\{ \hat f
\in \dom (\hatt H_\alpha)
\,\vert
\smallint_\bbR d \omega_\alpha (\lambda )
\hat f (\lambda)=0 \}, \quad
\hatt{\dot H}_\alpha=\hatt H_\alpha
\big\vert_{\dom (\hat{\dot H}_\alpha)}. \lb{4.5}
\end{equation}
Then
\begin{equation}
\ker ({\hatt{\dot H}}^*-z)=\{
c(\cdot -z)^{-1} \in \hatt{\calH}_\alpha
\,\vert\, c \in \bbC \}
\lb{4.6}
\end{equation}
and the pair $(\dot H, H_\alpha)$ in $\calH$
is unitarily equivalent to the pair
 $(\hatt{\dot H}_\alpha, \hatt H_\alpha)$ in $
\hatt{\calH}_\alpha$ (cf. Theorem \ref{t3.6}).
This representation of $(\dot H, H_\alpha)$
in terms of  $(\hatt{\dot H}_\alpha,
\hatt H_\alpha)$
has the advantage of very simple definitions
of $\hatt H_\alpha$ and $\hatt{\dot H}_\alpha$,
however,
one has to pay a price since different
$H_\alpha, \hatt{\dot H}_\alpha$ act in
different Hilbert spaces
$\hatt{\calH}_\alpha$. Hence it is desirable
to determine the
expression for all
$H_\alpha$, $\alpha \in [0, \pi )$
in connection with one fixed $\alpha$ say,
 $\alpha_0\in [0, \pi )$, in the corresponding
fixed Hilbert space
$\hatt{\calH}_{\alpha_0}=L^2 (\bbR;
d \omega_{\alpha_0})$
and we turn our attention to this task next.
\begin{lemma}\label{l4.1}
Fix $\alpha_0 \in [0, \pi )$ and define
\begin{equation}
U_{\alpha_0}:\hatt{\calH}_{\alpha_0}
\longrightarrow \calH, \quad
\hat f\to U_{\alpha_0} \hat f=
\slim_{N\to \infty} \int_{-N}^N
d (E_{H_{\alpha_0}}(\lambda)u_+) (\lambda - i )
\hat f (\lambda). \lb{4.7}
\end{equation}
Then  $U_{\alpha_0}$ is a unitary operator
 from ${\hatt{\calH}}_{\alpha_0} $ to
$\calH$ and
\begin{equation}
\dot H= U_{\alpha_0} \hatt{\dot H}
_{\alpha_0}U_{\alpha_0}^{-1},
\quad
H_{\alpha_0} = U_{\alpha_0} {\hatt H}
_{\alpha_0} U_{\alpha_0}^{-1}.
\lb{4.8}
\end{equation}
Moreover,
\begin{align}
&\hat u_+(\lambda)= (U^{-1}_{\alpha_0} u_+)
(\lambda)=(\lambda-i)^{-1}, \lb{4.9} \\
&\hat u_-(\lambda)= (U^{-1}_{\alpha_0} u_-)
(\lambda)=-e^{-2i\alpha_0}(\lambda+i)^{-1},
\quad \lambda\in \bbR, \lb{4.9a}
\end{align}
and hence
\begin{align}
&(U^{-1}_{\alpha_0}( u_++e^{2i\alpha}u_-))
(\lambda)=
2ie^{i(\alpha-\alpha_0 )}
(1+\lambda^2)^{-1}(
-\lambda \sin (\alpha-\alpha_0)+
\cos (\alpha-\alpha_0)), \quad
\no \\
& \hspace*{8.5cm}\alpha \in [0, \pi ),
\, \lambda \in \bbR. \lb{4.10}
\end{align}
\end{lemma}
\begin{proof}
\eqref{4.7} and \eqref{4.8} have been discussed
in Theorem \ref{t3.6}, \eqref{4.9} is clear from
\eqref{4.7}. From
\begin{align}
&U^{-1}_{\alpha_0}H_{\alpha_0}
( u_++e^{2i\alpha_0}u_-)=
U^{-1}_{\alpha_0}{\dot H}^*
( u_++e^{2i\alpha_0}u_-)=
i \hat u_+ - i e^{2i\alpha_0}\hat u_-,
\lb{4.11} \\
&\hatt H_{\alpha_0}
( \hat u_++e^{2i\alpha_0} \hat u_-)=\lambda
(\hat u_++e^{2i\alpha_0}\hat u_-),
\lb{4.12}
\end{align}
and \eqref{4.9} one infers
\begin{equation}
i(\lambda-i)^{-1}
-ie^{2i\alpha_0}
 \hat u_- (\lambda)=\lambda (\lambda-i)^{-1}
+e^{2i\alpha_0}
 \lambda \hat u_-(\lambda)
\lb{4.13}
\end{equation}
and hence \eqref{4.9a}. Equation \eqref{4.10} then
immediately follows from \eqref{4.9} and
\eqref{4.9a}.
\end{proof}

Equation  \eqref{4.10} confirms the fact that
any two different self-adjoint extensions
 of $\dot H$ are relatively prime
\begin{equation}
\dom (H_\alpha)\cap \dom (H_\beta) =\dom (\dot H),
\quad \alpha, \beta \in [0, \pi), \quad
\alpha\ne \beta
\lb{4.14}
\end{equation}
since
$\int_\bbR d \omega_{\alpha_0} (\lambda)=\infty$
 and hence
\begin{equation}
\int_\bbR d \omega_{\alpha_0}(\lambda) \lambda^2
\vert U_{\alpha_0}^{-1}(u_++e^{2i\alpha} u_-)
(\lambda)\vert^2=\infty \text{ for all }
 \alpha\ne \alpha_0.
\lb{4.15}
\end{equation}
This is of course an artifact of our
 special hypothesis $\text{def} (\hatt H)=(1,1).$

\vspace*{2mm}

Next, consider the normalized element
(cf. \eqref{4.10}
for $\alpha=\alpha_0$)
\begin{align}
&\hat g_\alpha \in \big(\ker (\hatt{\dot H}^*-i)\dot
+ \ker
(\hatt{\dot H}^*+i)\big)\cap \dom (\hatt H_\alpha),
\no \\
&{\hat g}_\alpha(\lambda)=\bigg(\int_\bbR
d \omega_\alpha (\nu)
 (1+\nu^2)^{-2}\bigg)^{-1/2}
(1+\lambda^2)^{-1}, \quad
\vert \vert  \hat g_\alpha \vert
\vert_{\hat{\calH}_\alpha}
=1. \lb{4.16}
\end{align}
Then
\begin{equation}
\dom (\hatt H_\alpha)=\text{lin.span}
\{ {\hat g}_\alpha \} \dot +
\dom (\hatt{\dot H}_\alpha)
\lb{4.17}
\end{equation}
by von Neumann's theory of self-adjoint extensions
of symmetric operators (cf., e.g., \cite{AG93}, Ch.~VII,
\cite{DS88}, Sect.~II.4, \cite{Na68}, Sect.~14,
\cite{RS75}, Sect.~X.1, \cite{Ne29}) and we may
consider the linear functional
$\ell_{\hat g_\alpha}$ on
$\dom (\hatt H_\alpha)$ defined by
\begin{equation}
\ell_{{\hat g}_\alpha}:
\dom ({\hatt H}_\alpha) \to \bbC, \quad
 \ell_{{\hat g}_\alpha}(\hat f)=c,
\lb{4.18}
\end{equation}
where
\begin{equation}
\hat f \in \dom(\hatt H_\alpha),
\quad \hat f =c {\hat g}_\alpha +\hat h,
\quad \hat h \in \dom (\hatt{\dot H}_\alpha).
\lb{4.19}
\end{equation}
\begin{lemma}\label{l4.2}
Let $\alpha \in [0, \pi  )$. Then
\begin{equation}
\sup_{\hat f \in \dom (\hat H_\alpha)}
\bigg(\frac{\vert \ell_{\hat g_\alpha}
(\hat f)\vert^2 }
{\vert\vert \hat f\vert\vert^2
_{\hat{\calH}_\alpha}+
\vert\vert \hatt H_\alpha \hat f
\vert\vert^2_{\hat{\calH}_\alpha}
}\bigg)=\int_\bbR d \omega_\alpha
(\lambda)(1+\lambda^2)^{-2}.
\lb{4.20}
\end{equation}
\end{lemma}
\begin{proof}
By \eqref{4.3} and \eqref{4.5} one computes
\begin{equation}
\int_\bbR d \omega_\alpha (\lambda)
\hat f (\lambda)=
c\int_\bbR d \omega_\alpha (\lambda)
\hat g_\alpha (\lambda)=
\ell_{{\hat g}_{\alpha}} (\hat f)
\bigg(\int_\bbR d \omega_{\alpha} (\lambda )
(1+\lambda^2)^{-2}\bigg)^{-1/2}
\lb{4.21}
\end{equation}
and hence the Cauchy-Schwarz inequality
applied to
\begin{align}
&\bigg| \int_\bbR d \omega_\alpha (\lambda)
\hat f(\lambda)
\bigg| \le
\bigg(\int_\bbR d \omega_\alpha (\lambda)
(1+\lambda^2)
\vert \hat f(\lambda) \vert^2\bigg)^{1/2}
\bigg(\int_\bbR d \omega_\alpha (\lambda)
(1+\lambda^2)^{-1}\bigg)^{1/2} \no \\
&=(\vert\vert\hat f
\vert\vert_{\hat{\calH}_\alpha}^2
+\vert\vert
{\hatt H}_\alpha  \hat f
\vert\vert_{\hat{\calH}_\alpha}^2)^{1/2}
\lb{4.22}
\end{align}
yields
\begin{equation}
\frac{\vert \ell_{{\hat g}_\alpha } (\hat f)\vert^2}
{\vert\vert\hat f
\vert\vert_{\hat{\calH}_\alpha}^2+
\vert\vert{\hatt H}_\alpha  \hat f
\vert\vert_{\hat{\calH}_\alpha}^2}
\le \int_\bbR d \omega_\alpha
(\lambda)(1+\lambda^2)^{-2}.
\lb{4.23}
\end{equation}
Since inequality \eqref{4.23} saturates
for ${\hat f}_0(\lambda)=(1+\lambda^2)^{-1}$,
${\hat f}_0 \in \dom ({\hatt H}_\alpha)$,
\eqref{4.20} is proved.
\end{proof}
Introducing the Donoghue-type $m$-function
\begin{equation}
m_\alpha^D (z)=\int_\bbR d \omega_\alpha (\lambda)
((\lambda - z )^{-1} -\lambda (1+\lambda^2)^{-1}),
 \quad \alpha \in [0, \pi ), \,\, z\in \bbC_+,
\lb{4.24}
\end{equation}
the analog of \eqref{3.13}, one can prove the
following result.
\begin{lemma}\label{l4.3} \mbox{\rm (Donoghue
\cite{Do65}.)}
\begin{equation}
m_\beta^D(z)=
\frac{
-\sin(\beta-\alpha)+\cos(\beta-\alpha)m_\alpha^D(z)}
{\cos(\beta-\alpha)+\sin(\beta-\alpha)m_\alpha^D(z)
}, \quad \alpha, \beta \in [0, \pi), \,\, z
\in \bbC_+. \lb{4.25}
\end{equation}
\end{lemma}

Next we turn to Schr\"odinger operator on the
half-line $[0, \infty )$.
 Let $q\in L^1([0, R])$ for all $R>0$, $q$
real-valued and
introduce the fundamental system
$\phi_\gamma(z,x)$, $\theta_\gamma(z,x)$, $z \in
\bbC$ of solutions of
\begin{equation}
-\psi''(z,x) +(q(x)-z)\psi(z,x)=0, \quad x>0
\lb{4.26}
\end{equation}
($\,'$ denotes $d/dx$) satisfying
\begin{equation}
\phi_\gamma(z, 0_+)=-
\theta'_\gamma(z, 0_+)=-\sin(\gamma), \quad
\phi'_\gamma(z, 0_+)=
\theta_\gamma(z, 0_+)=\cos(\gamma),
 \quad \gamma\in [0, \pi ).
\lb{4.27}
\end{equation}
Assuming that $-\frac{d^2}{dx^2}+q$ is in the limit
 point case at $\infty$, let $\psi_\gamma(z,x)$ be
the unique solution of \eqref{4.26} satisfying
\begin{align}
&\psi_\gamma(z, \cdot) \in L^2([0, \infty );
dx), \quad
\sin(\gamma) \psi_\gamma '(z,0_+)+
\cos(\gamma) \psi_\gamma (z,0_+)=1, \lb{4.28} \\
& \hspace*{7.2cm} \gamma\in [0, \pi ),
\,\, z \in \bbC_+. \no
\end{align}
Then $\psi_\gamma(z,  x)$ is of the form (see,
e.g., the discussion of Weyl's theory
in Appendix~A of  \cite{GS96})
\begin{equation}
\psi_\gamma(z,x)=\theta_\gamma(z,x)
+m_\gamma^W(z) \phi_\gamma(z,x),
\quad\gamma\in [0, \pi), \,\, z\in \bbC_+,
\lb{4.29}
\end{equation}
where $m_\gamma^W (z)$ denotes the Weyl-Titchmarsh
$m$-function \cite{Ti62}, Chs.~II, III, \cite{We10}
(as opposed
to Donoghue's $m$-function $m_\alpha^D(z)$ in
\eqref{4.24}) corresponding to the operator
$\wti H_\gamma$
 in $L^2([0, \infty); dx)$ defined by
\begin{align}
&(\wti H_\gamma f)(x)=-f''(x) +q(x)f(x),
\quad x>0, \no \\
&f \in \dom (\wti H_\gamma )=
\{ g \in L^2([0, \infty ); dx)\,\vert
\, g, g' \in AC([0, R])
\text{ for all } R>0; \lb{4.30} \\
& \hspace*{1.35cm} \sin (\gamma) g'(0_+)+
\cos (\gamma) g(0_+)=0; \,
-g''+qg\in L^2([0, \infty);dx) \},
\quad \gamma\in [0, \pi). \no
\end{align}
The family $\wti H_\gamma$, $\gamma\in
[0, \pi)$ represents
all self-adjoint extensions of the
densely defined closed prime symmetric operator
$\dot{\wti H}$ in
$L^2([0, \infty); dx)$ of deficiency indices
$(1,1)$,
\begin{align}
&(\dot{\wti H}f)(x)=-f''(x)+q(x)f(x),
\quad x>0, \no \\
&f \in \dom (\dot{\wti H}_\gamma )=
\{ g \in L^2([0, \infty ); dx)) \,\vert \,g,
g' \in AC([0, R])\text{ for all } R>0;
\lb{4.30a} \\
& \hspace*{3.4cm} g'(0_+)=g(0_+)=0; \,
-g''+qg\in L^2([0, \infty);dx) \}. \no
\end{align}
(Here $AC([a,b])$ denotes the set of absolutely
continuous functions on $[a,b]$.)
 Weyl's $m$-function is a Herglotz function
with representation
\begin{align}
&  m_\gamma^W(z)=
\begin{cases}
c_\gamma+\int_\bbR
d\omega_\gamma^W(\lambda)((\lambda-z)^{-1}
-\lambda(1+\lambda^2)^{-1}), \,
& \gamma \in [0, \pi),   \lb{4.31} \\
\cot (\gamma) +\int_\bbR
d \omega_\gamma^W(\lambda)
(\lambda-z)^{-1}, &\gamma \in (0, \pi ),
\end{cases}
\end{align}
for some $c_\gamma\in\bbR$, where
\begin{align}
&  \int_\bbR d \omega_\gamma^W(\lambda)
(1+\vert \lambda \vert)^{-1}
\begin{cases} < \infty,
& \gamma \in (0,\pi),  \lb{4.32}\\
=\infty, &\gamma= 0. \end{cases}
\end{align}
Moreover, one can prove the following result.
\begin{lemma}\label{l4.4} \mbox{\rm (See, e.g.,
Aronszajn \cite{Ar57}, \cite{EK82}, Sect.~2.5.)}
\begin{equation}
m_\delta^W (z)=\frac{-\sin(\delta-\gamma)
+\cos(\delta-\gamma)m_\gamma^W(z)}
{\cos(\delta-\gamma)+
\sin(\delta-\gamma)m_\gamma^W(z)
}, \quad \delta, \gamma \in [0, \pi),
\,\, z \in \bbC_+.
\lb{4.33}
\end{equation}
Moreover,
\begin{align}
&  m_\gamma^W(z)\underset{z\to i \infty}{=}
\begin{cases}
\cot(\gamma) + O (z^{-1/2}),
& \gamma \in [0,\pi),   \lb{4.34} \\
iz^{1/2} +o(1)
, &\gamma=0.  \end{cases}
\end{align}
\end{lemma}

In the following we denote by
$H_\alpha$
in $L^2([0, \infty);dx)$ the Schr\"odinger
 operator on $[0, \infty)$ defined as in
\eqref{4.1}
 but with $\dot H$ replaced by $\dot{\wti H}$
in \eqref{4.30a}.
The connection between $H_\alpha$ and
${\wti H}_\gamma$ and
$m_\alpha^D(z)$ and $m_\gamma^W(z)$ is then
determined as follows.
\begin{theorem}\lb{t4.5}
Suppose $\gamma(\alpha) \in [0, \pi )$ satisfies
\begin{equation}
\cot (\gamma(\alpha)) =
-\Re (m_0^W(i)) - \Im (m_0^W(i))\tan (\alpha),
\quad \alpha \in [0, \pi ).
\lb{4.35}
\end{equation}
Then
\begin{equation}
H_\alpha =\wti H_{\gamma(\alpha)}, \quad
\alpha\in [0,\pi). \lb{4.36}
\end{equation}
and
\begin{equation}
m_\alpha^D(z)=
(m_{\gamma (\alpha)}^W(z)-
\Re (m_{\gamma(\alpha)}^W (i)) /
\Im (m_{\gamma(\alpha)}^W(i)), \quad \alpha\in [0,\pi),
\, z \in \bbC_+.
\lb{4.37}
\end{equation}
\end{theorem}
\begin{proof}
Since $\psi_\gamma(z, x)$ are just constant
multiples of $\psi_0(z,x)$,
it suffices to focus on
$\psi_0(z,x)$. In order to prove \eqref{4.36},
subject to \eqref{4.35}, we need
\begin{equation}
\eta_\alpha=\vert\vert \psi_0(i)
\vert\vert^{-1}_{L^2([0,
\infty );dx)}\psi_0(i)
+\vert\vert \psi_0(-i)\vert
\vert^{-1}_{L^2([0, \infty );dx)}
e^{2i\alpha}\psi_0(-i)
 \in \dom (H_\alpha )
\lb{4.38}
\end{equation}
according to \eqref{4.1} and the fact
(cf.~\eqref{4.28})
\begin{equation}
u_\pm =\vert\vert \psi_0(\pm i)\vert\vert^{-1}
_{L^2([0, \infty );dx)}
\psi_0(\pm i).
\lb{4.39}
\end{equation}
Since it is known (see, e.g.,~\cite{CL85}, Sect.~9.2,
\cite{EK82}, Sect.~2.2) that
\begin{equation}
\vert\vert \psi_\gamma(z)\vert\vert^2
_{L^2([0, \infty );dx)}=
\Im (m_\gamma^W(z))/\Im (z), \quad z
\in \bbC\backslash \bbR,
\lb{4.40}
\end{equation}
one obtains from \eqref{4.47} and \eqref{4.29}
\begin{equation}
-\cot (\gamma(\alpha))=
\eta_\alpha'(0_+)/\eta_\alpha(0_+)=
(1+e^{2i\alpha})^{-1} (m_0^W(i)
+e^{2i\alpha}m_0^W(-i)),
\lb{4.41}
\end{equation}
which yields \eqref{4.35} and at the same
time proves \eqref{4.36}. By \eqref{4.24} and
\eqref{4.31},
\begin{equation}
m_\alpha^D(z)=A_\alpha m_{\gamma(\alpha)}^W(z)
+B_\alpha, \quad
\alpha \in [0, \pi ), \,\, z \in \bbC_+
\lb{4.42}
\end{equation}
for some $A_\alpha>0$ and $B_\alpha \in \bbR$.
The fact
\begin{equation}
m_\alpha^D (i)=i, \quad \alpha\in [0, \pi )
\lb{4.43}
\end{equation}
(use \eqref{3.5} or combine the normalization
$\int_\bbR d \omega_\alpha(\lambda)
(1+\lambda^2)^{-1}=1$ with \eqref{4.24}) immediately
yields \eqref{4.37}.
\end{proof}
\begin{corollary}\lb{c4.6}
Assume in addition that $\dot H \ge 0$. Then the
Friedrichs extension $H_F$ of $\dot H$ corresponds
to
\begin{equation}
\alpha=\alpha_F=\pi / 2 \quad  \text{and}\quad
\gamma=\gamma_F=0
\lb{4.44}
\end{equation}
and the Krein extension $H_K$ of $\dot H$ corresponds
to
\begin{equation}
\tan (\alpha ) =\tan (\alpha_K) = m_{\pi / 2 }^D(0_-)
\quad \text{and}
\quad\cot (\gamma)=\cot (\gamma_K)=-m_0^W(0_-)
\lb{4.45}
\end{equation}
in \eqref{4.1} and \eqref{4.30}.~The right-hand
sides in \eqref{4.45} are simultaneously infinite
if and only if $H_F=H_K$.
\end{corollary}
\begin{proof}
Since
$\lim_{\lambda \downarrow -\infty}m_0^W(\lambda)
=-\infty$
by \eqref{4.34}, \eqref{4.44} follows from
Corollary \ref{c3.8}\,(i).
 Similarly, \eqref{4.45} follows from \eqref{4.33}
 (replacing $\delta\to \gamma$ and
$\gamma\to 0$) and Corollary \ref{c3.8} (ii).
\end{proof}

Finally we return to the functional
 $\ell_{{\hat g}_\alpha}$ in \eqref{4.18} and
establish its properties in
 connection with the Schr\"odinger
operator $\wti {H}_\gamma$ on $[0, \infty).$
\begin{lemma}\label{l4.7}
Define $\hat g_\alpha$ by
\begin{equation}
U_\alpha^{-1} {\hat g}_\alpha=
\vert\vert\psi_0(i)+e^{2i\alpha}
\psi_0(-i)) \vert\vert^{-1}_{L^2([0,
\infty ); dx)}
(\psi_0(i)+e^{2i\alpha} \psi_0(-i)),
\quad \alpha \in [0, \pi ).
\lb{4.46}
\end{equation}
Then
\begin{align}
&  \ell_{{\hat g}_\alpha} (\hat f) \no \\
&=\begin{cases} (2 i \Im (m_0^W(i)))^{-1}
\vert\vert\psi_0(i)- \psi_0(-i)
\vert\vert_{L^2([0, \infty ); dx)}
 (U^{-1}_{\pi /2} \hat f)' (0_+),
&\alpha=\frac{\pi }{2},   \no \\
(1+e^{2i\alpha})^{-1}
\vert\vert\psi_0(i)+e^{2i\alpha} \psi_0(-i)
\vert\vert_{L^2([0, \infty ); dx)}
 (U^{-1}_{\alpha} \hat f) (0_+),
& \alpha\in [0, \pi )\backslash
\{ \frac{\pi}{2} \},
\end{cases} \no \\
&\hspace*{9.2cm} \hat f\in \dom (\hatt H_\alpha).
\lb{4.47}
\end{align}
\end{lemma}
\begin{proof}
 By \eqref{4.38} and \eqref{4.40},
$$
\psi_0(i)+e^{2i\alpha}\psi_0(-i) \in
\dom (H_\alpha).
$$
Hence
\begin{align}
&f=c\vert\vert\psi_0(i)+e^{2i\alpha}\psi_0(-i)
\vert\vert^{-1}_{L^2([0, \infty);dx)}
(\psi_0(i)+e^{2i\alpha}\psi_0(-i)) +h,
\lb{4.48} \\
&\hspace*{5.5cm} f\in \dom (H_\alpha),
\,\, h\in \dom (\dot H) \no
\end{align}
and
\begin{equation}
\ell_{\hat g_\alpha} (\hat f)=c, \quad \hat f
 \in \dom (\hatt H_\alpha ).
\lb{4.49}
\end{equation}
Since by \eqref{4.30},
\begin{equation}
h'(0_+)=h(0_+)=0,
\lb{4.50}
\end{equation}
one computes in the case $\alpha=\pi/2$
\begin{align}
f'(0_+)&=c\vert\vert\psi_0(i) -\psi_0(-i)
\vert\vert^{-1}_{L^2([0, \infty); dx)}
(\psi_0'(i, 0_+) -\psi_0'(-i, 0_+)) \no \\
&=c\vert\vert
\psi_0(i) -\psi_0(-i)
\vert\vert^{-1}_{L^2([0, \infty); dx)}
2i \Im (m_0^W(i)), \quad
f \in \dom (H_{\pi/ 2})
\lb{4.51}
\end{align}
using \eqref{4.27} and \eqref{4.29}.
 Similarly, for $\alpha\in [0, \pi )
\backslash \{ \pi /2 \}$
 one computes
\begin{align}
f(0_+)&=c\vert\vert
\psi_0(i) +e^{2i\alpha}\psi_0(-i)
\vert\vert^{-1}_{L^2([0, \infty); dx)}
(\psi_0(i, 0_+)
+e^{2i\alpha}\psi_0(-i, 0_+)) \no \\
&=c\vert\vert
\psi_0(i) +e^{2i\alpha}\psi_0(-i)
\vert\vert^{-1}_{L^2([0, \infty); dx)}
(1+e^{2i\alpha}), \lb{4.52} \\
& \hspace*{4cm} f \in \dom (H_{\pi/2}),
\,\, \alpha\in [0, \pi ) \backslash \{\pi/2 \}, \no
\end{align}
since $\psi_0(z, 0_+)=1, z\in \bbC\backslash \bbR$
by \eqref{4.27} and \eqref{4.29}.
Combining \eqref{4.49} and \eqref{4.51},
\eqref{4.52} proves
\eqref{4.47}.
\end{proof}

Lemmas \ref{l4.2} and \ref{l4.3} then yield the
principal result of this section:
\begin{theorem}\label{t4.8}
Let $\alpha\in [0,\pi).$ Then
\begin{align}
&\sup_{f\in \dom (H_{\pi / 2})}
\bigg(\frac{\vert f'(0_+) \vert^2}
{\vert\vert f\vert
\vert^2_{L^2([0, \infty); dx)}+
\vert\vert H_{\pi /2}f
\vert\vert^2_{L^2([0, \infty); dx)}}\bigg)
= \Im (m_0^W(i)),
\lb{4.53} \\
&\sup_{f\in \dom (H_\alpha)}
\bigg(\frac{\vert f(0_+) \vert^2}
{\vert\vert f\vert
\vert^2_{L^2([0, \infty); dx)}+
\vert\vert H_\alpha f\vert
\vert^2_{L^2([0, \infty); dx)}}\bigg)=
\frac{\cos^2(\alpha)}{\Im (m_0^W(i))}. \lb{4.54}
\end{align}
\end{theorem}
\begin{proof}
Consider $\alpha=\pi /2$ first. Then
Lemma \ref{l4.2} combined with
\eqref{4.8}, \eqref{4.39}, and \eqref{4.47}
yields
\begin{align}
&\sup_{f\in \dom (H_{\pi / 2})}
\bigg(\frac{\vert f'(0_+) \vert^2}
{\vert\vert f\vert
\vert^2_{L^2([0, \infty); dx)}+
\vert\vert H_{\pi /2}f\vert
\vert^2_{L^2([0, \infty); dx)}}\bigg) \no \\
&=\frac{
4 \vert \Im (m_0^W(i))\vert^2}
{\vert\vert \psi_0(i)
\vert\vert^2_{L^2([0, \infty); dx)}
\vert\vert u_+ -u_-
\vert\vert^2_{L^2([0, \infty); dx)}}
\int_\bbR d \omega_{\pi / 2} (\lambda)
(1+\lambda^2)^{-2}.
\lb{4.55}
\end{align}
Since
\begin{equation}
\vert\vert
u_+ -u_-
\vert\vert^2_{L^2([0, \infty); dx)}
=
\vert\vert
\hat u_+ -\hat u_-
\vert\vert^2_{\hat {\calH}_{\pi/2}}=
4\int_\bbR d \omega_{\pi / 2} (1 +\lambda^2)^{-2}
\lb{4.56}
\end{equation}
by \eqref{4.9} (taking $\alpha_0=\pi / 2 $) and
\begin{equation}
\vert\vert
\psi_0(i)
\vert\vert^2_{L^2([0, \infty); dx)}=
\Im (m_0^W(i))
\lb{4.57}
\end{equation}
by \eqref{4.40}, the right-hand side of
\eqref{4.55} coincides
with that in \eqref{4.53}. Similarly, one computes
from Lemma \ref{l4.2},
\eqref{4.8}, \eqref{4.39}, and  \eqref{4.47},
\begin{align}
&\sup_{f\in \dom (H_\alpha)}
\bigg(\frac{\vert f(0_+) \vert^2}
{\vert\vert f\vert
\vert^2_{L^2([0, \infty); dx)}+
\vert\vert H_\alpha f\vert
\vert^2_{L^2([0, \infty); dx)}}\bigg) \no \\
&=\frac{4\cos^2 (\alpha)}
{\vert\vert \psi_0(i)
\vert\vert^2_{L^2([0, \infty ); dx)}
\vert\vert u_+ +e^{2i\alpha}u_-
\vert\vert^2_{L^2([0, \infty ); dx)}}
\int_\bbR d \omega_{\alpha} (\lambda)
(1+\lambda^2)^{-2}. \lb{4.58}
\end{align}
Because of \eqref{4.57} and
\begin{equation}
\vert\vert
u_+ +e^{2i\alpha}u_-
\vert\vert^2_{L^2([0, \infty); dx)}=
\vert\vert
\hat u_++ e^{2i\alpha}\hat u_-
\vert\vert^2_{\hat {\calH}_\alpha }
=4 \int_\bbR d \omega_\alpha (\lambda)
(1+\lambda^2)^{-2},
\lb{4.59}
\end{equation}
\eqref{4.58} coincides with \eqref{4.54}.
\end{proof}
\begin{remark}\label{r4.9}

(i) In the special case $q(x)=0$, $x\ge 0$ one has
\begin{equation}
m_0^W(z)=i(z)^{1/2}
\lb{4.60}
\end{equation}
(using the branch with $\Im((z)^{1/2})\geq 0$,
$z\in\bbC$)
and hence \eqref{4.53} yields
\begin{equation}
\vert f'(0_+)\vert \le 2^{-1/4}
\bigg(\int_0^\infty dx (\vert f(x)\vert^2 +\vert
 f''
(x) \vert^2\bigg)^{1/2}, \quad
f \in H_0^{2,2} ((0, \infty)),
\lb{4.61}
\end{equation}
with $2^{-1/4}$ best possible and
\begin{align}
&H_0^{2,2}((0, \infty))=
\{ f \in L^2([0, \infty); dx) \,| \, f,
f'\in AC([0, R])
 \text{ for all } R>0; \, \no \\
& \hspace*{5.2cm} f(0_+)=0;\, f, f''\in
L^2([0, \infty ); dx ) \}
\lb{4.62}
\end{align}
the familiar Sobolev space.

(ii) Multiplying the two results \eqref{4.53}
and \eqref{4.54} reveals the curious fact,
\begin{align}
&\sup_{f\in \dom (H_{\pi / 2})}
\bigg(\frac{\vert f'(0_+) \vert^2}
{\vert\vert f\vert
\vert^2_{L^2([0, \infty); dx)}+
\vert\vert H_{\pi /2}f\vert
\vert^2_{L^2([0, \infty); dx)}}\bigg) \times
\lb{4.63} \\
& \times
\sup_{f\in \dom (H_\alpha)}
\bigg(\frac{\vert f(0_+) \vert^2}
{\vert\vert f\vert
\vert^2_{L^2([0, \infty); dx)}+
\vert\vert H_\alpha f\vert
\vert^2_{L^2([0, \infty); dx)}}\bigg)
=\cos^2(\alpha), \,\,\, \alpha \in
[0, \pi ). \no
\end{align}
\end{remark}

\vspace*{3mm}

Finally, we
conclude with two illustrations of a well-known result
of Livsic \cite{Li60} on quasi-hermitian
extensions in the special case of densely
defined closed prime symmetric operators with
deficiency indices $(1,1)$.

Following Livsic \cite{Li60} one defines a closed
operator
$H$ to be a {\it quasi-hermitian} extension of a densely
defined closed
prime symmetric operator $\dot H$ with deficiency
indices $(1,1)$ if
\begin{equation}
\dot H \varsubsetneqq H \varsubsetneqq
{\dot H}^* \lb{4.64}
\end{equation}
and $H$ is not self-adjoint.

A typical example of a quasi-hermitian
extension is obtained as follows.

Let $\dot T$ denote the following first-order
differential
operator on the interval $[0, 2a]$, $a>0$,
\begin{align}
&(\dot T f)(x)=-i f'(x), \quad \xi \in (0, 2a), \no \\
&f\in\dom (\dot T)=\{g\in L^2([0,2a])\,|\,g\in AC([0, 2a]);
\, g(0_+)=g(2a_-)=0;  \lb{4.65} \\
&\hspace{8.2cm} g'\in L^2([0, 2a]) \}. \no
\end{align}
Then for $\rho \in \bbC\cup\{\infty\}$,
$|\rho | \ne 1$ the
operator $T_\rho$
\begin{align}
&( T_\rho f)(x)=-i f'(x), \quad \xi\in (0,2a), \no \\
&f\in\dom ( T_\rho)=\{g\in L^2([0,2a]) \, | \,
g\in AC([0, 2a]);\, g(0_+)=\rho g(2a_-); \lb{4.66} \\
&\hspace*{7.83cm} g'\in L^2([0, 2a]) \} \no
\end{align}
is a quasi-hermitian extension of $\dot T.$ (Here
$\rho=\infty$ in \eqref{4.66}, in obvious notation,
denotes
the boundary condition $g(2a_-)=0$.)
Among all quasi-hermitian extensions of $\dot T$ there
are two exceptional ones that have empty spectrum. In fact,
the operator
$T_0$ corresponding to the value $\rho=0$ in \eqref{4.66}
as well as  its adjoint, $T_0^*=T_\infty$, have empty
spectra,
that is,
\begin{equation}
\spec(T_0)=\spec(T_\infty)=\emptyset. \lb{4.67}
\end{equation}

The following theorem proven by Livsic in 1946  provides
an interesting characterization of this example.
\begin{theorem}\lb{t4.10} \mbox{\rm (Livsic \cite{Li60}.)}
For a densely defined closed prime symmetric operator with
deficiency indices $(1,1)$ to be unitarily equivalent to
the differentiation operator $\dot T$ in $L^2([0, 2a])$
for some $a>0$ it is necessary and sufficient that
it admits a quasi-hermitian extension with empty spectrum.
\end{theorem}

Using Livsic's result we are able to characterize the
model representation for the pair $(\dot H, H)$,
where $\dot H$ is a densely defined prime closed symmetric
operator with deficiency indices $(1,1)$ which admits
a quasi-hermitian extension
with empty spectrum, and $H$ a self-adjoint extension of
$\dot H$.

\begin{theorem}\lb{t4.11}
Let $\omega $ be a Borel measure on $\bbR$ such that
\begin{equation}
\int_\bbR\frac{d \omega(\lambda)}{1+\lambda^2}=1, \quad
\int_\bbR d\omega(\lambda) =\infty, \lb{4.67a}
\end{equation}
$H$ the self-adjoint operator of multiplication by
$\lambda$ in $L^2(\bbR;d\omega)$,
\begin{equation}
(Hf)(\lambda)=\lambda f(\lambda), \quad
f\in\dom (H)= L^2(\bbR;(1+\lambda^2)d\omega). \lb{4.68}
\end{equation}
Define $\dot H$ to be the densely defined closed prime
symmetric restriction of $H$,
\begin{equation}
\dot H=H\big\vert_{\dom (\dot H)}, \quad
\dom (\dot H)= \{ f\in \dom (H)\, \vert  \int_\bbR
d \omega(\lambda)f(\lambda)=0 \}, \lb{4.69}
\end{equation}
with deficiency indices $(1,1).$ Then $\dot H$ admits
a quasi-hermitian extension
with empty spectrum if and only if for some $a>0$ and
 some $\alpha\in [0, \pi)$
the following representation holds
\begin{align}
&\int_\bbR  d \omega (\lambda)((z-\lambda)^{-1}
-\lambda(1+\lambda^2)^{-1}) =\frac{
\sin(\alpha)-\cos(\alpha) (\cot (az)/\coth(a))}
{\cos(\alpha)+\sin(\alpha) (\cot (az)/\coth(a))},
\lb{4.70} \\
& \hspace*{9.8cm} z\in \bbC \backslash \bbR. \no
\end{align}
In this case the measure $\omega$ is a
pure point measure,
\begin{equation}
\omega=
\frac{\coth (a)(1+\cot^2(\alpha))}
{a(1+\cot^2(\alpha)\coth^2(a))}
\sum_{n\in \bbZ}
\mu_{\{ (\beta+\pi n)/a \} },
\lb{4.71}
\end{equation}
where $\mu_{\{x \}}$ denotes the  Dirac measure
supported at $\xi\in \bbR$ with mass one
and $\beta=\beta(\alpha, a) \in [0, \pi)$ is the
solution of the equation
\begin{equation}
\cot(\beta)+ \cot (\alpha) \coth (a)=1 \text{ if }
\alpha\in (0,\pi) \text{ and } \beta =0
\text{ if } \alpha =0.
\lb{4.72}
\end{equation}
Moreover, the self-adjoint operator
$H$ given by \eqref{4.68} is unitarily equivalent
to the differentiation operator $T_\rho$ in
\eqref{4.66} with
\begin{equation}
\rho=e^{2i\beta}.
\lb{4.73}
\end{equation}
\end{theorem}
\begin{proof}
That $\dot H$ is a densely defined closed prime
symmetric operator with deficiency indices $(1,1)$
is proven in \cite{Do65}. By Livsic's theorem,
Theorem~\ref{t4.10}, the pair
$(\dot H, H)$ is unitarily equivalent to the
pair $(\dot T, T_\rho)$, where $\dot T$ is the operator
\eqref{4.65} in $L^2([0, 2a])$ for some $a>0$
and $T_\rho$ is some
self-adjoint extension  of $\dot T$  given by
\eqref{4.66} for some $\rho$, $\vert \rho\vert =1$.
By \eqref{3.5} and \eqref{3.13} (cf.~also \eqref{4.24}) we
conclude
\begin{align}
m_{T_\rho}^D(z)&=\int_\bbR  d \omega (\lambda)
((z-\lambda)^{-1}-\lambda(1+\lambda^2)^{-1}) \lb{4.74} \\
&=z+(1+z^2) (u_+,(T_\rho-z)^{-1}u_+)_{L^2(\bbR;d\omega)},
\lb{4.75} \\
&\hspace*{2.2cm}u_+\in \ker ({\dot T}^*-i), \,
\|u_+\|_{L^2(\bbR;d\omega)}=1,
\no
\end{align}
where $m_{T_\rho}^D(z)$ denotes the Donoghue
Weyl $m$-function of the operator $T_\rho$.\\
Let $\wti T$ be the self-adjoint extension of
$\dot T$ corresponding to periodic boundary conditions,
\begin{equation}
\dom (\wti T)=\{g\in L^2([0,2a])\,|\,g\in AC([0, 2a]); \,
g(0_+)=g(2a_-); \, g'\in L^2([0, 2a]) \}. \lb{4.76}
\end{equation}
By Lemma \ref{l4.3} there exists an $\alpha
\in [0, \pi)$ such that
\begin{equation}
m_{T_\rho}^D(z)=\frac{\sin(\alpha)
+\cos(\alpha)m_{\ti T}^D(z)}
{\cos(\alpha)-\sin(\alpha) m_{\ti T}^D(z)}, \lb{4.77}
\end{equation}
where $ m_{\ti T}^D(z)$ is the Donoghue Weyl
$m$-function of
the extension $\wti T$
\begin{equation}
 m_{\ti T}^D(z)=z+(1+z^2)
(u_+,(\wti T-z)^{-1}u_+)_{L^2([0,2a];dx)},\quad z
\in \bbC\backslash \bbR.
\lb{4.78}
\end{equation}
The assertion \eqref{4.70} then follows from the fact
\begin{equation}
m_{\ti T}^D(z)=-\frac{\cot (az)}{\coth(a)}. \lb{4.79}
\end{equation}
Next we prove \eqref{4.79}.~First, we note that the
resolvent of the operator
$\wti T$ can be explicitly computed as
\begin{equation}
((\wti T -z)^{-1} f)(x)=
ie^{izx}\bigg(\int_0^xe^{-izt} f(t) dt
+
\frac{e^{2iza}}
{1-e^{2iza}}
\int_0^{2a}e^{-izt} f(t) dt\bigg), \quad
z\in\bbC\backslash\bbR. \lb{4.80}
\end{equation}
Next we calculate the quadratic form of the resolvent of
$\wti T$ on the element
$u_+(x)= 2^{1/2}(1-e^{-4a})^{-1/2}\exp(-x)$ generating
$\ker ({\dot T}^*-i)$. By \eqref{4.80} we have
\begin{equation}
((\wti T -z)^{-1}u_+)(x)=\frac{2^{1/2}
(1-e^{-4a})^{-1/2}}{i-z}
\bigg(e^{-x}- e^{izx}\frac{1-e^{-2a}}{1-e^{2iza}}\bigg)
\lb{4.81}
\end{equation}
and therefore,
\begin{equation}
(u_+,(\wti T -z)^{-1}u_+)_{L^2([0,2a];dx)}=
\frac{1}{i-z}\bigg(1+
\frac{2(1-e^{-2a})(1-e^{2iaz-2a})}
{(iz-1)(1-e^{-4a})(1-e^{2iaz})}\bigg). \lb{4.82}
\end{equation}
Equations \eqref{4.78} and \eqref{4.82} then prove
 \eqref{4.79}.\\
In order to prove \eqref{4.71} we note that the
right-hand
side of \eqref{4.70} is a periodic Herglotz
 function with period $\pi/a$. Such Herglotz
functions have simple poles at the points
$\{ (\beta +\pi n)/a \}_{n\in \bbZ}$ with
residues
\begin{align}
&\underset{z=(\beta +\pi n)/a }{\text{Res}}
\bigg( \frac{\sin(\alpha)-\cos(\alpha)
(\cot (az)/\coth(a))}
{\cos(\alpha)+\sin(\alpha)(\cot (az)/\coth(a))}
\bigg) \no \\
&=-\frac{\coth(a)(1+\cot^2(\alpha))}
{a(1+\cot^2(\alpha)\coth^2(a))} , \quad n\in \bbZ,
\lb{4.82a}
\end{align}
proving \eqref{4.71}.\\
The last assertion of the theorem follows from the fact
that the support of the measure $\omega$
coincides with the spectrum of $H$ and therefore with
the one of the operator $T_\rho$ which is unitarily
equivalent
to $H$. The spectrum of the
self-adjoint operator  $T_\rho$ can explicitly be
computed as
\begin{equation}
\spec (T_\rho)=\{\frac{1}{2a} \text{arg} \,
\rho +\frac{\pi}{a} n \}_{n\in \bbZ}. \lb{4.83}
\end{equation}
Since the sets \eqref{4.83} and $\supp (\omega)$
coincide we conclude \eqref{4.71}.
\end{proof}
\begin{remark}\lb{r4.12}
We note that the weak limit as $a \to \infty$ of
the measures $\omega=\omega (\alpha, a)$
(with $\alpha$ fixed) given by \eqref{4.71} coincides with
$\pi^{-1} d \lambda$, where $d\lambda$ denotes the Lebesgue
measure on $\bbR$.
\end{remark}

The next result shows that this limiting case
$d \omega=\pi^{-1} d \lambda$ is also rather exotic.

\begin{theorem}\lb{t4.13}
Let $\omega $ be a Borel measure on $\bbR$ such that
\begin{equation}
\int_\bbR\frac{d \omega(\lambda)}{1+\lambda^2}=1, \quad
 \int_\bbR d\omega(\lambda)=\infty, \lb{4.85}
\end{equation}
$H$ the self-adjoint operator of multiplication by
$\lambda$
in $L^2(\bbR;d\omega)$,
\begin{equation}
(Hf)(\lambda)=\lambda f(\lambda), \quad
f\in\dom (H)=L^2(\bbR;(1+\lambda^2)d\omega). \lb{4.86}
\end{equation}
Define $\dot H$ to be the densely defined closed prime
symmetric restriction of $H$,
\begin{equation}
\dot H=H\big\vert_{\dom (\dot H)}, \quad
\dom (\dot H)= \{f\in \dom (H) \, \vert \,
\int_\bbR f(\lambda) d \omega(\lambda)=0 \}, \lb{4.87}
\end{equation}
with deficiency indices $(1,1).$
Then $\dot H$ admits a quasi-hermitian extension
with pure point spectrum the open upper
(lower) half-plane and spectrum the closed upper
(lower) half-plane if and only if
the following representation holds,
\begin{equation}
\int_\bbR  d \omega (\lambda)((z-\lambda)^{-1}
-\lambda(1+\lambda^2)^{-1}) =
\begin{cases} i,& \Im (z)>0, \\
  -i,&\Im  (z)< 0. \end{cases} \lb{4.88}
\end{equation}
In this case
\begin{equation}
d\omega=\pi^{-1}d \lambda. \lb{4.89}
\end{equation}
\end{theorem}
\begin{proof}
The setup in \eqref{4.85}--\eqref{4.87} is identical to
that in Theorem~\ref{t4.11} and hence needs no further
comments.
The fact that  $\dot H$ is unitarily equivalent
to the differentiation operator $\dot T $ acting
in $L^2(\bbR;dx)$,
\begin{align}
&(\dot T f)(x)=-i f'(x), \quad
\xi \in \bbR, \no \\
&f\in\dom (\dot T)=\{g\in L^2(\bbR;dx) \,|\, g\in
AC(\bbR); \, g(0)=0; \, g'\in L^2(\bbR;dx) \} \lb{4.90}
\end{align}
goes back to Livsic (see, e.g., Appendix I.5 in
\cite{AG93}).~In fact, the quasi-hermitian extension $T$ of
$\dot T$ defined by
\begin{align}
&(Tf)(x)=-if'(x), \quad \xi\in\bbR\backslash\{0\},
\lb{4.91} \\
&f\in\dom (T)=\{g\in L^2(\bbR;dx) \,|\,
g\in AC([-R,0])\cup
AC([0,R]) \text{ for all }R>0; \no \\
&\hspace*{7.4cm} g(0_-)=0; \, g'\in L^2(\bbR;dx) \}. \no
\end{align}
(and its adjoint $T^*$ with corresponding boundary
condition
$g(0_+)=0$) has spectrum the closed upper (lower)
half-plane
with pure point spectrum the open upper (lower) half-plane,
respectively. This is easily verified from an alternative
expression for $T$ given by
\begin{equation}
T=\dot T_- \oplus T_+ \text{ in }
L^2(\bbR;dx)=L^2((-\infty,0];dx) \oplus L^2([0,\infty);dx),
\lb{4.91a}
\end{equation}
where
\begin{align}
&(\dot T_-f)(x)=-if'(x), \quad x<0, \lb{4.91b} \\
&f\in\dom (\dot T_-)=\{g\in L^2((-\infty,0];dx)
\,|\, g\in AC([-R,0]) \text{ for all }R>0; \no \\
&\hspace*{5.66cm} g(0_-)=0; \, g'\in L^2((-\infty,0];dx) \},
\no \\
&(T_+f)(x)=-if'(x), \quad x>0, \lb{4.91c} \\
&f\in\dom (T_+)=\{g\in L^2([0,\infty);dx)
\,|\, g\in AC([0,R]) \text{ for all }R>0; \no \\
&\hspace*{7.1cm} g'\in L^2([0,\infty);dx) \}. \no
\end{align}
The explicit expressions for the resolvents of
$\dot T_-$ and $T_+$ (see, e.g., \cite{Ka80}, Example
III.6.9) then show that both operators have spectrum the
closed upper half-plane, that is,
\begin{equation}
\spec (\dot T_-)=\spec (T_+)=\ol{\bbC_+}. \lb{4.91d}
\end{equation}
Together with the aforementioned result of Livsic, this
shows that the  pair $(\dot H, H)$ is
unitarily equivalent to the
pair $(\dot T, T_\rho) $, where $T_\rho,$ $|\rho|=1$
is some
self-adjoint extension  of $\dot T$ in $L^2(\bbR;dx)$,
\begin{align}
&(T_\rho f)(x)=-if'(x),
\quad \xi\in\bbR\backslash\{0\}, \,
|\rho|=1, \lb{4.92} \\
&f\in\dom (T_\rho)=\{g\in L^2(\bbR;dx) \, | \,
g\in AC([-R,0])\cup AC([0,R]) \text{ for all }R>0; \no \\
&\hspace*{6.7cm} g(0_-)=\rho g(0_+); \,
g'\in L^2(\bbR;dx) \}. \no
\end{align}
Since the pair
$(\dot T, T_1)$ is unitarily equivalent to the model pair
$(\dot H, H)$ in \eqref{4.86} and \eqref{4.87} (it
suffices applying
the Fourier transform), where $d\omega=\pi^{-1} d \lambda$,
 we can immediately compute the Donoghue Weyl $m$-function
$m^D_{T_1}(z)$ of the self-adjoint extension $T_1$,
\begin{equation} \lb{4.93}
m^D_{T_1}(z)=\frac{1}{\pi}
\int_\bbR  d\,\lambda((z-\lambda)^{-1}
-\lambda(1+\lambda^2)^{-1}) =
\begin{cases} i,& \Im (z)>0, \\
  -i,&\Im (z)< 0. \end{cases}
\end{equation}
Since
\begin{equation}
\pm i=\frac{\sin(\alpha)+\cos(\alpha)(\pm i)}
{\cos(\alpha)-\sin(\alpha)(\pm i)} \text{ for all }
\alpha\in [0,\pi), \lb{4.94}
\end{equation}
Lemma \ref{l4.3} implies that the Donoghue Weyl
$m$-function
 $m^D_{T_\rho}(z)$ of the extension $T_\rho$ is
independent of $\rho,$ $\vert \rho \vert=1$ and hence
$m^D_{T_1}(z)=m^D_\rho(z)$.
Therefore, the model representation
for the pair $(\dot T, T_\rho)$ is given by
\eqref{4.85}--\eqref{4.87} with $d\omega=\pi^{-1}
d \lambda$, proving \eqref{4.89}. Finally,
\eqref{4.88} follows from \eqref{4.93}.
\end{proof}

\medskip

{\bf Acknowledgments.}~Eduard Tsekanovskii would like to
thank Daniel Alpay and Victor Vinnikov for support
and their kind invitation
to a wonderful conference organized in honor of the
80th birthday of Moshe Livsic, his dear and inspiring
teacher.



\begin{thebibliography}{99}
%
\bibitem{AS72} M.~Abramowitz and I.~A.~Stegun, {\it
Handbook of Mathematical
Functions}, Dover, New York, 1972.
%
\bi{AL95} V.~M.~Adamjan and H.~Langer, {\it Spectral
Properties of a class of rational operator valued
functions}, J. Operator Th. {\bf 33}, 259--277 (1995).
%
\bibitem{AG93} N.~I.~Akhiezer and I.~M.~Glazman, {\it
Theory of Linear
Operators in Hilbert Space}, Dover, New York, 1993.
%
\bi{ABN96} S.~Albeverio, J.~Brasche, and H.~Neidhardt,
{\it On inverse spectral theory for self-adjoint
extensions: mixed types of spectra}, preprint, 1996.
%
\bibitem{AGHKH88} S.~Albeverio, F.~Gesztesy,
R.~H{\o}egh-Krohn, H.~Holden,
{\it Solvable Models in Quantum Mechanics}, Springer,
Berlin, 1988.
%
\bi{AS80} A.~Alonso and B.~Simon, {\it The
Birman-Krein-Vishik theory of self-adjoint extensions of
semibounded operators}, J. Operator Th. {\bf 4}, 251--270
(1980).
%
\bi{AN70} T.~Ando and K.~Nishio, {\it Positive
selfadjoint extensions of positive symmetric operators},
T\^ohoku Math. J. {\bf 22}, 65--75 (1970).
%
\bi{AT82} Yu.~M.~Arlinskii and E.~R.~Tsekanovskii,
{\it Nonselfadjoint
contractive extensions of a Hermitian contraction and
theorems of Krein}, Russ. Math. Surv. {\bf 37:1},
151--152 (1982).
%
\bibitem{Ar57} N.~Aronszajn, {\it On a problem of Weyl
in the theory of
singular Sturm-Liouville equations}, Am. J. Math.
{\bf 79}, 597-610
(1957).
%
\bi{BGK79} H.~Bart, I.~Gohberg, and M.~A.~Kaashoek, {\it
Minimal Factorizations of Matrix and Operator Functions},
Operator Theory: Advances and Applications, Vol.~1,
Birkh\"auser, Basel, 1979.
%
\bi{BW83} H.~Baumg\"artel and M.~Wollenberg, {\it
Mathematical Scattering Theory}, Birkh\"auser, Basel,
1983.
%
\bi{BS87} M.~S.~Birman and M.~Z.~Solomjak, {\it Spectral
Theory of Self-Adjoint Operators in Hilbert Space},
Reidel, Dordrecht, 1987.
%
\bi{BN94} J.~F.~Brasche and H.~Neidhardt, {\it Some
remarks on Krein's extension theory}, Math. Nachr.
{\bf 165}, 159--181 (1994).
%
\bi{BN95} J.~Brasche and H.~Neidhardt, {\it On the
absolutely continuous spectrum of self-adjoint
extensions}, J. Funct. Anal. {\bf 131}, 364--385 (1995).
%
\bi{BN96} J.~Brasche and H.~Neidhardt, {\it On the
singularly continuous spectrum of self-adjoint
extensions}, Math. Z. {\bf 222}, 533--542 (1996).
%
\bi{BNW93} J.~Brasche, H.~Neidhardt, and J.~Weidmann,
{\it On the point spectrum of self-adjoint extensions},
Math. Z. {\bf 214}, 343--355 (1993).
%
\bi{Br71} M.~S.~Brodskii, {\it Triangular and Jordan
Representations of Linear Operators}, Amer. Math. Soc.,
Providence, RI, 1971.
%
\bi{Bu97} D.~Buschmann, {\it Spektraltheorie
verallgemeinerter Differentialausdr\"ucke - Ein
neuer Zugang}, Ph.D. Thesis, University
of Frankfurt, Germany, 1997.
%
\bi{Ca76} R.~W.~Carey, {\it A unitary invariant for
pairs of
self-adjoint
operators}, J. Reine Angew. Math. {\bf 283}, 294--312
(1976).
%
\bibitem{CL85}
E.~A.~Coddington and N.~Levinson,{\it Theory of Ordinary
Differential
Equations}, Krieger, Malabar, 1985.
%
\bi{DM87} V.~A.~Derkach and M.~M.~Malamud, {\it On the
Weyl function and Hermitian operators with gaps}, Sov.
Math. Dokl. {\bf 35}, 393--398 (1987).
%
\bi{DM91} V.~A.~Derkach and M.~M.~Malamud, {\it
Generalized
resolvents
and the boundary value problems for Hermitian operators
with
gaps}, J.
Funct. Anal. {\bf 95}, 1--95 (1991).
%
\bi{DM95} V.~A.~Derkach and M.~M.~Malamud, {\it The
extension
theory of
Hermitean operators and the moment problem}, J.~Math.~Sci.
{\bf 73},
141--242 (1995).
%
\bibitem{DMT88} V.~A.~Derkach, M.~M.~Malamud, and
E.~R.~Tsekanovskii, {\it
Sectorial extensions of a positive operator, and the
characteristic
function}, Sov. Math. Dokl.  {\bf 37}, 106-110 (1988).
%
\bibitem{Do65} W.~F.~Donoghue, {\it On the perturbation
of spectra},
Commun. Pure Appl. Math. {\bf 18}, 559-579 (1965).
%
\bi{DS88} N.~Dunford and J.~T.~Schwartz, {\it Linear
Operators Part II:
Spectral Theory}, Interscience, New York, 1988.
%
\bi{EK82} M.~S.~P.~Eastham and H.~Kalf, {\it
Schr\"{o}dinger--Type
Operators with Continuous Spectra}, Pitman, Boston, 1982.
%
\bi{GMT97} F.~Gesztesy, K.A. Makarov, E.~Tsekanovskii,
{\it An Addendum to Krein's formula}, preprint, 1997.
%
\bibitem{GS96} F.~Gesztesy and B.~Simon, {\it Uniqueness
theorems in
inverse spectral theory for one-dimensinal Schr\"odinger
operators},
Trans. Amer. Math. Soc. {\bf 348}, 349-373 (1996).
%
\bi{GT97} F.~Gesztesy and E.~Tsekanovskii, {\it On
matrix-valued
Herglotz functions}, preprint, 1997.
%
\bi{Gr60} D.~S.~Greenstein, {\it On the analytic
continuation of
functions which map the upper half plane into itself},
J. Math.
Anal. Appl. {\bf 1}, 355--362 (1960).
%
\bi{HS97a} S.~Hassi and H.~de Snoo, {\it One-dimensional
graph
perturbations of selfadjoint relations}, Ann. Acad. Sci.
Fenn. A I Math. {\bf 22}, 123--164 (1997).
%
\bi{HS97b} S.~Hassi and H.~de Snoo, {\it On rank one
perturbations of selfadjoint operators}, Integral Eqs.
Operator Th. {\bf 29}, 288--300 (1997).
%
\bibitem{HKS97a} S.~Hassi, M.~Kaltenb\"ack, and H.~de
Snoo,
{\it Triplets of Hilbert spaces and Friedrichs extensions
associated with the subclass $N_1$ of Nevanlinna
functions}, J. Operator Th. {\bf 37}, 155--181 (1997).
%
\bi{HKS97b} S.~Hassi, M.~Kaltenb\"ack, and H.~de Snoo,
{\it Generalized Krein-von Neumann extensions and
associated operator models}, preprint, 1997.
%
\bi{HKS98} S.~Hassi, M.~Kaltenb\"ack, and H.~de Snoo,
{\it Generalized finite rank perturbations associated
with Kac classes of matrix Nevanlinna functions}, in
preparation.
%
\bibitem{HSW98} S.~Hassi, H.~S.~V.~de Snoo, and
A.~D.~I.~Willemsma,
{\it Smooth rank one perturbations of selfadjoint
operators}, Proc. Amer. Math. Soc., to appear.
%
\bi{Ka80} T.~Kato, {\it Perturbation Theory for Linear
Operators}, corr. 2nd ed., Springer, Berlin, 1980.
%
\bibitem{Kr47} M.~G.~Krein, {\it The theory of
self-adjoint extensions
of semibounded Hermitian transformations and its
applications}, I. Mat.
Sb. {\bf 20}, 431-495 (1947) (Russian).
%
\bi{KO77} M.~G.~Krein and I.~E.~Ovcharenko, {\it
$Q$-functions and sc-resolvents of nondensely defined
hermitian contractions}, Sib. Math. J. {\bf 18},
728--746 (1977).
%
\bi{KO78}  M.~G.~Krein and I.~E.~Ov\v carenko,
{\it Inverse problems
for $Q$-functions and resolvent matrices of positive
hermitian operators},
Sov. Math. Dokl. {\bf 19}, 1131--1134 (1978).
%
\bibitem{KS74} M.~G.~Krein and Ju.~L.~Smul'jan,
{\it On linear-fractional
transformations with operator coefficients}, Amer.
Math. Soc. Transl.~(2)
{\bf 103}, 125-152 (1974).
%
\bibitem{LT77} H.~Langer, B.~Textorius,
{\it On generalized resolvents and $Q$-functions
of symmetric
linear relations (subspaces) in Hilbert space},
Pacific J. Math. {\bf 72}, 135--165 (1977).
%
\bi{Li60} M.~S.~Livsic, {\it On a class of
linear operators in Hilbert space}, Amer. Math. Soc.
Transl.~(2) {\bf 13}, 61--83 (1960).
%
\bi{Ma92a} M.~M.~Malamud, {\it Certain classes of
extensions of a lacunary Hermitian operator},
Ukrain. Math. J. {\bf 44}, 190--204 (1992).
%
\bi{Ma92b} M.~M.Malamud, {\it On a formula of the
generalized resolvents
of a nondensely defined hermitian operator},
Ukrain. Math. J. {\bf 44}, 1522--1547 (1992).
%
\bi{MS96} R.~Mennicken and A.~A.~Shkalikov, {\it
Spectral decomposition of symmetric operator matrices},
Math. Nachr. {\bf 179}, 259--273 (1996).
%
\bi{Na89} S.~N.~Naboko, {\it Boundary values of
analytic operator functions with a positive imaginary
part}, J. Sov. Math. {\bf 44}, 786--795 (1989).
%
\bi{Na90} S.~N.~Naboko, {\it Nontangential boundary
values of operator-valued $R$-functions in a
half-plane}, Leningrad Math. J. {\bf 1}, 1255--1278
(1990).
%
\bi{Na91} S.~N.~Naboko, {\it Structure of the
singularities of operator functions with a positive
imaginary part}, Funct. Anal. Appl. {\bf 25},
243--253 (1991).
%
\bi{Na96} S.~N.~Naboko, {\it Zygmund's theorem and
the boundary behavior of operator $R$-functions},
Funct. Anal. Appl. {\bf 30}, 211--213 (1996).
%
\bi{Na40} M.~A.~Naimark, {\it Spectral functions of a
symmetric operator}, Izv. Akad. Nauk SSSR {\bf 4},
227--318 (1940). (Russian.)
%
\bi{Na43} M.~A.~Naimark, {\it On a representation of
additive operator set functions}, Dokl. Akad. Nauk SSSR
{\bf 41}, 359--361 (1943). (Russian.)
%
\bibitem{Na68} M.~A.~Naimark, {\it Linear Differential
Operators II},
F. Ungar Publ., New York, 1968.
%
\bi{RS75} M.~Reed and B.~Simon, {\it Methods of Modern
Mathematical Physics, Vol.~II: Fourier Analysis,
Self-Adjointness}, Academic Press, New York, 1975.
%
\bi{Re98} C.~Remling, {\it Spectral analysis of higher
order differential operators I: General properties of
the $M$-function}, J. London Math. Soc., to appear.
%
\bi{RR97} M. Rosenblum, J. Rovnyak,
{\it Hardy Classes and Operator Theory}, Dover, New
York, 1997.
%
\bibitem{Sa65} Sh.~N.~Saakjan,
{\it On the theory of the resolvents of a symmetric
operator with
infinite deficiency indices},
Dokl. Akad. Nauk Arm. SSR {\bf 44},  193--198 (1965).
(Russian.)
%
\bibitem{Si95} B.~Simon, {\it Spectral analysis of
rank one
perturbations and applications}, CRM Proceedings and
Lecture
Notes {\bf
8}, 109-149 (1995).
%
\bi{SW86} B.~Simon and T.~Wolff, {\it Singular continuous
spectrum
under rank one perturbations and localization for random
Hamiltonians},
Commun. Pure Appl. Math. {\bf 39}, 75-90 (1986).
%
\bibitem{Sk79} C.~F.~Skau, {\it Positive self-adjoint
extensions of
operators affiliated with a von Neumann algebra}, Math.
Scand. {\bf 44},
171-195 (1979).
%
\bi{Sh71} Yu.~L.~Shmul'yan, {\it On operator R-functions},
Sib. Math. J. {\bf 12}, 315--322 (1971).
%
\bi{SF70} B.~Sz.-Nagy and C.~Foias, {\it Harmonic Analysis
of Operators on Hilbert Space}, North-Holland, Amsterdam,
1970.
%
\bi{Ti62} E.~C.~Titchmarsh, {\it Eigenfunction Expansions
Associated with
Second-Order Differential Equations}, Part I, 2nd ed.,
Oxford
University
Press, Oxford, 1962.
%
\bi{Ts80} E.~R.~Tsekanovskii, {\it Non-self-adjoint
accretive extensions of
positive operators and theorems of
Friedrichs-Krein-Phillips}, Funct. Anal.
Appl. {\bf 14}, 156--157 (1980).
%
\bi{Ts81} E.~R.~Tsekanovskii, {\it Friedrichs and
Krein extensions of
positive operators and holomorphic contraction semigroups},
Funct. Anal.
Appl. {\bf 15}, 308--309 (1981).
%
\bibitem{Ts92} E.~R.~Tsekanovskii, {\it Accretive
extensions and problems
on the Stieltjes operator-valued functions realizations},
in {\it Operator Theory and Complex Analysis},
T.~Ando and I.~Gohberg (eds.), Operator
Theory: Advances and Applications, Vol.~59, Birkh\"auser,
Basel, 1992, pp.~328--347.
%
\bi{TS77} E.~R.~Tsekanovskii and Yu.~L.~Shmul'yan, {\it The
theory of
bi-extensions of operators on rigged Hilbert spaces.
Unbounded
operator
colligations and characteristic functions}, Russ. Math.
Surv.
{\bf 32:5}, 73--131 (1977).
%
\bibitem{Ne29} J.~von Neumann, {\it Allgemeine
Eigenwerttheorie
hermitescher Funktionaloperatoren}, Math. Ann.
{\bf 102}, 49--131
(1929--30).
%
\bi{We10} H.~Weyl, {\it \"Uber gew\"ohnliche
Differentialgleichungen
mit Singularit\"aten und die zugeh\"origen Entwicklungen
willk\"urlicher
Funktionen}, Math. Ann. {\bf 68}, 220--269 (1910).
%

\end{thebibliography}
\end{document}